\input amstexl


\catcode`\@=11
\ifx\amstexloaded@\relax\else
 \errmessage{AmS-TeX must be loaded before LamS-TeX}\fi
\ifx\laxread@\undefined\else\catcode`\@=\active \fi
\def\err@#1{\errmessage{LamS-TeX error: #1}}
\def^^L{\par}
\let\+\tabalign
\def\newcount{\alloc@0\count\countdef\insc@unt}
\def\newdimen{\alloc@1\dimen\dimendef\insc@unt}
\def\newskip{\alloc@2\skip\skipdef\insc@unt}
\def\newmuskip{\alloc@3\muskip\muskipdef\@cclvi}
\def\newbox{\alloc@4\box\chardef\insc@unt}
\let\newtoks\relax
\def\newhelp#1#2{\newtoks#1#1\expandafter{\csname#2\endcsname}}
\def\newtoks{\alloc@5\toks\toksdef\@cclvi}
\def\newread{\alloc@6\read\chardef\sixt@@n}
\def\newwrite{\alloc@7\write\chardef\sixt@@n}
\def\newfam{\alloc@8\fam\chardef\sixt@@n}
\def\newlanguage{\alloc@9\language\chardef\@cclvi}
\def\newinsert#1{\global\advance\insc@unt by\m@ne
  \ch@ck0\insc@unt\count
  \ch@ck1\insc@unt\dimen
  \ch@ck2\insc@unt\skip
  \ch@ck4\insc@unt\box
  \allocationnumber=\insc@unt
  \global\chardef#1=\allocationnumber
  \wlog{\string#1=\string\insert\the\allocationnumber}}
\def\newif#1{\count@\escapechar \escapechar\m@ne
  \expandafter\expandafter\expandafter
   \edef\@if#1{true}{\let\noexpand#1=\noexpand\iftrue}%
  \expandafter\expandafter\expandafter
   \edef\@if#1{false}{\let\noexpand#1=\noexpand\iffalse}%
  \@if#1{false}\escapechar\count@}

\def\Err@#1{\errhelp\defaulthelp@\err@{#1}}
{\catcode`\@=\active
 \edef\next{\gdef\noexpand@{\futurelet\noexpand\next
  \csname at\string@\endcsname}}
 \next
}
\def\at@{\ifcat\noexpand\next a\let\next@\at@@\else
 \ifcat\noexpand\next0\let\next@\at@@\else
 \ifcat\noexpand\next\relax\let\next@\at@@\else
 \let\next@\at@@@\fi\fi\fi\next@}
\def\at@@@{\errhelp\athelp@\err@{Invalid use of @}}
\def\at@@#1{\expandafter
 \ifx\csname\string#1@at\endcsname\relax\let\next@\at@@@\else
 \DN@{\csname\string#1@at\endcsname}\fi\next@}
\def\atdef@#1{\expandafter\def\csname\string#1@at\endcsname}
\newif\iftest@
\def\tagin@#1{\tagin@false
 \DN@##1\tag##2##3\next@{\test@true\ifx\tagin@##2\test@false\fi}%
 \next@#1\tag\tagin@\next@\tagin@false\iftest@\tagin@true\fi}
\let\lkerns@\relax
\def\nolinebreak{\RIfM@\mathmodeerr@\nolinebreak\else
 \ifhmode\saveskip@\lastskip\unskip
 \nobreak\ifdim\saveskip@>\z@\hskip\saveskip@\fi\lkerns@
 \else\vmodeerr@\nolinebreak\fi\fi}
\def\allowlinebreak{\RIfM@\mathmodeerr@\allowlinebreak\else
 \ifhmode\saveskip@\lastskip\unskip
 \allowbreak\ifdim\saveskip@>\z@\hskip\saveskip@\fi\lkerns@
 \else\vmodeerr@\allowlinebreak\fi\fi}
\def\linebreak{\RIfM@\mathmodeerr@\linebreak\else
 \ifhmode\unskip\unkern\break\lkerns@
 \else\vmodeerr@\linebreak\fi\fi}
\let\nkerns@\relax
\def\newline{\RIfM@\mathmodeerr@\newline\else
 \ifhmode\unskip\unkern\null\hfill\break\nkerns@
 \else\vmodeerr@\newline\fi\fi}%
\def\newbox@{\alloc@@4\box\chardef\insc@unt}
\def\newcount@{\alloc@@0\count\countdef\insc@unt}
\def\accentedsymbol#1#2{\expandafter\newbox@\csname\exstring@#1@box\endcsname
 \setbox\csname\exstring@#1@box\endcsname\hbox{$\m@th#2$}%
 \define#1{\copy\csname\exstring@#1@box\endcsname{}}}
\def\rightadd@#1\to#2{\toks@{\\#1}\toks@@\expandafter{#2}\xdef#2{\the\toks@@
 \the\toks@}\toks@{}\toks@@{}}
\def\fontlist@{\\\tenrm\\\sevenrm\\\fiverm\\\teni\\\seveni\\\fivei
 \\\tensy\\\sevensy\\\fivesy\\\tenex\\\tenbf\\\sevenbf\\\fivebf
 \\\tensl\\\tenit}
\def\font@#1=#2 {\rightadd@#1\to\fontlist@\font#1=#2 }
\def\ismember@#1#2{\global\let\Next@ F\let\next@= #2%
 {\def\\##1{\let\nextii@##1\ifx\nextii@\next@\global\let\Next@ T\fi}#1}%
 \test@false\ifx\Next@ T\test@true\fi\let\next@\relax}
\def\FNSS@#1{\let\FNSS@@#1\FN@\FNSS@@@}
\def\FNSS@@@{\ifx\next\space@\def\FNSS@@@@. {\FN@\FNSS@@@}\else
 \def\FNSS@@@@.{\FNSS@@}\fi\FNSS@@@@.}
\atdef@"{\unskip
 \DN@{\ifx\next`\DN@`{\FN@\nextii@}%
  \else\ifx\next\lq\DN@\lq{\FN@\nextii@}%
  \else\DN@####1{\FN@\nextiii@}\fi\fi
  \next@}%
 \DNii@{\ifx\next`\DN@`{\sldl@``}%
  \else\ifx\next\lq\DN@\lq{\sldl@``}%
  \else\DN@{\dlsl@`}\fi\fi\next@}%
 \def\nextiii@{\ifx\next'\DN@'{\srdr@''}%
  \else\ifx\next\rq\DN@\rq{\srdr@''}%
  \else\DN@{\drsr@'}\fi\fi\next@}%
 \FNSS@\next@}
\def\root{%
  \DN@{\ifx\next\uproot\let\next@\nextii@\else
   \ifx\next\leftroot\let\next@\nextiii@\else
   \let\next@\plainroot@\fi\fi\next@}%
  \DNii@\uproot##1{\uproot@##1\relax\FNSS@\nextiv@}%
  \def\nextiv@{\ifx\next\leftroot\let\next@\nextv@\else
   \let\next@\plainroot@\fi\next@}%
  \def\nextv@\leftroot##1{\leftroot@##1\relax\plainroot@}%
  \def\nextiii@\leftroot##1{\leftroot@##1\relax\FNSS@\nextvi@}%
  \def\nextvi@{\ifx\next\uproot\let\next@\nextvii@\else
   \let\next@\plainroot@\fi\next@}%
  \def\nextvii@\uproot##1{\uproot@##1\relax\plainroot@}%
  \bgroup\uproot@\z@\leftroot@\z@
 \FNSS@\next@}
\def\loop#1\repeat{\def\iterate{#1\relax\expandafter\iterate\fi}%
 \iterate\let\iterate\relax}
\def\gloop@#1\repeat{\gdef\iterate@{#1\relax\expandafter\iterate@\fi}%
 \iterate@\global\let\iterate@\relax}
\def\printoptions{\W@{Do you want S(yntax check),
  G(alleys) or P(ages)?^^JType S, G or P, follow by <return>: }\loop
 \read\m@ne to\ans@
 \edef\next@{\def\noexpand\Ans@{\ans@}}\uppercase\expandafter{\next@}%
 \ifx\Ans@\S@\test@true\syntax\else
 \ifx\Ans@\G@\test@true\galleys\else
 \ifx\Ans@\P@\test@true\else
 \test@false\fi\fi\fi
 \iftest@\else\W@{Type S, G or P, follow by <return>: }%
 \repeat}
\expandafter\let\csname A@;\endcsname;
\expandafter\let\csname A@:\endcsname:
\expandafter\let\csname A@?\endcsname?
\expandafter\let\csname A@!\endcsname!
\def\APdef#1{\def\next@{\expandafter\let\csname A@\string#1\endcsname#1}%
 \afterassignment\next@\def#1}
\let\fextra@\,
\def\tdots@{\unskip
 \DN@{$\m@th\mathinner{\ldotp\ldotp\ldotp}\,
   \ifx\next,\,$\else\ifx\next.\,$\else
   \ifx\next;\,$\else
   \expandafter\ifx\csname A@\string;\endcsname\next\fextra@$\else
   \ifx\next:\,$\else
   \expandafter\ifx\csname A@\string:\endcsname\next\fextra@$\else
   \ifx\next?\,$\else
   \expandafter\ifx\csname A@\string?\endcsname\next\fextra@$\else
   \ifx\next!\,$\else
   \expandafter\ifx\csname A@\string!\endcsname\next\fextra@$\else
   $ \fi\fi\fi\fi\fi\fi\fi\fi\fi\fi}%
 \ \FN@\next@}
\def\extrap@#1{%
 \ifx\next,\DN@{#1\,}\else
 \ifx\next;\DN@{#1\,}\else
 \expandafter\ifx\csname A@\string;\endcsname\next\DN@{#1\fextra@}\else
 \ifx\next.\DN@{#1\,}\else\extra@
 \ifextra@\DN@{#1\,}\else
 \let\next@#1\fi\fi\fi\fi\fi\next@}
\def\dotsc{\DN@{\ifx\next;\plainldots@\,\else
 \expandafter\ifx\csname A@\string;\endcsname\next\plainldots@\fextra@\else
 \ifx\next.\plainldots@\,\else\extra@\plainldots@
 \ifextra@\,\fi\fi\fi\fi}%
 \FN@\next@}
\def\keybin@{\keybin@true
 \ifx\next+\else\ifx\next=\else\ifx\next<\else\ifx\next>\else\ifx\next-\else
 \ifx\next*\else\ifx\next:\else
 \expandafter\ifx\csname A@\string;\endcsname\next\else
 \keybin@false\fi\fi\fi\fi\fi\fi\fi\fi}
\def\boldkey#1{\ifcat\noexpand#1A%
  \ifcmmibloaded@{\fam\cmmibfam#1}\else
   \Err@{First bold symbol font not loaded}\fi
 \else
 \let\next=#1%
 \ifx#1!\mathchar"5\bffam@21 \else
 \expandafter\ifx\csname A@\string!\endcsname\next\mathchar"5\bffam@21 \else
 \ifx#1(\mathchar"4\bffam@28 \else\ifx#1)\mathchar"5\bffam@29 \else
 \ifx#1+\mathchar"2\bffam@2B \else\ifx#1:\mathchar"3\bffam@3A \else
 \expandafter\ifx\csname A@\string:\endcsname\next\mathchar"3\bffam@3A \else
 \ifx#1;\mathchar"6\bffam@3B \else
 \expandafter\ifx\csname A@\string;\endcsname\next\mathchar"6\bffam@3B \else
 \ifx#1=\mathchar"3\bffam@3D \else
 \ifx#1?\mathchar"5\bffam@3F \else
 \expandafter\ifx\csname A@\string?\endcsname\next\mathchar"5\bffam@3F \else
 \ifx#1[\mathchar"4\bffam@5B \else
 \ifx#1]\mathchar"5\bffam@5D \else
 \ifx#1,\mathchari@63B \else
 \ifx#1-\mathcharii@200 \else
 \ifx#1.\mathchari@03A \else
 \ifx#1/\mathchari@03D \else
 \ifx#1<\mathchari@33C \else
 \ifx#1>\mathchari@33E \else
 \ifx#1*\mathcharii@203 \else
 \ifx#1|\mathcharii@06A \else
 \ifx#10\bold0\else\ifx#11\bold1\else\ifx#12\bold2\else\ifx#13\bold3\else
 \ifx#14\bold4\else\ifx#15\bold5\else\ifx#16\bold6\else\ifx#17\bold7\else
 \ifx#18\bold8\else\ifx#19\bold9\else
  \Err@{\noexpand\boldkey can't be used with #1}%
 \fi\fi\fi\fi\fi\fi\fi\fi\fi\fi\fi\fi\fi\fi\fi
 \fi\fi\fi\fi\fi\fi\fi\fi\fi\fi\fi\fi\fi\fi\fi\fi\fi\fi}
\def\arabic#1{#1}
\def\alph#1{\count@#1\relax\advance\count@96 \ifnum\count@>122
 \Err@{\noexpand\alph invalid for numbers > 26}\else\char\count@\fi}
\def\Alph#1{\count@#1\relax\advance\count@64 \ifnum\count@>90
 \Err@{\noexpand\Alph invalid for numbers > 26}\else\char\count@\fi}

\def\Roman#1{\uppercase\expandafter{\romannumeral#1}}
\def\fnsymbol#1{\count@#1\relax
 \count@@\count@
 \advance\count@\m@ne\divide\count@7
 \count@@@\count@\advance\count@@@\@ne
 \multiply\count@7 \advance\count@@-\count@
 \count@\count@@@
 {\loop
  \ifcase\count@@\or*\or\dag\or\ddag\or\P\or\S\or\text{$\|$}\or\#\fi
  \advance\count@\m@ne\ifnum\count@>\z@\repeat}}
\def\cardnine@#1{\ifcase#1\or one\or two\or three\or four\or five\or
 six\or seven\or eight\or nine\fi}
\let\alloc@\alloc@@
\newcount\ten@
\ten@10
\def\cardinal#1{\count@#1\relax
 \ifnum\count@>99 \number\count@
 \else
  \ifnum\count@=\z@ zero%
  \else
   \ifnum\count@<\ten@\cardnine@\count@
   \else
    \ifnum\count@<20
     \advance\count@-\ten@
     \ifcase\count@ ten\or eleven\or twelve\or thirteen\or fourteen\or
      fifteen\or sixteen\or seventeen\or eighteen\or nineteen\fi
    \else
     \count@@\count@\count@@@\count@@
     \divide\count@\ten@\multiply\count@\ten@
     \advance\count@@@-\count@\divide\count@\ten@
     \ifcase\count@\or\or twenty\or thirty\or forty\or fifty\or sixty\or
      seventy\or eighty\or ninety\fi
     \ifnum\count@@@=\z@\else-\cardnine@\count@@@\fi
    \fi
   \fi
  \fi
 \fi}
\def\ordnine@#1{\ifcase#1\or first\or second\or third\or fourth\or fifth\or
 sixth\or seventh\or eighth\or ninth\fi}
\newcount\count@@@@
\def\ordsuffix@{\count@@@@\count@
 \divide\count@\ten@
 \count@@@\count@\count@@\count@
 \divide\count@@\ten@\multiply\count@@\ten@
 \advance\count@@@-\count@@
 \ifnum\count@@@=\@ne th%
 \else
  \count@@@\count@@@@
  \count@@\count@@@@
  \divide\count@@\ten@\multiply\count@@\ten@
  \advance\count@@@-\count@@
  \ifcase\count@@@ th\or st\or nd\or rd\else th\fi
 \fi}
\def\nordinal#1{\count@#1\relax\number\count@\ordsuffix@}
\def\spordinal#1{\count@#1\relax\number\count@$^{\text{\ordsuffix@}}$}
\def\ordinal#1{\count@#1\relax
 \ifnum\count@>99 \number\count@\ordsuffix@
 \else
   \ifnum\count@=\z@ zeroth%
  \else
    \ifnum\count@<\ten@\ordnine@\count@
    \else
     \ifnum\count@<20 \advance\count@-\ten@
      \ifcase\count@ tenth\or eleventh\or twelfth\or thirteenth\or
       fourteenth\or fifteenth\or sixteenth\or seventeenth\or eighteenth\or
       nineteenth\fi
     \else
      \count@@\count@
      \divide\count@\ten@\multiply\count@\ten@
      \count@@@\count@@\advance\count@@@-\count@
      \divide\count@\ten@
      \ifcase\count@\or\or twent\or thirt\or fort\or fift\or sixt\or sevent\or
       eight\or ninet\fi
      \ifnum\count@@@=\z@ ieth\else y-\ordnine@\count@@@\fi
     \fi
    \fi
  \fi
 \fi}
\font@\tensmc=cmcsc10
\textonlyfont@\smc\tensmc
\newtoks\noexpandtoks@
\noexpandtoks@{\let\arabic\relax\let\alph\relax\let\Alph\relax
 \let\Roman\relax\let\fnsymbol\relax\let\rm\relax
 \let\it\relax\let\bf\relax\let\sl\relax\let\smc\relax
 \let\/\relax\let\null\relax}
\def\noexpands@{\the\noexpandtoks@}
\def\Nonexpanding#1{\global\noexpandtoks@
 \expandafter{\the\noexpandtoks@\let#1\relax}}
\def\prevanish@{\saveskip@\z@\ifhmode\saveskip@\lastskip\unskip\fi}
\def\postvanish@{\ifdim\saveskip@>\z@\hskip\saveskip@\fi\FN@\postvanish@@}
\def\postvanish@@{\DN@.{}%
 \ifx\next\space@\ifdim\saveskip@>\z@\DN@. {}\fi\fi\next@.}
\def\invisible#1{\prevanish@\ignorespaces#1\unskip\postvanish@}
\def\vanishlist@{\\\invisible}
\let\noindent@\noindent
\def\noindent{\par\noindent@\FN@\pretendspace@}
\def\pretendspace@{\ismember@\vanishlist@\next
 \iftest@\nobreak\hskip-\p@\hskip\p@\fi}

\newtoks\everypartoks@
\def\noindent@@{\par\everypartoks@\expandafter{\the\everypar}\everypar{}%
 \noindent@\everypar\expandafter{\the\everypartoks@}}
\def\page{\Err@{\noexpand\page has no meaning by itself}}
\let\page@C\pageno
\let\page@P\empty
\let\page@Q\empty
\def\page@S#1{#1\/}
\def\page@F{\rm}
\def\page@N{\arabic}   
\newif\ifindexing@
\def\indexfile{\ifindexing@\else
 \alloc@@7\write\chardef\sixt@@n\ndx@
 \immediate\openout\ndx@=\jobname.ndx
 \global\indexing@true\fi}
\global\advance\insc@unt\m@ne
\ch@ck0\insc@unt\count
\ch@ck1\insc@unt\dimen
\ch@ck2\insc@unt\skip
\ch@ck4\insc@unt\box
\allocationnumber\insc@unt
\global\chardef\margin@\allocationnumber
\dimen\margin@\maxdimen
\count\margin@\z@
\skip\margin@\z@
\newif\ifindexproofing@
\def\indexproofing{\indexproofing@true}
\def\noindexproofing{\indexproofing@false}
\def\unmacro@#1:#2->#3\unmacro@{\def\macpar@{#2}\def\macdef@{#3}}
\def\starparts@#1{\def\stari@{#1}\def\starii@{#1}\let\stariii@\empty
 \test@false
 \DN@##1*##2##3\next@{\ifx\starparts@##2\test@false\else\test@true\fi}%
 \next@#1*\starparts@\next@
 \iftest@\DN@{\starparts@@#1\starparts@@}\else\let\next@\relax\fi\next@}
\def\starparts@@#1*#2\starparts@@{\def\starii@{#1}\def\stariii@{*#2}}
\def\windex@{\ifindexing@
 \expandafter\unmacro@\meaning\stari@\unmacro@
 \edef\macdef@{\string"\macdef@\string"}%
 \edef\next@{\write\ndx@{\macdef@}}\next@
 \write\ndx@{{\number\pageno}{\page@N}{\page@P}{\page@Q}}%
 \fi
 \ifindexproofing@
  \ifx\stariii@\empty\else
   \expandafter\unmacro@\meaning\stariii@\unmacro@\fi
  \insert\margin@{\hbox{\rm\vrule\height9\p@\depth2\p@\width\z@\starii@
  \ifx\stariii@\empty\else\tt\macdef@\fi}}\fi}
\catcode`\"=\active
\def"{\FN@\quote@}
\def\quote@{\ifx\next"\expandafter\quote@@\else\expandafter\quote@@@\fi}
\def\quote@@@#1"{\starparts@{#1}\starii@\windex@}
\def\quote@@"#1"{\prevanish@\starparts@{#1}\windex@\FN@\quote@@@@}
\def\quote@@@@{\ifx\next"\DN@"{\postvanish@}\else
 \let\next@\postvanish@\fi\next@}
\rightadd@"\to\vanishlist@
\def\idefine#1{\DN@{#1}\DNii@{\noexpand#1}%
 \afterassignment\idefine@\def\nextiii@}
\def\idefine@{\ifindexing@
 \expandafter\let\next@\nextiii@
 \expandafter\unmacro@\meaning\nextiii@\unmacro@
 \immediate\write\ndx@{\noexpand\define\nextii@\macpar@{\macdef@}}\fi}
\def\iabbrev*#1#2{\ifindexing@\toks@{#2}%
 \immediate\write\ndx@{\noexpand\abbrev*\noexpand#1{\the\toks@}}\fi}
\newread\laxread@
\newwrite\laxwrite@
\let\fnpages@\empty
\def\Finit@#1#2\Finit@{\let\nextii@#1\def\nextiii@{#2}}
\catcode`\~=11
\def\getparts@ @#1~#2~#3~#4~#5~#6{\def\nextiv@{#1}%
 \def\nextiii@{#2~#3~#4~#5~}\count@#6\relax}
\newif\ifdocument@
\def\document{\ifdocument@\else\global\document@true
 \let\fontlist@\empty
 \immediate\openin\laxread@=\jobname.lax\relax
 {\endlinechar\m@ne\noexpands@\catcode`\@=11 \catcode`\~=11
  \loop\ifeof\laxread@\else
   \read\laxread@ to\next@
   \ifx\next@\empty
   \else
    \expandafter\Finit@\next@\Finit@
    \if\nextii@ F%
     \expandafter\rightadd@\nextiii@\to\fnpages@
    \else
     \expandafter\getparts@\next@
     \edef\next@{\gdef\csname\nextiv@ @L\endcsname{\nextiii@\number\count@}}%
     \next@
    \fi
   \fi
  \repeat}%
 \immediate\closein\laxread@
 \immediate\openout\laxwrite@=\jobname.lax\relax\fi}
\let\thelabel@\relax
\def\thelabels@{\thelabel@ ~\thelabel@@ ~\thelabel@@@ ~\thelabel@@@@ ~}
\def\label#1{\prevanish@
 \ifx\thelabel@\relax
  \Err@{There's nothing here to be labelled}%
 \else
  {\noexpands@
  \expandafter\ifx\csname#1@L\endcsname\relax
   \expandafter\xdef\csname#1@L\endcsname{\thelabels@0}%
   \immediate\write\laxwrite@{@#1~\thelabels@1}%
  \else
   \edef\next@{@~\csname#1@L\endcsname}%
    \expandafter\getparts@\next@
    \ifodd\count@
    \expandafter\xdef\csname#1@L\endcsname{\thelabels@0}%
    \immediate\write\laxwrite@{@#1~\thelabels@1}%
   \else
    \Err@{Label #1 already used}%
   \fi
  \fi
  }%
 \fi
 \postvanish@}
\rightadd@\label\to\vanishlist@
\def\thepages@{\page@N{\number\page@C}~%
 \page@S{\page@P\page@N{\number\page@C}\page@Q}~%
 \number\page@C ~\page@P\page@N{\number\page@C}\page@Q ~}
\def\pagelabel#1{\prevanish@
 \expandafter\ifx\csname#1@L\endcsname\relax
  {\noexpands@
  \expandafter\xdef\csname#1@L\endcsname{\thepages@2}}%
  \write\laxwrite@{@#1~\thepages@3}%
 \else
  {\noexpands@
  \edef\next@{@~\csname#1@L\endcsname}%
  \expandafter\getparts@\next@
  \ifodd\count@
   \ifnum\count@=\@ne
    \expandafter\xdef\csname#1@L\endcsname{\thelabels@2}%
   \fi
   \write\laxwrite@{@#1~\thepages@3}%
  \else
   \Err@{Label #1 already used}%
  \fi
  }%
 \fi
 \postvanish@}
\rightadd@\pagelabel\to\vanishlist@
\newif\ifreferr@
\referr@true
\def\RefErrors{\global\referr@true}
\def\RefWarnings{\global\referr@false}
\setbox\z@\hbox{\global\count@=`^^30}
\ifnum\count@=48 \let\versionthree@\relax\fi
\def\nolabel@#1#2#3{\expandafter\ifx\csname#2@L\endcsname\relax
 \ifreferr@\Err@{No \noexpand\label found for #2}\else
 \W@{Warning: No \noexpand\label found for #2.}%
 \ifx\versionthree@\relax\W@{l.\number\inputlineno\space ... \string#1{#2}}\fi
 \fi#3\else}
\def\csL@#1{{\noexpands@\xdef\Next@{\csname#1@L\endcsname}}}
\def\ref#1{\nolabel@\ref{#1}\relax
 \DNii@##1~##2\nextii@{##1}%
 \csL@{#1}\expandafter\nextii@\Next@\nextii@\fi}
\def\Ref#1{\nolabel@\Ref{#1}\relax
 \DNii@##1~##2~##3\nextii@{##2}%
 \csL@{#1}\expandafter\nextii@\Next@\nextii@\fi}
\def\nref#1{\nolabel@\nref{#1}\relax
 \DNii@##1~##2~##3~##4\nextii@{##3}%
 \csL@{#1}\expandafter\nextii@\Next@\nextii@\fi}
\def\pref#1{\nolabel@\pref{#1}\relax
 \DNii@##1~##2~##3~##4~##5\nextii@{##4}%
 \csL@{#1}\expandafter\nextii@\Next@\nextii@\fi}
\let\pref@\pref
\def\Evaluatenref#1{\nolabel@\Evaluatenref{#1}{\gdef\Nref{-10000 }}%
 \DNii@##1~##2~##3~##4\nextii@{\DNii@{##3}}%
 \csL@{#1}\expandafter\nextii@\Next@\nextii@
 \xdef\Nref{\nextii@}\fi}
\def\Evaluatepref#1{\nolabel@\Evaluatepref{#1}{\global\let\Pref\empty}%
 \DNii@##1~##2~##3~##4~##5\nextii@{\DNii@{##4}}%
 \csL@{#1}\expandafter\nextii@\Next@\nextii@
 \xdef\Pref{\nextii@}\fi}
\def\readlax#1{\immediate\openin\laxread@=#1.lax\relax
 \ifeof\laxread@\W@{}\W@{File #1.lax not found.}\W@{}\fi
 {\endlinechar\m@ne\noexpands@\catcode`\@=11 \catcode`\~=11
  \loop\ifeof\laxread@\else
   \read\laxread@ to\nextv@
   \ifx\nextv@\empty
   \else
    \expandafter\Finit@\nextv@\Finit@
    \ifx\nextii@ F%
    \else
     \expandafter\getparts@\nextv@
     \expandafter\ifx\csname\nextiv@ @L\endcsname\relax
      \edef\next@{\gdef\csname\nextiv@ @L\endcsname
       {\nextiii@\ifnum\count@=\@ne0\else2\fi}}%
      \next@
     \else
      \Err@{Label \nextiv@\space in #1.lax already used}%
     \fi
    \fi
   \fi
  \repeat}%
 \immediate\closein\laxread@}
\catcode`\~=\active

\def\FNSSP@{\FNSS@\pretendspace@}
\everydisplay{\csname displaymath \endcsname}
\expandafter\def\csname displaymath \endcsname#1$${#1$$\FNSSP@}
\def\locallabel@{\let\thelabel@\Thelabel@\let\thelabel@@\Thelabel@@
 \let\thelabel@@@\Thelabel@@@\let\thelabel@@@@\Thelabel@@@@}
\newcount\tag@C
\tag@C\z@
\let\tag@P\empty
\let\tag@Q\empty
\def\tag@S#1{{\rm(}{#1\/}{\rm)}}
\let\tag@N\arabic
\def\tag@F{\rm}
\def\maketag@{\FN@\maketag@@}
\def\maketag@@{\ifx\next\relax\DN@\relax{\FN@\maketag@@}\else
 \ifx\next"\let\next@\maketag@@@\else
 \let\next@\maketag@@@@\fi\fi\next@}
\def\xdefThelabel@#1{\xdef\Thelabel@{#1{\Thelabel@@@}}}
\def\xdefThelabel@@#1{\xdef\Thelabel@@{#1{\Thelabel@@@@}}}
\def\maketag@@@@#1\maketag@{\global\advance\tag@C\@ne
 {\noexpands@
  \xdef\Thelabel@@@{\number\tag@C}%
  \xdefThelabel@\tag@N
  \xdef\Thelabel@@@@{\ifmathtags@$\tag@P\Thelabel@\tag@Q$\else
   \tag@P\Thelabel@\tag@Q\fi}%
  \xdefThelabel@@\tag@S
  }%
 \locallabel@
 \hbox{\tag@F\thelabel@@}%
 #1}
\def\Qlabel@#1{{\noexpands@\xdef\Thelabel@@{#1}%
 \let\style\empty\xdef\Thelabel@@@@{#1}%
 \let\pre\empty\let\post\empty\xdef\Thelabel@{#1}%
 \let\numstyle\empty\xdef\Thelabel@@@{#1}}}
\def\maketag@@@"#1"#2\maketag@{%
 {\let\pre\tag@P\let\post\tag@Q\let\style\tag@S\let\numstyle\tag@N
  \hbox{\tag@F#1}%
  \noexpands@
  \Qlabel@{#1}%
  }%
 \locallabel@
 #2}
\def\align@{\inalign@true\inany@true
 \vspace@\allowdisplaybreak@\displaybreak@\intertext@
 \def\tag{\global\tag@true\ifnum\and@=\z@
  \DN@{&\omit\global\rwidth@\z@&\relax}\else
  \DN@{&\relax}\fi\next@}%
 \iftagsleft@\DN@{\csname align \endcsname}\else
  \DN@{\csname align \space\endcsname}\fi\next@}
\def\noset@{\def\Offset##1##2{\prevanish@\postvanish@}%
 \def\Reset##1##2{\prevanish@\postvanish@}}
\def\measure@#1\endalign{\global\lwidth@\z@\global\rwidth@\z@
 \global\maxlwidth@\z@\global\maxrwidth@\z@
 \global\and@\z@
 \setbox\z@\vbox
  {\noset@\everycr{\noalign{\global\tag@false\global\and@\z@}}\Let@
  \halign{\setboxz@h{$\m@th\displaystyle{\@lign##}$}%
   \global\lwidth@\wdz@
   \ifdim\lwidth@>\maxlwidth@\global\maxlwidth@\lwidth@\fi
   \global\advance\and@\@ne
   &\setboxz@h{$\m@th\displaystyle{{}\@lign##}$}\global\rwidth@\wdz@
   \ifdim\rwidth@>\maxrwidth@\global\maxrwidth@\rwidth@\fi
   \global\advance\and@\@ne
   &\Tag@\eat@{##}\crcr#1\crcr}}%
 \totwidth@\maxlwidth@\advance\totwidth@\maxrwidth@}
\def\prepost@{\global\let\tag@P@\tag@P\global\let\tag@Q@\tag@Q}
\def\reprepost@{\let\tag@P\tag@P@\let\tag@Q\tag@Q@}
\expandafter\def\csname align \space\endcsname#1\endalign
 {\measure@#1\endalign\global\and@\z@
 \ifingather@\everycr{\noalign{\global\and@\z@}}\else\displ@y@\fi
 \Let@\tabskip\centering@
 \halign to\displaywidth
  {\hfil\strut@\setboxz@h{$\m@th\displaystyle{\@lign##\prepost@}$}%
  \boxz@\global\advance\and@\@ne
  \tabskip\z@skip
  &\setboxz@h{$\m@th\displaystyle{{}\@lign##\prepost@}$}%
  \global\rwidth@\wdz@\boxz@\hfil\global\advance\and@\@ne
  \tabskip\centering@
  &\setboxz@h{\@lign\strut@\reprepost@\maketag@##\maketag@}%
  \dimen@\displaywidth\advance\dimen@-\totwidth@
  \divide\dimen@\tw@\advance\dimen@\maxrwidth@\advance\dimen@-\rwidth@
  \ifdim\dimen@<\tw@\wdz@\llap{\vtop{\normalbaselines\null\boxz@}}%
  \else\llap{\boxz@}\fi
  \tabskip\z@skip
  \crcr#1\crcr
  \black@\totwidth@}}
\expandafter\def\csname align \endcsname#1\endalign{\measure@#1\endalign
 \global\and@\z@
 \ifdim\totwidth@>\displaywidth\let\displaywidth@\totwidth@\else
  \let\displaywidth@\displaywidth\fi
 \ifingather@\everycr{\noalign{\global\and@\z@}}\else\displ@y@\fi
 \Let@\tabskip\centering@\halign to\displaywidth
  {\hfil\strut@\setboxz@h{$\m@th\displaystyle{\@lign##\prepost@}$}%
  \global\lwidth@\wdz@\global\lineht@\ht\z@
  \boxz@\global\advance\and@\@ne
  \tabskip\z@skip&\setboxz@h{$\m@th\displaystyle{{}\@lign##\prepost@}$}%
  \ifdim\ht\z@>\lineht@\global\lineht@\ht\z@\fi
  \boxz@\hfil\global\advance\and@\@ne
  \tabskip\centering@&\kern-\displaywidth@
  \setboxz@h{\@lign\strut@\reprepost@\maketag@##\maketag@}%
  \dimen@\displaywidth\advance\dimen@-\totwidth@
  \divide\dimen@\tw@\advance\dimen@\maxlwidth@\advance\dimen@-\lwidth@
  \ifdim\dimen@<\tw@\wdz@
   \rlap{\vbox{\normalbaselines\boxz@\vbox to\lineht@{}}}\else
   \rlap{\boxz@}\fi
  \tabskip\displaywidth@\crcr#1\crcr\black@\totwidth@}}
\def\attag@#1{\let\Maketag@\maketag@\let\TAG@\Tag@
 \let\Prepost@\prepost@\let\Reprepost@\reprepost@
 \let\Tag@\relax\let\maketag@\relax
 \let\prepost@\relax\let\reprepost@\relax
 \ifmeasuring@
  \def\llap@##1{\setboxz@h{##1}\hbox to\tw@\wdz@{}}%
  \def\rlap@##1{\setboxz@h{##1}\hbox to\tw@\wdz@{}}%
 \else\let\llap@\llap\let\rlap@\rlap\fi
 \toks@{\hfil\strut@
  $\m@th\displaystyle{\@lign\the\hashtoks@\prepost@}$%
  \tabskip\z@skip\global\advance\and@\@ne&
  $\m@th\displaystyle{{}\@lign\the\hashtoks@\prepost@}$\hfil
  \ifxat@\tabskip\centering@\fi\global\advance\and@\@ne}%
 \iftagsleft@
  \toks@@{\tabskip\centering@&\Tag@\kern-\displaywidth
   \rlap@{\@lign\reprepost@\maketag@\the\hashtoks@\maketag@}%
   \global\advance\and@\@ne\tabskip\displaywidth}\else
  \toks@@{\tabskip\centering@&\Tag@\llap@{\@lign\reprepost@\maketag@
   \the\hashtoks@\maketag@}\global\advance\and@\@ne\tabskip\z@skip}\fi
 \atcount@#1\relax\advance\atcount@\m@ne
 \loop\ifnum\atcount@>\z@
  \toks@\expandafter{\the\toks@&\hfil$\m@th\displaystyle{\@lign
  \the\hashtoks@\prepost@}$\global\advance\and@\@ne
  \tabskip\z@skip
  &$\m@th\displaystyle{{}\@lign\the\hashtoks@\prepost@}$\hfil\ifxat@
  \tabskip\centering@\fi\global\advance\and@\@ne}\advance\atcount@\m@ne
 \repeat
 \edef\preamble@{\the\toks@\the\toks@@}%
 \edef\preamble@@{\preamble@}%
 \let\maketag@\Maketag@\let\Tag@\TAG@
 \let\prepost@\Prepost@\let\reprepost@\Reprepost@}
\def\unlabel@{\def\label##1{\prevanish@\postvanish@}%
 \def\pagelabel##1{\prevanish@\postvanish@}}
\newcount\tag@CC
\expandafter\def\csname alignat \endcsname#1#2\endalignat
 {\inany@true\xat@false
 \def\tag{\global\tag@true
  \count@#1\relax\multiply\count@\tw@\advance\count@\m@ne
  \gdef\tag@{&}%
  \loop\ifnum\count@>\and@\xdef\tag@{&\omit\tag@}%
  \advance\count@\m@ne\repeat
  \tag@\relax}%
 \vspace@\allowdisplaybreak@\displaybreak@\intertext@
 \displ@y@\measuring@true\tag@CC\tag@C
 \setbox\savealignat@\hbox{\noset@\unlabel@$\m@th\displaystyle\Let@
  \attag@{#1}\vbox{\halign{\span\preamble@@\crcr#2\crcr}}$}%
 \measuring@false
 \Let@\attag@{#1}\tag@C\tag@CC
 \tabskip\centering@\halign to\displaywidth
  {\span\preamble@@\crcr#2\crcr\black@{\wd\savealignat@}}}
\expandafter\def\csname xalignat \endcsname#1#2\endxalignat
 {\inany@true\xat@true
 \def\tag{\global\tag@true
  \count@#1\relax\multiply\count@\tw@\advance\count@\m@ne
  \gdef\tag@{&}%
  \loop\ifnum\count@>\and@\xdef\tag@{&\omit\tag@}%
  \advance\count@\m@ne\repeat
  \tag@\relax}%
 \vspace@\allowdisplaybreak@\displaybreak@\intertext@
 \displ@y@\measuring@true\tag@CC\tag@C
 \setbox\savealignat@\hbox{\noset@\unlabel@$\m@th\displaystyle\Let@
  \attag@{#1}\vbox{\halign{\span\preamble@@\crcr#2\crcr}}$}%
 \measuring@false\Let@\attag@{#1}\tag@C\tag@CC
 \tabskip\centering@\halign to\displaywidth
 {\span\preamble@@\crcr#2\crcr\black@{\wd\savealignat@}}}
\def\gather{\RIfMIfI@\DN@{\onlydmatherr@\gather}\else
 \ingather@true\inany@true\def\tag{&\relax}%
 \vspace@\allowdisplaybreak@\displaybreak@\intertext@
 \displ@y\Let@
 \iftagsleft@\DN@{\csname gather \endcsname}\else
  \DN@{\csname gather \space\endcsname}\fi\fi
 \else\DN@{\onlydmatherr@\gather}\fi\next@}
\def\exstring@{\expandafter\eat@\string}
\def\newcounter#1{\define#1{}%
 \edef\next@{\def\noexpand#1{\futurelet\noexpand\next
  \csname\exstring@#1@Z\endcsname}}\next@
 \edef\next@{\def\csname\exstring@#1@Z\endcsname
  {\global\advance\csname\exstring@#1@C\endcsname\@ne
  {\csname\exstring@#1@F\endcsname\csname\exstring@#1@S\endcsname
   {\csname\exstring@#1@P\endcsname\csname\exstring@#1@N\endcsname
   {\noexpand\number\csname\exstring@#1@C\endcsname}%
   \csname\exstring@#1@Q\endcsname}}%
  \noexpand\ifx\noexpand\next\noexpand\label
   \def\noexpand\next@\noexpand\label########1{{\noexpand\noexpands@
    \xdef\noexpand\Thelabel@{\csname\exstring@#1@N\endcsname
     {\noexpand\number\csname\exstring@#1@C\endcsname}}%
    \xdef\noexpand\Thelabel@@@{\noexpand\number
     \csname\exstring@#1@C\endcsname}%
    \xdef\noexpand\Thelabel@@{\csname\exstring@#1@S\endcsname
     {\csname\exstring@#1@P\endcsname
     \csname\exstring@#1@N\endcsname
     {\noexpand\number\csname\exstring@#1@C\endcsname}%
     \csname\exstring@#1@Q\endcsname}}%
    \xdef\noexpand\Thelabel@@@@{\csname\exstring@#1@P\endcsname
     \csname\exstring@#1@N\endcsname
     {\noexpand\number\csname\exstring@#1@C\endcsname}%
     \csname\exstring@#1@Q\endcsname}}%
    {\noexpand\locallabel@\noexpand\label{########1}}}%
   \noexpand\else\let\noexpand\next@\relax\noexpand\fi\noexpand\next@}}\next@
 \expandafter\newcount@\csname\exstring@#1@C\endcsname
 \expandafter\let\csname\exstring@#1@N\endcsname\arabic
 \expandafter\def\csname\exstring@#1@S\endcsname##1{##1\/}%
 \expandafter\let\csname\exstring@#1@P\endcsname\empty
 \expandafter\let\csname\exstring@#1@Q\endcsname\empty
 \expandafter\def\csname\exstring@#1@F\endcsname{\rm}%
 }
\def\HASH@#1#2{\ifnum#2=\z@\else
 \edef\next@{\toks@{\the\toks@\the\hashtoks@#2}%
 \toks@@{\the\toks@@{\the\hashtoks@#2}}}\next@\expandafter\HASH@\fi}
\def\HASH@@{\toks@{}\toks@@{}\expandafter\HASH@\macpar@00}
\def\usecounter#1#2{\expandafter\ifx\csname\exstring@#1@Z\endcsname
 \relax\Err@{\noexpand#1not created with \string\newcounter}\fi
 \expandafter\let\csname\exstring@#1@@Z\endcsname\relax
 \expandafter\let\csname\exstring@#1@@Z@\endcsname\relax
 \expandafter\let\csname\exstring@#1@@Z@@\endcsname\relax
 \edef\next@{\def\noexpand#2{\futurelet\noexpand\next
  \csname\exstring@#1@@Z\endcsname}}\next@
 \edef\next@{\def\csname\exstring@#1@@Z\endcsname{\noexpand\ifx
  \noexpand\next\noexpand\label\def\noexpand\next@\noexpand\label
   ########1{\csname\exstring@#1@@Z@\endcsname
   {\noexpand#1\noexpand\label{########1}}}%
   \noexpand\else\noexpand\ifx\noexpand\next
   \noexpand"\def\noexpand\next@\noexpand"########1\noexpand"%
   {\csname\exstring@#1@@Z@\endcsname{{\expandafter\noexpand
   \csname\exstring@#1@F\endcsname
   \let\noexpand\pre\expandafter\noexpand\csname\exstring@#1@P\endcsname
   \let\noexpand\post\expandafter\noexpand\csname\exstring@#1@Q\endcsname
   \let\noexpand\style\expandafter\noexpand\csname\exstring@#1@S\endcsname
   \let\noexpand\numstyle\expandafter\noexpand\csname\exstring@#1@N\endcsname
   ########1}}}\noexpand\else
   \def\noexpand\next@{\csname\exstring@#1@@Z@\endcsname{\noexpand#1}}%
   \noexpand\fi\noexpand\fi\noexpand\next@}}\next@
 \def\next@{\expandafter\expandafter\expandafter\unmacro@\expandafter
  \meaning\csname\exstring@#1@@Z@@\endcsname\unmacro@
  \HASH@@
  \edef\next@{\def\csname\exstring@#1@@Z@\endcsname\the\toks@{%
   \expandafter\noexpand\csname\exstring@#1@@Z@@\endcsname\the\toks@@
   \noexpand\FNSSP@}}\next@}%
 \afterassignment\next@
 \expandafter\def\csname\exstring@#1@@Z@@\endcsname}
\def\listbi@{\penalty50 \medskip}
\def\listbii@{\penalty100 \smallskip}
\let\listbiii@\relax
\let\listbiv@\relax
\let\listbv@\relax
\def\listmi@{\advance\leftskip30\p@\relax}
\let\listmii@\listmi@
\let\listmiii@\listmi@
\let\listmiv@\listmi@
\let\listmv@\listmi@
\def\itemi@#1{\noindent@@\llap{#1\hskip5\p@}}
\let\itemii@\itemi@
\let\itemiii@\itemi@
\let\itemiv@\itemi@
\let\itemv@\itemi@
\def\liste@{\penalty-50 \medskip}
\def\listei@{\penalty-100 \smallskip}
\let\listeii@\relax
\let\listeiii@\relax
\let\listeiv@\relax
\expandafter\newcount\csname list@C1\endcsname
\csname list@C1\endcsname\z@
\expandafter\newcount\csname list@C2\endcsname
\csname list@C2\endcsname\z@
\expandafter\newcount\csname list@C3\endcsname
\csname list@C3\endcsname\z@
\expandafter\newcount\csname list@C4\endcsname
\csname list@C4\endcsname\z@
\expandafter\newcount\csname list@C5\endcsname
\csname list@C5\endcsname\z@
\expandafter\let\csname list@P1\endcsname\empty
\expandafter\let\csname list@P2\endcsname\empty
\expandafter\let\csname list@P3\endcsname\empty
\expandafter\let\csname list@P4\endcsname\empty
\expandafter\let\csname list@P5\endcsname\empty
\expandafter\let\csname list@Q1\endcsname\empty
\expandafter\let\csname list@Q2\endcsname\empty
\expandafter\let\csname list@Q3\endcsname\empty
\expandafter\let\csname list@Q4\endcsname\empty
\expandafter\let\csname list@Q5\endcsname\empty
\expandafter\def\csname list@S1\endcsname#1{{\rm(}{#1\/}{\rm)}}
\expandafter\def\csname list@S2\endcsname#1{{\rm(}{#1\/}{\rm)}}
\expandafter\def\csname list@S3\endcsname#1{{\rm(}{#1\/}{\rm)}}
\expandafter\def\csname list@S4\endcsname#1{{\rm(}{#1\/}{\rm)}}
\expandafter\def\csname list@S5\endcsname#1{{\rm(}{#1\/}{\rm)}}
\expandafter\let\csname list@N1\endcsname\arabic
\expandafter\let\csname list@N2\endcsname\arabic
\expandafter\let\csname list@N3\endcsname\arabic
\expandafter\let\csname list@N4\endcsname\arabic
\expandafter\let\csname list@N5\endcsname\arabic
\expandafter\def\csname list@F1\endcsname{\rm}
\expandafter\def\csname list@F2\endcsname{\rm}
\expandafter\def\csname list@F3\endcsname{\rm}
\expandafter\def\csname list@F4\endcsname{\rm}
\expandafter\def\csname list@F5\endcsname{\rm}
\newcount\listlevel@
\listlevel@\z@
\def\list@@C{\csname list@C\number\listlevel@\endcsname}
\def\list@@P{\csname list@P\number\listlevel@\endcsname}
\def\list@@Q{\csname list@Q\number\listlevel@\endcsname}
\def\list@@S{\csname list@S\number\listlevel@\endcsname}
\def\list@@N{\csname list@N\number\listlevel@\endcsname}
\def\list@@F{\csname list@F\number\listlevel@\endcsname}
\newif\iffirstitemi@
\newif\iffirstitemii@
\newif\iffirstitemiii@
\newif\iffirstitemiv@
\newif\iffirstitemv@
\def\Firstitem@true{\csname firstitem\romannumeral\listlevel@
 @true\endcsname}
\def\Firstitem@false{\csname firstitem\romannumeral\listlevel@
 @false\endcsname}
\def\Listm@{\csname listm\romannumeral\listlevel@ @\endcsname}
\def\Item@{\csname item\romannumeral\listlevel@ @\endcsname}
\def\Liste@{\csname liste\romannumeral\listlevel@ @\endcsname}
\newif\iflistcontinue@
\def\keepitem{\listcontinue@true}
\newcount\list@C@
\def\list{%
 \iflistcontinue@\csname list@C1\endcsname\csname list@C@\endcsname\fi
 \global\csname list@C2\endcsname\z@
 \global\csname list@C3\endcsname\z@
 \global\csname list@C4\endcsname\z@
 \global\csname list@C5\endcsname\z@
 \begingroup
 \firstitemi@true
 \listlevel@\@ne
 \def\item{\FN@\item@}%
 \FN@\list@}
\Invalid@\runinitem
\def\list@{\ifx\next\par
 \DN@\par{\FN@\list@}\else
 \ifx\next\runinitem
  \DN@\runinitem{\FN@\runinitem@}\else
  \DN@{\par\dimen@\parskip\parskip\dimen@}\fi\fi\next@}
\newif\ifoutlevel@
\newif\ifrunin@
\def\item@{%
 \ifoutlevel@\Liste@\outlevel@false\fi
 \ifrunin@\runin@false\par
  \dimen@\parskip\parskip\dimen@
  \Listm@\fi
 \iffirstitemi@\listbi@\listmi@\firstitemi@false\else\par\fi
 \iffirstitemii@\listbii@\listmii@\firstitemii@false\else\par\fi
 \iffirstitemiii@\listbiii@\listmiii@\firstitemiii@false\else\par\fi
 \iffirstitemiv@\listbiv@\listmiv@\firstitemiv@false\else\par\fi
 \iffirstitemv@\listbv@\listmv@\firstitemv@false\else\par\fi
 \DN@"##1"{{\let\pre\list@@P\let\post\list@@Q
  \let\style\list@@S\let\numstyle\list@@N
  \vskip-\parskip
  \Item@{\list@@F##1}%
  \noexpands@
  \Qlabel@{##1}}%
  \locallabel@
  \FNSSP@}%
 \DNii@{\global\advance\list@@C\@ne
  {\noexpands@
   \xdef\Thelabel@@@{\number\list@@C}%
   \xdefThelabel@\list@@N
   \xdef\Thelabel@@@@{\list@@P\Thelabel@\list@@Q}%
   \xdefThelabel@@\list@@S
  }%
  \locallabel@
  \vskip-\parskip
  \Item@{\list@@F\thelabel@@}%
  \FN@\pretendspace@}%
 \ifx\next"\expandafter\next@\else\expandafter\nextii@\fi}
\def\runinitem@{%
  \runin@true
  \Firstitem@false
  \DN@"##1"{{\let\pre\list@@P\let\post\list@@Q
   \let\style\list@@S\let\numstyle\list@@N
   \unskip\space{\list@@F##1} %
   \noexpands@
   \Qlabel@{##1}}%
   \locallabel@
   \ignorespaces}%
  \DNii@{\global\advance\list@@C\@ne
   {\noexpands@
    \xdef\Thelabel@@@{\number\list@@C}%
    \xdefThelabel@\list@@N
    \xdef\Thelabel@@@@{\list@@P\Thelabel@\list@@Q}%
    \xdefThelabel@@\list@@S
   }%
   \locallabel@
   \unskip\space{\list@@F\thelabel@@} }%
  \ifx\next"\expandafter\next@\else\expandafter\nextii@\fi}
\def\inlevel{\ifnum\listlevel@=5
 \DN@{\Err@{Already 5 levels down}}\else
 \DN@{\begingroup\advance\listlevel@\@ne
 \Firstitem@true\FN@\inlevel@}\fi\next@}
\def\inlevel@{\ifx\next\par
 \DN@\par{\FN@\inlevel@}\else
 \ifx\next\runinitem
  \DN@\runinitem{\FN@\runinitem@}\else
  \let\next@\relax\fi\fi\next@}
\def\outlevel{\ifnum\listlevel@=\@ne
 \Err@{At top level}\else
 \par\global\list@@C\z@\endgroup\outlevel@true\fi}
\def\endlist{%
 \expandafter\global\csname list@C@\endcsname\csname list@C1\endcsname
 \par
 \global\toks\@ne{}\count@\listlevel@
 {\loop
  \ifnum\count@>\z@\global\toks\@ne\expandafter{\the\toks\@ne\endgroup}%
  \advance\count@\m@ne
  \repeat}%
 \the\toks\@ne
 \liste@
 \listcontinue@false\global\csname list@C1\endcsname\z@
 \vskip-\parskip
 \noindent@@
 \FN@\pretendspace@}
\newif\iffirstdescribe@
\def\describe{\par
 \begingroup\firstdescribe@true
 \def\item##1{%
  \iffirstdescribe@\penalty50 \medskip\vskip-\parskip
  \firstdescribe@false\else\par\fi
  \noindent@@\hangindent2pc\hangafter\@ne
  {\bf##1}\hskip.5em}}

\Invalid@\pullin
\Invalid@\pullinmore
\newif\iffirstpull@
\def\margins{\par\begingroup\firstpull@true
 \def\pullin##1##2{\par
  \iffirstpull@\firstpull@false\else\endgroup\fi
  \begingroup\DN@{##1}%
  \ifx\next@\empty\leftskip\z@\else\ifx\next@\space\leftskip\z@
  \else\leftskip##1\fi\fi
  \DN@{##2}\ifx\next@\empty\rightskip\z@\else\ifx\next@\space
  \rightskip\z@\else\rightskip##2\fi\fi\ignorespaces}%
 \def\pullinmore##1##2{\par
  \xdef\Next@{\leftskip\the\leftskip\relax\rightskip\the\rightskip\relax}%
  \iffirstpull@\firstpull@false\else\endgroup\fi
  \begingroup\Next@
  \DN@{##1}%
  \ifx\next@\empty\else\ifx\next@\space\else\advance\leftskip##1\fi\fi
  \DN@{##2}\ifx\next@\empty\else\ifx\next@\space\else
  \advance\rightskip##2\fi\fi\ignorespaces}}

\newif\ifnopunct@
\newif\ifnospace@
\newif\ifoverlong@
\let\nofrillslist@\empty
\let\overlonglist@\empty
\def\nopunct{\nopunct@true\FN@\nopunct@}
\def\nospace{\nospace@true\FN@\nospace@}
\def\overlong{\overlong@true\FN@\overlong@}
\def\nopunct@{\ifx\next\nospace
 \DN@\nospace{\nospace@true\FN@\nopnos@}\else\ifx\next\overlong
 \DN@\overlong{\overlong@true\FN@\nopol@}\else
 \let\next@\nopunct@@\fi\fi\next@}
\def\nopunct@@#1{\ismember@\nofrillslist@#1%
 \iftest@\let\next@#1\else
 \DN@{\nopunct@false\Err@{\noexpand\nopunct can't be used with
 \string#1}#1}\fi\next@}
\def\nospace@{\ifx\next\nopunct
 \DN@\nopunct{\nopunct@true\FN@\nopnos@}\else\ifx\next\overlong
 \DN@\overlong{\overlong@true\FN@\nosol@}\else
 \let\next@\nospace@@\fi\fi\next@}
\def\nospace@@#1{\ismember@\nofrillslist@#1%
 \iftest@\let\next@#1\else
 \DN@{\nospace@false\Err@{\noexpand\nospace can't be used with
 \string#1}#1}\fi\next@}
\def\overlong@{\ifx\next\nopunct
 \DN@\nopunct{\nopunct@true\FN@\nopol@}\else\ifx\next\nospace
 \DN@\nospace{\nospace@true\FN@\nosol@}\else
 \let\next@\overlong@@\fi\fi\next@}
\def\overlong@@#1{\ismember@\overlonglist@#1%
 \iftest@\let\next@#1\else
 \DN@{\overlong@false\Err@{\noexpand\overlong can't be used with
 \string#1}#1}\fi\next@}
\def\nopnos@{\ifx\next\overlong
 \DN@\overlong{\overlong@true\nopnosol@}\else
 \let\next@\nopnos@@\fi\next@}
\def\nopol@{\ifx\next\nospace
 \DN@\nospace{\nospace@true\nopnosol@}\else
 \let\next@\nopol@@\fi\next@}
\def\nosol@{\ifx\next\nopunct
 \DN@\nopunct{\nopunct@true\nopnosol@}\else
 \let\next@\nosol@@\fi\next@}
\def\nopnos@@#1{\ismember@\nofrillslist@#1%
 \iftest@\let\next@#1\else
 \DN@{\nopunct@false\nospace@false
  \Err@{\noexpand\nopunct\noexpand\nospace
   can't be used with \string#1}#1}\fi\next@}
\def\testii@#1{\ismember@\nofrillslist@#1%
 \iftest@\let\nextiii@ T\else\let\nextiii@ F\fi
 \ismember@\overlonglist@#1%
 \iftest@\let\nextiv@ T\else\let\nextiv@ F\fi
 \test@false\if\nextiii@ T\if\nextiv@ T\test@true\fi\fi}
\def\nopol@@#1{\testii@{#1}%
 \iftest@\let\next@#1%
 \else\DN@{\if\nextiii@ T\else\nopunct@false\fi
  \if\nextiv@ T\else\overlong@false\fi
  \Err@{\if\nextiii@ T\else\noexpand\nopunct\fi
  \if\nextiv@ T\else\noexpand\overlong\fi can't be used
  with \string#1}#1}\fi\next@}
\def\nosol@@#1{\testii@{#1}%
 \iftest@\let\next@#1%
 \else\DN@{\if\nextiii@ T\else\nospace@false\fi
  \if\nextiv@ T\else\overlong@false\fi
  \Err@{\if\nextiii@ T\else\noexpand\nospace\fi
  \if\nextiv@ T\else\noexpand\overlong\fi can't be used
  with \string#1}#1}\fi\next@}
\def\nopnosol@#1{\testii@{#1}%
 \iftest@\let\next@#1%
 \else\DN@{\if\nextiii@ T\else\nopunct@false\nospace@false\fi
  \if\nextiv@ T\else\overlong@false\fi
  \Err@{\if\nextiii@ T\else\noexpand\nopunct\noexpand\nospace\fi
  \if\nextiv@ T\else\noexpand\overlong\fi can't be used
  with \string#1}#1}\fi\next@}
\def\punct@#1{\ifnopunct@\else#1\fi}
\def\addspace@#1{\ifnospace@\else#1\fi}
\def\hss@{\ifoverlong@\z@ plus\@m\p@ minus\@m\p@
 \else \z@ plus\@m\p@\fi}
\rightadd@\demo\to\nofrillslist@
\newif\ifclaim@
\def\exxx@{\expandafter\expandafter\expandafter\eat@\expandafter\string}
\let\colon@:
\def\demo#1{\ifclaim@
 \Err@{Previous \expandafter\noexpand\claimtype@ has
  no matching \string\end\exxx@\claimtype@}%
 \let\next@\relax
 \else
  \par
  \ifdim\lastskip<\smallskipamount\removelastskip\smallskip\fi
  \begingroup
  \noindent@@{\smc\ignorespaces#1\unskip
   \punct@{\null\colon@}\addspace@\enspace}%
  \nopunct@false\nospace@false
  \rm
  \DN@{\FNSSP@}%
 \fi
 \next@}
\def\enddemo{\par\endgroup\nopunct@false\nospace@false\smallskip}
\rightadd@\claim\to\nofrillslist@
\def\claim@F{\smc}
\def\claim@@@F{\csname\exxx@\claimtype@ @F\endcsname}
\def\claimformat@#1#2#3{%
 \medbreak\noindent@@{\smc#1 {\claim@@@F#2} #3%
 \punct@{\null.}\addspace@\enspace}\sl}
\def\claimformat@@#1#2{\claimformat@{\ignorespaces#1\unskip}%
 {\ifx\thelabel@@\empty\unskip\else\thelabel@@\fi}%
 {\ignorespaces#2\unskip}%
 \let\Claimformat@@\claimformat@@\FNSSP@}
\let\Claimformat@@\claimformat@@
\def\claim@@@P{\csname\exxx@\claimtype@ @P\endcsname}
\def\claim@@@Q{\csname\exxx@\claimtype@ @Q\endcsname}
\def\claim@@@S{\csname\exxx@\claimtype@ @S\endcsname}
\def\claim@@@N{\csname\exxx@\claimtype@ @N\endcsname}
\def\claim@@@C{\csname claim@C\claimclass@\endcsname}
\newcount\claim@C
\claim@C\z@
\let\claim@P\empty
\let\claim@Q\empty
\def\claim@S#1{#1\/}
\let\claim@N\arabic
\def\claim{\claim@true\let\claimclass@\empty
 \def\claimtype@{\claim}\FN@\claim@}
\def\claim@{%
 \ifx\next\c
  \let\next@\claim@c
 \else
  \ifx\next"%
   \let\next@\claim@q
  \else
   \begingroup\global\advance\claim@C\@ne
   {\noexpands@
    \xdef\Thelabel@@@{\number\claim@C}%
    \xdefThelabel@\claim@N
    \xdef\Thelabel@@@@{\claim@P\Thelabel@\claim@Q}%
    \xdefThelabel@@\claim@S
   }%
   \locallabel@
   \let\next@\Claimformat@@
  \fi
 \fi
 \next@}
\def\claim@c\c#1{\claim@true\begingroup
 \expandafter
 \ifx\csname claim@C#1\endcsname\relax
  \expandafter\newcount@\csname claim@C#1\endcsname
  \global\csname claim@C#1\endcsname\@ne
 \else
  \global\advance\csname claim@C#1\endcsname\@ne
 \fi
 \def\claimclass@{#1}%
 {\noexpands@
  \xdef\Thelabel@@@{\number\claim@@@C}%
  \xdefThelabel@\claim@@@N
  \xdef\Thelabel@@@@{\claim@@@P\Thelabel@\claim@@@Q}%
  \xdefThelabel@@\claim@@@S
 }%
 \locallabel@
 \FNSS@\claim@c@}
\def\claim@q"#1"{\begingroup
 {\let\pre\claim@@@P\let\post\claim@@@Q
  \let\style\claim@@@S\let\numstyle\claim@@@N
  \noexpands@
  \Qlabel@{#1}}%
 \locallabel@
 \FNSS@\claim@q@}
\def\claim@c@{\ifx\next"%
 \global\advance\claim@@@C\m@ne\let\next@\claim@cq
 \else\let\next@\Claimformat@@\fi\next@}
\def\claim@cq"#1"{{\let\pre\claim@@@P\let\post\claim@@@Q
 \let\style\claim@@@S\let\numstyle\claim@@@N
 \noexpands@
 \Qlabel@{#1}}%
 \locallabel@
 \FNSS@\Claimformat@@}
\def\claim@q@{\ifx\next\c\expandafter\claim@qc
 \else\expandafter\Claimformat@@\fi}
\def\claim@qc\c#1{\expandafter\ifx\csname claim@C#1\endcsname\relax
 \expandafter\newcount@\csname claim@C#1\endcsname
 \global\csname claim@C#1\endcsname\z@\fi
 \FNSS@\Claimformat@@}
\def\endclaim{\endgroup\claim@false\nopunct@false\nospace@false
 \let\Claimformat@@\claimformat@@\medbreak}
\Invalid@\claimclause
\def\newclaim{\FN@\newclaim@}
\def\newclaim@{\ifx\next\claimclause
 \DN@\claimclause##1{\newclaim@@{##1}}\else
 \DN@{\newclaim@@\relax}\fi\next@}
\def\claimlist@{\\\claim}
\newtoks\claim@i
\newtoks\claim@v
\let\noclaimclause@=F
\def\newclaim@@#1#2#3\c#4#5{\define#2{}%
 \rightadd@#2\to\claimlist@\rightadd@#2\to\nofrillslist@%
 \expandafter\def\csname\exstring@#2@P\endcsname{\claim@P}%
 \expandafter\def\csname\exstring@#2@Q\endcsname{\claim@Q}%
 \expandafter\def\csname\exstring@#2@S\endcsname{\claim@S}%
 \expandafter\def\csname\exstring@#2@N\endcsname{\claim@N}%
 \expandafter\def\csname\exstring@#2@F\endcsname{\claim@F}%
 \expandafter\def\csname end\exstring@#2\endcsname{\endclaim}%
 \expandafter\ifx\csname claim@C#4\endcsname\relax
  \expandafter\newcount@\csname claim@C#4\endcsname
  \global\csname claim@C#4\endcsname\z@\fi
 \edef\next@{\let\csname\exstring@#2@C\endcsname
   \csname claim@C#4\endcsname}\next@
 \def#2{\ifx\noclaimclause@ T\else#1\fi
  \global\claim@i{#1}\gdef\claim@iv{#4}\global\claim@v{#5}%
  \def\claimtype@{#2}\def\Claimformat@@{\claimformat@@{#5}}\claim@c\c{#4}}}
\def\shortenclaim#1#2{\define#2{}%
 \ismember@\claimlist@#1%
 \iftest@
  \rightadd@#2\to\nofrillslist@%
  \expandafter\def\csname\exstring@#2@P\endcsname
   {\csname\exstring@#1@P\endcsname}%
  \expandafter\def\csname\exstring@#2@Q\endcsname
   {\csname\exstring@#1@Q\endcsname}%
  \expandafter\def\csname\exstring@#2@S\endcsname
   {\csname\exstring@#1@S\endcsname}%
  \expandafter\def\csname\exstring@#2@N\endcsname
   {\csname\exstring@#1@N\endcsname}%
  \expandafter\def\csname\exstring@#2@F\endcsname
   {\csname\exstring@#1@F\endcsname}%
  \expandafter\def\csname end\exstring@#2\endcsname{\endclaim}%
  \edef\next@{\let\csname\exstring@#2@C\endcsname
    \csname claim\exstring@#1C\endcsname}\next@
  \setbox\z@\vbox{\let\noclaimclause@ T#1""\relax\endgroup}%
  \edef#2{\the\claim@i
   \def\noexpand\claimtype@{\noexpand#2}%
   \def\noexpand\Claimformat@@{\noexpand\claimformat@@{\the\claim@v}\relax}%
   \noexpand\claim@c\noexpand\c{\claim@iv}}%
 \else
  \Err@{\noexpand#1not yet created by \string\newclaim}%
 \fi}
\def\classtest@#1{\DN@{#1}\ifx\next@\claimclass@
 \test@true\else\test@false\fi}
\def\typetest@#1{\DN@{#1}\ifx\next@\claimtype@\test@true\else
  \test@false\fi}
\newif\iftoc@
\def\tocfile{\iftoc@\else\alloc@@7\write\chardef\sixt@@n\toc@
 \immediate\openout\toc@=\jobname.toc
 \alloc@@7\write\chardef\sixt@@n\tic@
 \immediate\openout\tic@=\jobname.tic
 \global\toc@true\fi}
\rightadd@\hl\to\nofrillslist@
\rightadd@\HL\to\overlonglist@
\def\HL@@C{\csname HL@C\HLlevel@\endcsname}
\def\HL@@P{\csname HL@P\HLlevel@\endcsname}
\def\HL@@Q{\csname HL@Q\HLlevel@\endcsname}
\def\HL@@S{\csname HL@S\HLlevel@\endcsname}
\def\HL@@N{\csname HL@N\HLlevel@\endcsname}
\def\HL@@F{\csname HL@F\HLlevel@\endcsname}
\def\HL@@@C{\csname\exxx@\HLtype@ @C\endcsname}
\def\HL@@@P{\csname\exxx@\HLtype@ @P\endcsname}
\def\HL@@@Q{\csname\exxx@\HLtype@ @Q\endcsname}
\def\HL@@@S{\csname\exxx@\HLtype@ @S\endcsname}
\def\HL@@@N{\csname\exxx@\HLtype@ @N\endcsname}
\def\HL#1{\expandafter
 \ifx\csname HL@C#1\endcsname\relax
  \DN@{\Err@{\string\HL#1 not defined in this style}}%
 \else
  \DN@{\gdef\HLlevel@{#1}\def\HLname@{\HL{#1}}\let\HLtype@\relax\FNSS@\HL@}%
 \fi
 \next@}%
\newif\ifquoted@
\let\aftertoc@\relax
\def\HL@{%
 \DN@"##1"##2\endHL{\def\entry@{##2}\quoted@true
  {\noexpands@
  \ifx\HLtype@\relax
   \let\pre\HL@@P\let\post\HL@@Q\let\style\HL@@S\let\numstyle\HL@@N
  \else
   \let\pre\HL@@@P\let\post\HL@@@Q\let\style\HL@@@S\let\numstyle\HL@@@N
  \fi
  \Qlabel@{##1}\let\style\relax\xdef\Qlabel@@@@{##1}%
  \xdef\Thepref@{\Thelabel@@@@}}%
  \csname HL@\HLlevel@\endcsname##2\endHL
  \let\pref\Thepref@
  \csname HL@I\HLlevel@\endcsname
  \csname HL@J\HLlevel@\endcsname
  \let\pref\pref@
  \HLtoc@	
  \aftertoc@
  \let\aftertoc@\relax\overlong@false}%
 \DNii@##1\endHL{\def\entry@{##1}\quoted@false
  {\noexpands@
  \ifx\HLtype@\relax
   \global\advance\HL@@C\@ne
   \xdef\Thelabel@@@{\number\HL@@C}%
   \xdefThelabel@{\HL@@N}%
   \xdef\Thelabel@@@@{\HL@@P\Thelabel@\HL@@Q}%
   \xdefThelabel@@{\HL@@S}%
  \else
   \global\advance\HL@@@C\@ne
   \xdef\Thelabel@@@{\number\HL@@@C}%
   \xdefThelabel@{\HL@@@N}%
   \xdef\Thelabel@@@@{\HL@@@P\Thelabel@\HL@@@Q}%
   \xdefThelabel@@{\HL@@@S}%
  \fi
  \xdef\Thepref@{\Thelabel@@@@}}%
  \csname HL@\HLlevel@\endcsname##1\endHL
  \let\pref\Thepref@
  \csname HL@I\HLlevel@\endcsname
  \csname HL@J\HLlevel@\endcsname
  \let\pref\pref@
  \HLtoc@
  \aftertoc@
  \let\aftertoc@\relax\overlong@false}%
 \ifx\next"\expandafter\next@\else\expandafter\nextii@\fi}%
\Invalid@\endHL
\def\hl@@C{\csname hl@C\hllevel@\endcsname}
\def\hl@@P{\csname hl@P\hllevel@\endcsname}
\def\hl@@Q{\csname hl@Q\hllevel@\endcsname}
\def\hl@@S{\csname hl@S\hllevel@\endcsname}
\def\hl@@N{\csname hl@N\hllevel@\endcsname}
\def\hl@@F{\csname hl@F\hllevel@\endcsname}
\def\hl@@@C{\csname\exxx@\hltype@ @C\endcsname}
\def\hl@@@P{\csname\exxx@\hltype@ @P\endcsname}
\def\hl@@@Q{\csname\exxx@\hltype@ @Q\endcsname}
\def\hl@@@S{\csname\exxx@\hltype@ @S\endcsname}
\def\hl@@@N{\csname\exxx@\hltype@ @N\endcsname}
\def\hl#1{\expandafter
 \ifx\csname hl@C#1\endcsname\relax
  \DN@{\Err@{\string\hl#1 not defined in this style}}%
 \else
  \DN@{\gdef\hllevel@{#1}\def\hlname@{\hl{#1}}\let\hltype@\relax\FNSS@\hl@}%
 \fi
 \next@}
\def\hl@{%
 \DN@"##1"##2{\def\entry@{##2}\quoted@true
  {\noexpands@
  \ifx\hltype@\relax
   \let\pre\hl@@P\let\post\hl@@Q\let\style\hl@@S\let\numstyle\hl@@N
  \else
   \let\pre\hl@@@P\let\post\hl@@@Q\let\style\hl@@@S\let\numstyle\hl@@@N
  \fi
  \Qlabel@{##1}\let\style\relax\xdef\Qlabel@@@@{##1}%
  \xdef\Thepref@{\Thelabel@@@@}}%
  \csname hl@\hllevel@\endcsname{##2}%
  \let\pref\Thepref@
  \csname hl@I\hllevel@\endcsname
  \csname hl@J\hllevel@\endcsname
  \let\pref\pref@
  \hltoc@
  \aftertoc@
  \let\aftertoc@\relax\nopunct@false\nospace@false\FNSSP@}%
 \DNii@##1{\def\entry@{##1}\quoted@false
  {\noexpands@
  \ifx\hltype@\relax
   \global\advance\hl@@C\@ne
   \xdef\Thelabel@@@{\number\hl@@C}%
   \xdefThelabel@{\hl@@N}%
   \xdef\Thelabel@@@@{\hl@@P\Thelabel@\hl@@Q}%
   \xdefThelabel@@{\hl@@S}%
  \else
   \global\advance\hl@@@C\@ne
   \xdef\Thelabel@@@{\number\hl@@@C}%
   \xdefThelabel@{\hl@@@N}%
   \xdef\Thelabel@@@@{\hl@@@P\Thelabel@\hl@@@Q}%
   \xdefThelabel@@{\hl@@@S}%
  \fi
  \xdef\Thepref@{\Thelabel@@@@}}%
  \csname hl@\hllevel@\endcsname{##1}%
  \let\pref\Thepref@
  \csname hl@I\hllevel@\endcsname
  \csname hl@J\hllevel@\endcsname
  \let\pref\pref@
  \hltoc@
  \aftertoc@
  \let\aftertoc@\relax\nopunct@false\nospace@false\FNSSP@}%
 \ifx\next"\expandafter\next@\else\expandafter\nextii@\fi}%
\def\six@#1#2 #3 #4 #5 #6 #7 {\DN@{#2}\ifx\next@\empty
 \DN@##1\six@{}\else
 \write#1{ #2 #3 #4 #5 #6 #7}\DN@{\six@#1}\fi
 \next@}
\def\Sixtoc@{\ifx\macdef@\empty\else
 \DN@##1##2\next@{\def\macdef@{##1##2}}%
 \expandafter\next@\macdef@\next@
 \edef\next@
  {\noexpand\six@\toc@\macdef@
  \space\space\space\space\space\space\space\space\space\space\space\space
  \noexpand\six@}%
 \next@\let\macdef@\relax\fi}
\def\QorThelabel@@@@{\ifquoted@
 \noexpand\noexpand\noexpand"\Qlabel@@@@\noexpand\noexpand\noexpand"\else
 \Thelabel@@@@\fi}
\def\HLtoc@{%
 \iftoc@
 \expandafter\expandafter\expandafter\unmacro@
  \expandafter\meaning\csname HL@W\HLlevel@\endcsname\unmacro@
  {\noexpands@\let\style\relax
   \edef\next@{\write\toc@{\noexpand\noexpand\expandafter\noexpand\HLname@
   {\macdef@}{\QorThelabel@@@@}}}%
  \next@}%
  \expandafter\unmacro@\meaning\entry@\unmacro@
  \Sixtoc@
  \write\toc@{\noexpand\Page{\number\pageno}{\page@N}%
   {\page@P}{\page@Q}^^J}%
 \fi}
\def\hltoc@{%
 \iftoc@
 \expandafter\expandafter\expandafter\unmacro@
  \expandafter\meaning\csname hl@W\hllevel@\endcsname\unmacro@
  {\noexpands@\let\style\relax
  \edef\next@{\write\toc@{%
   \ifnopunct@\noexpand\noexpand\noexpand\nopunct\fi
   \ifnospace@\noexpand\noexpand\noexpand\nospace\fi
   \noexpand\noexpand\expandafter\noexpand\hlname@
   {\macdef@}{\QorThelabel@@@@}}}%
  \next@}%
  \expandafter\unmacro@\meaning\entry@\unmacro@
  \Sixtoc@
  \write\toc@{\noexpand\Page{\number\pageno}{\page@N}%
   {\page@P}{\page@Q}^^J}%
 \fi}
\def\mainfile#1{\def\mainfile@{#1}}
\def\checkmainfile@{\ifx\mainfile@\undefined
 \Err@{No \noexpand\mainfile specified}\fi}
\expandafter\newcount@\csname HL@C1\endcsname
\csname HL@C1\endcsname\z@
\expandafter\def\csname HL@S1\endcsname#1{#1\null.}
\expandafter\let\csname HL@N1\endcsname\arabic
\expandafter\let\csname HL@P1\endcsname\empty
\expandafter\let\csname HL@Q1\endcsname\empty
\expandafter\def\csname HL@F1\endcsname{\bf}
\expandafter\let\csname HL@W1\endcsname\empty
\expandafter\newcount@\csname hl@C1\endcsname
\csname hl@C1\endcsname\z@
\expandafter\def\csname hl@S1\endcsname#1{#1\/}
\expandafter\let\csname hl@N1\endcsname\arabic
\expandafter\let\csname hl@P1\endcsname\empty
\expandafter\let\csname hl@Q1\endcsname\empty
\expandafter\def\csname hl@F1\endcsname{\bf}
\expandafter\let\csname hl@W1\endcsname\empty
\expandafter\def\csname HL@1\endcsname#1\endHL{\bigbreak
 {\locallabel@
  \global\setbox\@ne\vbox{\Let@\tabskip\hss@
  \halign to\hsize{\bf\hfil\ignorespaces##\unskip\hfil\cr
  \expandafter\ifx\csname HL@W1\endcsname\empty\else
   \csname HL@W1\endcsname\space\fi
  {\HL@@F\ifx\thelabel@@\empty\else\thelabel@@\space\fi}%
  \ignorespaces#1\crcr}}%
  }%
 \unvbox\@ne\nobreak\medskip}
\expandafter\def\csname hl@1\endcsname#1{\medbreak\noindent@@
 {\locallabel@
 \bf{\hl@@F\ifx\thelabel@@\empty\else\thelabel@@\space\fi}%
 \ignorespaces#1\unskip\punct@{\null.}\addspace@\enspace}}
\expandafter\def\csname HL@I1\endcsname{\Reset\hl1{1}%
 \ifx\pref\empty\newpre\hl1{}\else\newpre\hl1{\pref.}\fi}
\def\NameHL#1#2{\define#2{}%
 \expandafter\ifx\csname HL@R#1\endcsname\relax
 \else
  \def\nextiv@{\let\nextiii@}%
  \expandafter\nextiv@\csname HL@R#1\endcsname
  \expandafter\let\nextiii@\undefined
  \expandafter\let\csname\exxx@\nextiii@ @C\endcsname\relax
  \expandafter\let\csname\exxx@\nextiii@ @P\endcsname\relax
  \expandafter\let\csname\exxx@\nextiii@ @Q\endcsname\relax
  \expandafter\let\csname\exxx@\nextiii@ @S\endcsname\relax
  \expandafter\let\csname\exxx@\nextiii@ @N\endcsname\relax
  \expandafter\let\csname\exxx@\nextiii@ @F\endcsname\relax
  \expandafter\let\csname\exxx@\nextiii@ @W\endcsname\relax
  \expandafter\let\csname end\exxx@\nextiii@\endcsname\undefined
 \fi
 \expandafter\gdef\csname HL@R#1\endcsname{#2}%
 \expandafter\gdef\csname\exstring@#2@R\endcsname{{HL}{#1}}%
 \iftoc@\write\toc@{\noexpand\NameHL#1\noexpand#2^^J}\fi
 \rightadd@#2\to\overlonglist@
 \edef\next@{\let\csname\exstring@#2@C\endcsname\expandafter\noexpand
  \csname HL@C#1\endcsname}\next@
 \edef\next@{\let\csname\exstring@#2@P\endcsname\expandafter\noexpand
  \csname HL@P#1\endcsname}\next@
 \edef\next@{\let\csname\exstring@#2@Q\endcsname\expandafter\noexpand
  \csname HL@Q#1\endcsname}\next@
 \edef\next@{\let\csname\exstring@#2@S\endcsname\expandafter\noexpand
  \csname HL@S#1\endcsname}\next@
 \edef\next@{\let\csname\exstring@#2@N\endcsname\expandafter\noexpand
  \csname HL@N#1\endcsname}\next@
 \edef\next@{\let\csname\exstring@#2@F\endcsname\expandafter\noexpand
  \csname HL@F#1\endcsname}\next@
 \edef\next@{\let\csname\exstring@#2@W\endcsname\expandafter\noexpand
  \csname HL@W#1\endcsname}\next@
 \edef\next@{\def\noexpand#2####1\expandafter\noexpand
  \csname end\exstring@#2\endcsname
  {\def\noexpand\HLtype@{\noexpand#2}%
   \def\noexpand\HLname@{\noexpand#2}%
   \gdef\noexpand\HLlevel@{#1}%
   \noexpand\FNSS@\noexpand\HL@####1\noexpand\endHL}}%
  \next@
 \edef\next@{\noexpand\Invalid@\expandafter\noexpand
  \csname end\exstring@#2\endcsname}%
 \next@}
\def\Namehl#1#2{\define#2{}%
 \expandafter\ifx\csname hl@R#1\endcsname\relax
 \else
  \def\nextiv@{\let\nextiii@}%
  \expandafter\nextiv@\csname hl@R#1\endcsname
  \expandafter\let\nextiii@\undefined
  \expandafter\let\csname\exxx@\nextiii@ @C\endcsname\relax
  \expandafter\let\csname\exxx@\nextiii@ @P\endcsname\relax
  \expandafter\let\csname\exxx@\nextiii@ @Q\endcsname\relax
  \expandafter\let\csname\exxx@\nextiii@ @S\endcsname\relax
  \expandafter\let\csname\exxx@\nextiii@ @N\endcsname\relax
  \expandafter\let\csname\exxx@\nextiii@ @F\endcsname\relax
  \expandafter\let\csname\exxx@\nextiii@ @W\endcsname\relax
 \fi
 \expandafter\gdef\csname hl@R#1\endcsname{#2}%
 \expandafter\gdef\csname\exstring@#2@R\endcsname{{hl}{#1}}%
 \iftoc@\write\toc@{\noexpand\Namehl#1\noexpand#2^^J}\fi
 \rightadd@#2\to\nofrillslist@%
 \edef\next@{\let\csname\exstring@#2@C\endcsname\expandafter\noexpand
  \csname hl@C#1\endcsname}\next@
 \edef\next@{\let\csname\exstring@#2@P\endcsname\expandafter\noexpand
  \csname hl@P#1\endcsname}\next@
 \edef\next@{\let\csname\exstring@#2@Q\endcsname\expandafter\noexpand
  \csname hl@Q#1\endcsname}\next@
 \edef\next@{\let\csname\exstring@#2@S\endcsname\expandafter\noexpand
  \csname hl@S#1\endcsname}\next@
 \edef\next@{\let\csname\exstring@#2@N\endcsname\expandafter\noexpand
  \csname hl@N#1\endcsname}\next@
 \edef\next@{\let\csname\exstring@#2@F\endcsname\expandafter\noexpand
  \csname hl@F#1\endcsname}\next@
 \edef\next@{\let\csname\exstring@#2@W\endcsname\expandafter\noexpand
  \csname hl@W#1\endcsname}\next@
 \edef\next@{\def\noexpand#2{%
  \def\noexpand\hltype@{\noexpand#2}%
  \def\noexpand\hlname@{\noexpand#2}%
  \gdef\noexpand\hllevel@{#1}%
  \noexpand\FNSS@\noexpand\hl@}}%
 \next@}%
\def\Initialize{\FN@\Init@}
\def\Init@{\ifx\next\HL\let\next@\InitH@\else\ifx\next\hl\let\next@\InitH@
  \else\let\next@\InitS@\fi\fi\next@}
\def\InitH@#1#2{\expandafter\ifx\csname\exstring@#1@C#2\endcsname\relax
 \DN@{\Err@{\noexpand#1level #2 not defined in this style}}\else
 \DN@{\expandafter\gdef\csname\exstring@#1@J#2\endcsname}\fi\next@}
\def\InitC@#1#2{\edef\nextii@{\expandafter\noexpand\csname#1\endcsname{#2}}}
\def\InitS@#1{\expandafter\ifx\csname\exstring@#1@R\endcsname\relax
 \Err@{\noexpand#1not defined in this style}\let\next@\relax\else
 \DN@{\let\next@}\expandafter\next@\csname\exstring@#1@R\endcsname
 \expandafter\InitC@\next@
 \DN@{\expandafter\InitH@\nextii@}\fi\next@}
\def\value#1{\expandafter
 \ifx\csname\exstring@#1@C\endcsname\relax
  \expandafter\ifx\csname\exstring@#1@C1\endcsname\relax
   \DN@{\Err@{\noexpand\value can't be used with \string#1}}%
  \else
   \DN@{\value@#1}%
  \fi
 \else
  \DN@{\number\csname\exstring@#1@C\endcsname\relax}%
 \fi
 \next@}
\def\value@#1#2{\expandafter
 \ifx\csname\exstring@#1@C#2\endcsname\relax
  \DN@{\Err@{\string\value\string#1 can't be followed by \string#2}}%
 \else
  \DN@{\number\csname\exstring@#1@C#2\endcsname\relax}%
 \fi
 \next@}
\newcount\Value
\def\Evaluate#1{\expandafter
 \ifx\csname\exstring@#1@C\endcsname\relax
  \expandafter\ifx\csname\exstring@#1@C1\endcsname\relax
   \DN@{\Err@{\noexpand\Evaluate can't be used with \string#1}}%
  \else
   \DN@{\Evaluate@#1}%
  \fi
 \else
  \DN@{\global\Value\csname\exstring@#1@C\endcsname}%
 \fi
 \next@}
\def\Evaluate@#1#2{\expandafter
 \ifx\csname\exstring@#1@C#2\endcsname\relax
  \DN@{\Err@{\string\Evaluate\string#1 can't be followed by \string#2}}%
 \else
  \DN@{\global\Value\csname\exstring@#1@C#2\endcsname}%
 \fi\next@}
\def\pre#1{\expandafter
 \ifx\csname\exstring@#1@P\endcsname\relax
  \expandafter\ifx\csname\exstring@#1@P1\endcsname\relax
   \DN@{\Err@{\noexpand\pre can't be used with \string#1}}%
  \else
   \DN@{\pre@#1}%
  \fi
 \else
  \DN@{{\csname\exstring@#1@P\endcsname}}%
 \fi
 \next@}
\def\pre@#1#2{\expandafter
 \ifx\csname\exstring@#1@P#2\endcsname\relax
  \DN@{\Err@{\string\pre\string#1 can't be followed by \string#2}}%
 \else
  \DN@{{\csname\exstring@#1@P#2\endcsname}}%
 \fi
 \next@}
\def\post#1{\expandafter
 \ifx\csname\exstring@#1@Q\endcsname\relax
  \expandafter\ifx\csname\exstring@#1@Q1\endcsname\relax
   \DN@{\Err@{\noexpand\post can't be used with \string#1}}%
  \else
   \DN@{\post@#1}%
  \fi
 \else
  \DN@{{\csname\exstring@#1@Q\endcsname}}%
 \fi
 \next@}
\def\post@#1#2{\expandafter
 \ifx\csname\exstring@#1@Q#2\endcsname\relax
  \DN@{\Err@{\string\post\string#1 can't be followed by \string#2}}%
 \else
  \DN@{{\csname\exstring@#1@Q#2\endcsname}}%
 \fi
 \next@}
\def\style#1{\expandafter
 \ifx\csname\exstring@#1@S\endcsname\relax
  \expandafter\ifx\csname\exstring@#1@S1\endcsname\relax
   \DN@{\Err@{\noexpand\style can't be used with \string#1}}%
  \else
   \DN@{\style@#1}%
  \fi
 \else
  \DN@{\csname\exstring@#1@S\endcsname}%
 \fi
 \next@}
\def\style@#1#2{\expandafter
 \ifx\csname\exstring@#1@S#2\endcsname\relax
  \DN@{\Err@{\string\style\string#1 can't be followed by \string#2}}%
 \else
  \DN@{\csname\exstring@#1@S#2\endcsname}%
 \fi
 \next@}
\def\fontstyle#1{\expandafter
 \ifx\csname\exstring@#1@F\endcsname\relax
  \expandafter\ifx\csname\exstring@#1@F1\endcsname\relax
   \DN@{\Err@{\noexpand\fontstyle can't be used with \string#1}}%
  \else
   \DN@{\fontstyle@#1}%
  \fi
 \else
  \DN@##1{{\csname\exstring@#1@F\endcsname##1}}%
 \fi
 \next@}
\def\fontstyle@#1#2{\expandafter
 \ifx\csname\exstring@#1@F#2\endcsname\relax
  \DN@{\Err@{\string\fontstyle\string#1 can't be followed by \string#2}}%
 \else
  \DN@##1{{\csname\exstring@#1@F#2\endcsname##1}}%
 \fi
 \next@}
\def\Reset#1{\expandafter
 \ifx\csname\exstring@#1@C\endcsname\relax
  \expandafter\ifx\csname\exstring@#1@C1\endcsname\relax
   \DN@{\Err@{\noexpand\Reset can't be used with \string#1}}%
  \else
   \DN@{\Reset@#1}%
  \fi
 \else
  \DN@##1{\count@##1\relax\ifx#1\page\else\advance\count@\m@ne\fi
   \global\csname\exstring@#1@C\endcsname\count@}%
 \fi
 \next@}
\def\Reset@#1#2{\expandafter
 \ifx\csname\exstring@#1@C#2\endcsname\relax
  \DN@{\Err@{\string\Reset\string#1 can't be followed by \string#2}}%
 \else
  \DN@##1{\count@##1\relax\advance\count@\m@ne
   \global\csname\exstring@#1@C#2\endcsname\count@}%
 \fi
 \next@}
\def\Offset#1{\expandafter
 \ifx\csname\exstring@#1@C\endcsname\relax
  \expandafter\ifx\csname\exstring@#1@C1\endcsname\relax
   \DN@{\Err@{\noexpand\Offset can't be used with \string#1}}%
  \else
   \DN@{\Offset@#1}%
  \fi
 \else
  \DN@##1{\count@##1\relax\advance\count@\m@ne\global\advance
   \csname\exstring@#1@C\endcsname\count@}%
 \fi
 \next@}
\def\Offset@#1#2{\expandafter
 \ifx\csname\exstring@#1@C#2\endcsname\relax
  \DN@{\Err@{\string\Offset\string#1 can't be followed by \string#2}}%
 \else
  \DN@##1{\count@##1\relax\advance\count@\m@ne
   \global\advance\csname\exstring@#1@C#2\endcsname\count@}%
 \fi
 \next@}
\def\getR@#1#2{\def\nextiv@{\let\nextiii@}\expandafter\nextiv@
 \csname\exstring@#1@R#2\endcsname}
\def\letR@#1#2#3{\expandafter\let\csname#1@#3#2\endcsname\Next@}
\def\letR@@#1#2{\expandafter\let\csname\exstring@#1@#2\endcsname\Next@}
\def\newpre#1{\expandafter
 \ifx\csname\exstring@#1@P\endcsname\relax
  \expandafter\ifx\csname\exstring@#1@P1\endcsname\relax
   \DN@{\Err@{\noexpand\newpre can't be used with \string#1}}%
  \else
   \DN@{\newpre@#1}%
  \fi
 \else
  \DN@{%
   \DNii@{%
    \endgroup
    \expandafter\let\csname\exstring@#1@P\endcsname\Next@
    \expandafter\ifx\csname\exstring@#1@R\endcsname\relax\else
    \getR@#1{}\expandafter\letR@\nextiii@ P\fi
    }%
   \begingroup\noexpands@\afterassignment\nextii@\xdef\Next@}%
 \fi
 \next@}
\def\newpre@#1#2{\expandafter
 \ifx\csname\exstring@#1@P#2\endcsname\relax
  \DN@{\Err@{\string\newpre\string#1 can't be followed by \string#2}}%
 \else
  \DN@{%
   \DNii@{%
    \endgroup
    \expandafter\let\csname\exstring@#1@P#2\endcsname\Next@
    \expandafter\ifx\csname\exstring@#1@R#2\endcsname\relax\else
    \getR@#1{#2}\expandafter\letR@@\nextiii@ P\fi
    }%
   \begingroup\noexpands@\afterassignment\nextii@\xdef\Next@}%
 \fi
 \next@}
\def\newpost#1{\expandafter
 \ifx\csname\exstring@#1@Q\endcsname\relax
  \expandafter\ifx\csname\exstring@#1@Q1\endcsname\relax
   \DN@{\Err@{\noexpand\newpost can't be used with \string#1}}%
  \else
   \DN@{\newpost@#1}%
  \fi
 \else
  \DN@{%
   \DNii@{%
    \endgroup
    \expandafter\let\csname\exstring@#1@Q\endcsname\Next@
    \expandafter\ifx\csname\exstring@#1@R\endcsname\relax\else
    \getR@#1{}\expandafter\letR@\nextiii@ Q\fi
    }%
   \begingroup\noexpands@\afterassignment\nextii@\xdef\Next@}%
 \fi
 \next@}
\def\newpost@#1#2{\expandafter
 \ifx\csname\exstring@#1@Q#2\endcsname\relax
  \DN@{\Err@{\string\newpost\string#1 can't be followed by \string#2}}%
 \else
  \DN@{%
   \DNii@{%
    \endgroup
    \expandafter\let\csname\exstring@#1@Q#2\endcsname\Next@
    \expandafter\ifx\csname\exstring@#1@R#2\endcsname\relax\else
    \getR@#1{#2}\expandafter\letR@@\nextiii@ Q\fi
    }%
   \begingroup\noexpands@\afterassignment\nextii@\xdef\Next@}%
 \fi
 \next@}
\def\newstyle#1{\expandafter
 \ifx\csname\exstring@#1@S\endcsname\relax
  \expandafter\ifx\csname\exstring@#1@S1\endcsname\relax
   \DN@{\Err@{\noexpand\newstyle can't be used
    with \string#1}}%
  \else
   \DN@{\newstyle@#1}%
  \fi
 \else
  \DN@{%
   \DNii@{%
    \expandafter\let\csname\exstring@#1@S\endcsname\Next@
    \expandafter\ifx\csname\exstring@#1@R\endcsname\relax\else
    \getR@#1{}\expandafter\letR@\nextiii@ S\fi
    }%
   \afterassignment\nextii@\gdef\Next@}%
 \fi
 \next@}
\def\newstyle@#1#2{\expandafter
 \ifx\csname\exstring@#1@S#2\endcsname\relax
  \DN@{\Err@{\string\newstyle\string#1 can't be followed by
   \string#2}}%
 \else
  \DN@{%
   \DNii@{%
    \expandafter\let\csname\exstring@#1@S#2\endcsname\Next@
    \expandafter\ifx\csname\exstring@#1@R#2\endcsname\relax\else
    \getR@#1{#2}\expandafter\letR@@\nextiii@ S\fi
    }%
   \afterassignment\nextii@\gdef\Next@}%
 \fi
 \next@}
\def\newnumstyle#1{\expandafter
 \ifx\csname\exstring@#1@N\endcsname\relax
  \expandafter\ifx\csname\exstring@#1@N1\endcsname\relax
   \DN@{\Err@{\noexpand\newnumstyle can't be used with
    \string#1}}%
  \else
   \DN@{\newnumstyle@#1}%
  \fi
 \else
  \DN@##1{%
   \gdef\Next@{##1}%
    \expandafter\let\csname\exstring@#1@N\endcsname\Next@
    \expandafter\ifx\csname\exstring@#1@R\endcsname\relax\else
    \getR@#1{}\expandafter\letR@\nextiii@ N\fi
    }%
 \fi
 \next@}
\def\newnumstyle@#1#2{\expandafter
 \ifx\csname\exstring@#1@N#2\endcsname\relax
  \DN@{\Err@{\string\newnumstyle\string#1 can't be followed by
   \string#2}}%
 \else
  \DN@##1{%
   \gdef\Next@{##1}%
    \expandafter\let\csname\exstring@#1@N#2\endcsname\Next@
    \expandafter\ifx\csname\exstring@#1@R#2\endcsname\relax\else
    \getR@#1{#2}\expandafter\letR@@\nextiii@ N\fi
    }%
  \fi
 \next@}
\def\newfontstyle#1{\expandafter
 \ifx\csname\exstring@#1@F\endcsname\relax
  \expandafter\ifx\csname\exstring@#1@F1\endcsname\relax
   \DN@{\Err@{\noexpand\newfontstyle can't be used with
    \string#1}}%
  \else
   \DN@{\newfontstyle@#1}%
  \fi
 \else
  \DN@##1{%
   \gdef\Next@{##1}%
    \expandafter\let\csname\exstring@#1@F\endcsname\Next@
    \expandafter\ifx\csname\exstring@#1@R\endcsname\relax\else
    \getR@#1{}\expandafter\letR@\nextiii@ F\fi
    }%
 \fi
 \next@}
\def\newfontstyle@#1#2{\expandafter
 \ifx\csname\exstring@#1@F#2\endcsname\relax
  \DN@{\Err@{\string\newfontstyle\string#1 can't be followed by
   \string#2}}%
 \else
  \DN@##1{%
   \gdef\Next@{##1}%
    \expandafter\let\csname\exstring@#1@F#2\endcsname\Next@
    \expandafter\ifx\csname\exstring@#1@R#2\endcsname\relax\else
    \getR@#1{#2}\expandafter\letR@@\nextiii@ F\fi
    }%
 \fi
 \next@}
\def\word#1{\expandafter
 \ifx\csname\exstring@#1@W\endcsname\relax
  \expandafter\ifx\csname\exstring@#1@W1\endcsname\relax
   \DN@{\Err@{\noexpand\word can't be used with \string#1}}%
  \else
   \DN@{\word@#1}%
  \fi
 \else
  \DN@{{\csname\exstring@#1@W\endcsname}}%
 \fi
 \next@}
\def\word@#1#2{\expandafter
 \ifx\csname\exstring@#1@W#2\endcsname\relax
  \DN@{\Err@{\string\word\noexpand#1can't be followed by \string#2}}%
 \else
  \DN@{{\csname\exstring@#1@W#2\endcsname}}%
 \fi
 \next@}
\def\newword#1{\expandafter
 \ifx\csname\exstring@#1@W\endcsname\relax
  \expandafter\ifx\csname\exstring@#1@W1\endcsname\relax
   \DN@{\Err@{\noexpand\newword can't be used  with \string#1}}%
  \else
   \DN@{\newword@#1}%
  \fi
 \else
  \DN@{%
   \DNii@{%
    \expandafter\let\csname\exstring@#1@W\endcsname\Next@
    \expandafter\ifx\csname\exstring@#1@R\endcsname\relax\else
     \getR@#1{}\expandafter\letR@\nextiii@ W\fi
    }%
   \afterassignment\nextii@\gdef\Next@}%
 \fi
 \next@}
\def\newword@#1#2{\expandafter
 \ifx\csname\exstring@#1@W#2\endcsname\relax
  \DN@{\Err@{\string\newword\noexpand#1can't be followed by \string#2}}%
 \else
  \DN@{%
   \DNii@{%
    \expandafter\let\csname\exstring@#1@W#2\endcsname\Next@
    \expandafter\ifx\csname\exstring@#1@R#2\endcsname\relax\else
     \getR@#1{#2}\expandafter\letR@@\nextiii@ W\fi
    }%
   \afterassignment\nextii@\gdef\Next@}%
 \fi
 \next@}
\newif\iffn@
\newcount\footmark@C
\footmark@C\z@
\def\footmark@S#1{$^{#1}$}
\let\footmark@N\arabic
\def\footmark@F{\rm}
\def\foottext@S#1{$^{#1}$}
\def\foottext@F{\rm}
\let\modifyfootnote@\relax
\def\modifyfootnote#1{\def\modifyfootnote@{#1}}
\def\vfootnote@#1{\insert\footins
 \bgroup
 \floatingpenalty\@MM\interlinepenalty\interfootnotelinepenalty
 \leftskip\z@\rightskip\z@\spaceskip\z@\xspaceskip\z@
 \rm\splittopskip\ht\strutbox\splitmaxdepth\dp\strutbox
 \locallabel@\noindent@@{\foottext@F#1}\modifyfootnote@
 \footstrut\FN@\fo@t}
\def\fo@t{\ifcat\bgroup\noexpand\next\expandafter\f@@t\else
 \expandafter\f@t\fi}
\def\f@t#1{#1\@foot}
\def\f@@t{\bgroup\aftergroup\@foot\afterassignment\FNSSP@\let\next@}
\def\@foot{\unskip\lower\dp\strutbox\vbox to\dp\strutbox{}\egroup
 \iffn@\expandafter\fn@false\else
 \expandafter\postvanish@\fi}
\newif\ifplainfn@
\plainfn@true
\def\fancyfootnotes{\plainfn@false}
\newcount\fancyfootmarkcount@
\fancyfootmarkcount@\z@
\newcount\lastfnpage@
\lastfnpage@-\@M
\let\justfootmarklist@\empty
\def\footmark{\let\@sf\empty
 \ifhmode\edef\@sf{\spacefactor\the\spacefactor}\/\fi
 \DN@{\ifx"\next\expandafter\nextii@\else\expandafter\footmark@\fi}%
 \DNii@"##1"{%
  \iffirstchoice@
   {\let\style\footmark@S\let\numstyle\footmark@N
   \footmark@F##1%
   \noexpands@
   \let\style\foottext@S
   \Qlabel@{##1}%
   }%
   \iffn@\else
    {\noexpands@
    \xdef\Next@{{\Thelabel@}{\Thelabel@@}{\Thelabel@@@}{\Thelabel@@@@}}%
    }%
    \expandafter\rightappend@\Next@\to\justfootmarklist@
   \fi
  \fi
  \@sf\relax}%
 \FN@\next@}
\def\footmark@{%
 \iffirstchoice@
  \global\advance\footmark@C\@ne
  \ifplainfn@
   \xdef\adjustedfootmark@{\number\footmark@C}%
  \else
   {\let\\\or\xdef\Next@{\ifcase\number\footmark@C\fnpages@\else
     -\@M\fi}}%
   \ifnum\Next@=-\@M
    \xdef\adjustedfootmark@{\number\footmark@C}%
   \else
    \ifnum\Next@=\lastfnpage@
     \global\advance\fancyfootmarkcount@\@ne
    \else
     \global\fancyfootmarkcount@\@ne
     \global\lastfnpage@\Next@
    \fi
    \xdef\adjustedfootmark@{\number\fancyfootmarkcount@}%
   \fi
  \fi
  {\noexpands@
  \xdef\Thelabel@@@{\adjustedfootmark@}%
  \xdefThelabel@\footmark@N
  \xdef\Thelabel@@@@{\Thelabel@}%
  \xdefThelabel@@\foottext@S
  }%
  \iffn@\else
   {\noexpands@
   \xdef\Next@{{\Thelabel@}{\Thelabel@@}{\Thelabel@@@}{\Thelabel@@@@}}%
   }%
   \expandafter\rightappend@\Next@\to\justfootmarklist@
  \fi
  \ifplainfn@
  \else
   \edef\next@{\write\laxwrite@{F\noexpand\the\pageno}}\next@
  \fi
 \fi
 \footmark@S{\footmark@N{\adjustedfootmark@}}%
 \@sf\relax}
\def\foottext{\prevanish@
 \ifx\justfootmarklist@\empty
  \Err@{There is no \noexpand\footmark for this \string\foottext}\fi
 \DN@\\##1##2\next@{\DN@{##1}\gdef\justfootmarklist@{##2}}%
 \expandafter\next@\justfootmarklist@\next@
 \expandafter\foottext@\next@}
\def\foottext@#1#2#3#4{{\noexpands@
  \xdef\Thelabel@{#1}\xdef\Thelabel@@{#2}%
  \xdef\Thelabel@@@{#3}\xdef\Thelabel@@@@{#4}}%
  \vfootnote@{\thelabel@@}}
\rightadd@\foottext\to\vanishlist@
\def\footnote{\fn@true
 \let\@sf\empty
 \ifhmode\edef\@sf{\spacefactor\the\spacefactor}\/\fi
 \DN@{\ifx"\next\expandafter\nextii@\else\expandafter\nextiii@\fi}%
 \DNii@"##1"{\footmark"##1"\vfootnote@{\let\style\foottext@S
  \let\numstyle\footmark@N##1}}%
 \def\nextiii@{\footmark\vfootnote@{\foottext@S{\footmark@N
  {\adjustedfootmark@}}}}%
 \FN@\next@}
\newdimen\litindent
\litindent20\p@
\newbox\litbox@
\newbox\Litbox@
\newcount\interlitpenalty@
\interlitpenalty@\@M
\newcount\litlines@
{\obeyspaces\gdef\defspace@{\def {\allowbreak\hskip.5emminus.15em}}}
{\obeylines\gdef\letM@{\let^^M\CtrlM@}}
\def\CtrlM@{\egroup
 \ifcase\litlines@\advance\litlines@\@ne\or
 \box\litbox@\advance\litlines@\@ne\else
 \penalty\interlitpenalty@\box\litbox@\fi
 \Lit@}
\def\Lit@{\setbox\litbox@\hbox\bgroup\litdefs@\hskip\litindent}
\newcount\littab@
\littab@8
\def\littab#1{\littab@#1\relax}
{\catcode`\^^I=\active\gdef\letTAB@{\let^^I\TAB@}}
\def\TAB@{\egroup
 \dimen@\wd\litbox@
 \advance\dimen@-\litindent
 \setboxz@h{\tt0}%
 \dimen@ii\littab@\wdz@
 \divide\dimen@\dimen@ii
 \multiply\dimen@\dimen@ii
 \advance\dimen@\littab@\wdz@
 \advance\dimen@\litindent
 \setbox\litbox@\hbox\bgroup\litdefs@\hbox to\dimen@{\unhbox\litbox@\hfil}}
{\catcode`\`=\active\gdef`{\relax\lq}}
\let\litbs@\relax
\let\litbs@@\relax
\def\litbackslash#1{%
 \edef\litbs@{\catcode`\string#1=\z@
 \def\noexpand\litbs@@{\def\expandafter\noexpand\csname\string#1\endcsname
  {\char`\string#1}}}}
\def\litcodes@{\catcode`\\=12
 \catcode`\{=12 \catcode`\}=12
 \catcode`\$=12 \catcode`\&=12
 \catcode`\#=12
 \catcode`\^=12 \catcode`\_=12
 \catcode`\@=12 \catcode`\~=12 \catcode`\"=12
 \catcode`\;=12 \catcode`\:=12 \catcode`\!=12 \catcode`\?=12
 \catcode`\%=12 \litbs@\catcode`\`=\active\obeyspaces\defspace@}
\def\activate@#1#2{{\lccode`\~=`#2%
 \lowercase{%
  \if0#1%
  \gdef\Next@{\def~{\egroup\endgroup\bigskip\vskip-\parskip
   \def\next@{\noindent@@\FN@\pretendspace@}\FNSS@\next@}}\else
  \gdef\Next@{\def~{\egroup\egroup\endgroup}}\fi
  }%
 }}
\def\litdefs@{\let\0\empty\let\1\litdelim@\def\ {\char32 }\litbs@@}%
\def\litdelimiter#1{%
 \edef\litdelim@{\char`#1}%
 \def\lit#1{\leavevmode\begingroup\litcodes@\litdefs@
  \tt\hyphenchar\tentt\m@ne\lit@}%
 \def\lit@##1#1{##1\endgroup\null}%
 \def\Lit#1{\ifhmode$$\abovedisplayskip\bigskipamount
  \abovedisplayshortskip\bigskipamount
  \belowdisplayskip\z@\belowdisplayshortskip\z@
  \postdisplaypenalty\@M
  $$\vskip-\baselineskip\else\bigskip\fi
  \begingroup\litlines@\z@
  \catcode`#1=\active\activate@0#1\Next@
  \def\displaybreak{\egroup\break\litlines@\z@\Lit@}%
  \def\allowdisplaybreak{\egroup\allowbreak\litlines@\z@\Lit@}%
  \def\allowdisplaybreaks{\egroup\allowbreak\interlitpenalty@\z@
   \litlines@\z@\Lit@}%
  \litcodes@\tt\catcode`\^^I=\active\letTAB@
  \obeylines\letM@\Lit@}%
 \def\Litbox##1=#1{\begingroup\ifodd##1\relax\aftergroup\global\fi
  \aftergroup\setbox\aftergroup##1\aftergroup\box\aftergroup\Litbox@
  \def\allowdisplaybreak{\egroup\allowbreak\litlines@\z@\Lit@}%
  \def\allowdisplaybreaks{\egroup\allowbreak\interlitpenalty@\z@
   \litlines@\z@\Lit@}%
  \catcode`#1=\active\activate@1#1\Next@
  \litcodes@\tt\catcode`\^^I=\active\letTAB@
  \obeylines\letM@\global\setbox\Litbox@\vbox\bgroup\litindent\z@%
  \litlines@\z@\Lit@}%
}
\newbox\titlebox@
\setbox\titlebox@\vbox{}
\rightadd@\title\to\overlonglist@
\def\title{\begingroup\Let@
 \global\setbox\titlebox@\vbox\bgroup\tabskip\hss@
 \halign to\hsize\bgroup\bf\hfil\ignorespaces##\unskip\hfil\cr}
\def\endtitle{\crcr\egroup\egroup\endgroup\overlong@false}
\newbox\authorbox@
\rightadd@\author\to\overlonglist@
\def\author{\begingroup\Let@
 \global\setbox\authorbox@\vbox\bgroup\tabskip\hss@
 \halign to\hsize\bgroup\rm\hfil\ignorespaces##\unskip\hfil\cr}
\def\endauthor{\crcr\egroup\egroup\endgroup\overlong@false}
\newbox\affilbox@
\def\affil{\begingroup\Let@
 \global\setbox\affilbox@\vbox\bgroup\tabskip\hss@
 \halign to\hsize\bgroup\rm\hfil\ignorespaces##\unskip\hfil\cr}%
\def\endaffil{\crcr\egroup\egroup\endgroup\overlong@false}
\let\date@\relax
\def\date#1{\gdef\date@{\ignorespaces#1\unskip}}
\def\today{\ifcase\month\or January\or February\or March\or April\or May\or
 June\or July\or August\or September\or October\or November\or December\fi
 \space\number\day, \number\year}
\def\maketitle{\hrule\height\z@\vskip-\topskip
 \vskip24\p@ plus12\p@ minus12\p@
 \unvbox\titlebox@
 \ifvoid\authorbox@\else\vskip12\p@ plus6\p@ minus3\p@\unvbox\authorbox@\fi
 \ifvoid\affilbox@\else\vskip10\p@ plus5\p@ minus2\p@\unvbox\affilbox@\fi
 \ifx\date@\relax\else\vskip6\p@ plus2\p@ minus\p@\centerline{\rm\date@}\fi
 \vskip18\p@ plus12\p@ minus6\p@}
\def\cite{%
 \DNii@(##1)##2{{\rm[}{##2}, {##1\/}{\rm]}}%
 \def\nextiii@##1{{\rm[}{##1\/}{\rm]}}%
 \DN@{\ifx\next(\expandafter\nextii@\else\expandafter\nextiii@\fi}%
 \FN@\next@}
\def\makebib@W{Bibliography}
\def\makebib{\begingroup\rm\bigbreak\centerline{\smc\makebib@W}%
 \nobreak\medskip
 \sfcode`\.=\@m\everypar{}\parindent\z@
 \def\nopunct{\nopunct@true}\def\nospace{\nospace@true}%
 \nopunct@false\nospace@false
 \def\lkerns@{\null\kern\m@ne sp\kern\@ne sp}%
 \def\nkerns@{\null\kern-\tw@ sp\kern\tw@ sp}%
}
\let\endmakebib\endgroup
\newif\ifnoprepunct@
\newif\ifnoprespace@
\newif\ifnoquotes@
\def\noprepunct{\noprepunct@true}
\def\noprespace{\noprespace@true}
\def\noquotes{\noquotes@true}
\newbox\nobox@
\newbox\keybox@
\newbox\bybox@
\newbox\paperbox@
\newbox\paperinfobox@
\newbox\jourbox@
\newbox\volbox@
\newbox\issuebox@
\newbox\yrbox@
\newbox\pgbox@
\newbox\ppbox@
\newbox\bookbox@
\newbox\inbookbox@
\newbox\bookinfobox@
\newbox\publbox@
\newbox\publaddrbox@
\newbox\edbox@
\newbox\edsbox@
\newbox\langbox@
\newbox\translbox@
\newbox\finalinfobox@
\def\setbibinfo@#1{\edef\next@{\ifnopunct@1\else0\fi
 \ifnospace@1\else0\fi\ifnoprepunct@1\else0\fi\ifnoprespace@1\else0\fi
 \ifnoquotes@1\else0\fi}%
 \DNii@{00000}%
 \ifx\next@\nextii@\else\xdef\bibinfo@{\bibinfo@\the#1,\next@}%
 \fi}
\def\getbibinfo@#1{\ifx\bibinfo@\empty
 \let\next@0\let\nextii@0\let\nextiii@0\let\nextiv@0\let\nextv@0\else
 \edef\next@{\def
  \noexpand\next@####1\the#1,####2####3####4####5####6####7\noexpand\next@
  {\let\noexpand\next@####2\let\noexpand\nextii@####3%
  \let\noexpand\nextiii@####4\let\noexpand\nextiv@####5%
  \let\noexpand\nextv@####6}%
  \noexpand\next@\bibinfo@\the#1,00000\noexpand\next@}\next@
 \fi}
\newif\ifbookinquotes@
\def\bookinquotes{\bookinquotes@true}
\newif\ifpaperinquotes@
\def\paperinquotes{\paperinquotes@true}
\newif\ifininbook@
\def\ininbook{\ininbook@true}
\newif\ifopenquotes@
\def\closequotes@{\ifopenquotes@''\openquotes@false\fi}
\newif\ifbeginbib@
\newif\ifendbib@
\newif\ifprevjour@
\newif\ifprevbook@
\newdimen\bibindent@
\bibindent@20\p@
\def\bib{\global\let\bibinfo@\empty\global\let\translinfo@\relax\beginbib@true
 \begingroup\noindent@
 \hangindent\bibindent@\hangafter\@ne\bib@}
\def\v@id#1{\setbox#1\box\voidb@x}
\def\bib@{\v@id\nobox@\v@id\keybox@\v@id\bybox@\v@id\paperbox@
 \v@id\paperinfobox@\v@id\jourbox@\v@id\volbox@\v@id\issuebox@
 \v@id\yrbox@\v@id\pgbox@\v@id\ppbox@\v@id\bookbox@\v@id\inbookbox@
 \v@id\bookinfobox@\v@id\publbox@\v@id\publaddrbox@\v@id\edbox@
 \v@id\edsbox@\v@id\langbox@\v@id\translbox@\v@id\finalinfobox@
 \bgroup}
\def\Setnonemptybox@#1#2{\unskip\setbibinfo@#1\egroup#2%
 \def\aftergroup@{\ifdim\wd#1=\z@\setbox#1\box\voidb@x\fi}%
 \setbox#1\vbox\bgroup\aftergroup\aftergroup@\hsize\maxdimen\leftskip\z@
 \rightskip\z@\hbadness\@M\hfuzz\maxdimen\noindent}
\def\setnonemptybox@#1{\Setnonemptybox@#1\relax}
\def\no{\setnonemptybox@\nobox@}
\def\key{\setnonemptybox@\keybox@\bf}
\def\by{\setnonemptybox@\bybox@}
\def\bysame{\setnonemptybox@\bybox@\leaders\hrule\hskip3em\null}
\def\paper{\setnonemptybox@\paperbox@
 \ifpaperinquotes@\getbibinfo@\paperbox@
 \if\nextv@1\else``\fi\else\it\fi}
\def\paperinfo{\setnonemptybox@\paperinfobox@}
\def\jour{\Setnonemptybox@\jourbox@\prevjour@true}
\def\vol{\setnonemptybox@\volbox@\bf}
\def\issue{\setnonemptybox@\issuebox@}
\def\yr{\setnonemptybox@\yrbox@}

\def\pg{\setnonemptybox@\pgbox@}
\def\pp{\setnonemptybox@\ppbox@}
\def\book{\Setnonemptybox@\bookbox@\prevbook@true
 \ifbookinquotes@\getbibinfo@\bookbox@
 \if\nextv@1\else``\fi\else\it\fi}
\def\inbook{\Setnonemptybox@\inbookbox@\prevbook@true
 \ifininbook@ in \fi\ifbookinquotes@\getbibinfo@\inbookbox@
 \if\nextv@1\else``\fi\fi}
\def\bookinfo{\setnonemptybox@\bookinfobox@}
\def\publ{\setnonemptybox@\publbox@}
\def\publaddr{\setnonemptybox@\publaddrbox@}
\def\ed{\setnonemptybox@\edbox@}
\def\eds{\setnonemptybox@\edsbox@}
\def\lang{\setnonemptybox@\langbox@}
\def\finalinfo{\setnonemptybox@\finalinfobox@}
\def\setboxzl@{\setbox\z@\lastbox}
\def\getbox@#1{\setbox\z@\vbox{\vskip-\@M\p@
 \unvbox#1%
 \setboxzl@
 \global\setbox\@ne\hbox{\unhbox\z@\unskip\unskip\unpenalty}%
 \ifdim\lastskip=-\@M\p@\else
 \loop\ifdim\lastskip=-\@M\p@
 \else\unskip\unpenalty\setboxzl@
 \global\setbox\@ne\hbox{\unhbox\z@\unhbox\@ne}%
 \repeat\fi}%
 \unhbox\@ne}
\def\adjustpunct@#1{\count@\lastkern
 \ifnum\count@=\z@#1\closequotes@\else
 \ifnum\count@>\tw@#1\closequotes@\else
 \ifnum\count@<-\tw@#1\closequotes@\else
  \unkern\unkern\setboxzl@
  \skip@\lastskip\unskip
  \count@@\lastpenalty\unpenalty
  \ifnum\count@=\tw@\unskip\setboxzl@\fi
  \ifdim\skip@=\z@\else\hskip\skip@\fi
  #1\closequotes@
  \ifnum\count@=\tw@\null\hfill\fi
  \penalty\count@@
 \fi\fi\fi}
\def\prepunct@#1#2{\getbibinfo@#2%
 \ifnopunct@
 \else
  \if\nextiii@0\adjustpunct@#1\fi
 \fi
 \closequotes@
 \ifnospace@
 \else
  \if\nextiv@0\space\else\fi
 \fi
 \nopunct@false\nospace@false
 \if\next@1\nopunct@true\fi
 \if\nextii@1\nospace@true\fi}
\def\ppunbox@#1#2{\prepunct@{#1}#2%
 \getbox@#2}
\let\semicolon@;
\def\endbib@{%
 \ifbeginbib@
  \ifvoid\nobox@
   \ifvoid\keybox@\else\hbox to\bibindent@{[\getbox@\keybox@]\hss}\fi
  \else\hbox to\bibindent@{\hss\getbox@\nobox@. }\fi
  \ifvoid\bybox@\else\getbox@\bybox@\fi
 \else
  \nopunct@true
  \ifvoid\bybox@\else\ppunbox@\relax\bybox@\fi
 \fi
 \ifvoid\translbox@\else\ppunbox@,\translbox@\fi
 \ifvoid\paperbox@\else\ppunbox@,\paperbox@\ifpaperinquotes@
  \if\nextv@1\else\openquotes@true\fi\fi
 \fi
 \ifvoid\paperinfobox@\else\ppunbox@,\paperinfobox@\fi
 \test@false
 \ifvoid\jourbox@\else\test@true\ppunbox@,\jourbox@\fi
 \ifprevjour@\test@true\fi
 \iftest@
  \ifvoid\volbox@\else\ppunbox@\relax\volbox@\fi
  \ifvoid\issuebox@
   \else\prepunct@\relax\issuebox@ no.~\getbox@\issuebox@\fi
  \ifvoid\yrbox@\else\prepunct@\relax\yrbox@(\getbox@\yrbox@)\fi
  \ifvoid\ppbox@\else\ppunbox@,\ppbox@\fi
  \ifvoid\pgbox@\else\prepunct@,\pgbox@ p.~\getbox@\pgbox@\fi
 \fi
 \test@false
 \ifvoid\bookbox@\else\test@true\ppunbox@,\bookbox@\ifbookinquotes@
  \if\nextv@1\else\openquotes@true\fi\fi\fi
 \ifvoid\inbookbox@\else\test@true\ppunbox@,\inbookbox@\ifbookinquotes@
  \if\nextv@1\else\openquotes@true\fi\fi\fi
 \ifprevbook@\test@true\fi
 \iftest@
  \ifvoid\edbox@\else\prepunct@\relax\edbox@(\getbox@\edbox@, ed.)\fi
  \ifvoid\edsbox@\else\prepunct@\relax\edsbox@(\getbox@\edsbox@, eds.)\fi
  \ifvoid\bookinfobox@\else\ppunbox@,\bookinfobox@\fi
  \ifvoid\publbox@\else\ppunbox@,\publbox@\fi
  \ifvoid\publaddrbox@\else\ppunbox@,\publaddrbox@\fi
  \ifvoid\yrbox@\else\ppunbox@,\yrbox@\fi
  \ifvoid\ppbox@\else\prepunct@,\ppbox@ pp.~\getbox@\ppbox@\fi
  \ifvoid\pgbox@\else\prepunct@,\pgbox@ p.~\getbox@\pgbox@\fi
 \fi
 \ifvoid\finalinfobox@
  \ifendbib@
   \ifnopunct@\else.\closequotes@\fi
  \else
  \ifvoid\langbox@\else\space(\getbox@\langbox@)\fi
   \/\semicolon@\closequotes@
  \fi
 \else
  \ifendbib@
   \ppunbox@{.\spacefactor3000\relax}\finalinfobox@
    \ifnopunct@\else.\fi
  \else
   \ppunbox@,\finalinfobox@\/\semicolon@\fi
 \fi
 \ifvoid\langbox@\else\space(\getbox@\langbox@)\fi
}
\def\endbib{\unskip\egroup\endbib@true\endbib@\par\endgroup}
\def\morebib{\unskip\egroup
 \endbib@false\endbib@
 \global\let\bibinfo@\empty\beginbib@false
 \bib@}
\def\anotherbib{\unskip\egroup
 \endbib@false\endbib@
 \global\let\bibinfo@\empty\beginbib@false
 \prevjour@false\prevbook@false\bib@}
\def\transl{\unskip
 \xdef\translinfo@{\the\translbox@,\ifnopunct@1\else0\fi
 \ifnospace@1\else0\fi\ifnoprepunct@1\else0\fi\ifnoprespace@1\else0\fi0}%
 \egroup\endbib@false\endbib@
 \global\let\bibinfo@\translinfo@\beginbib@false
 \bib@
 \egroup
 \def\aftergroup@{\ifdim\wd\translbox@=\z@\setbox\translbox@\box\voidb@x\fi}%
 \setbox\translbox@\vbox\bgroup\aftergroup\aftergroup@
 \hsize\maxdimen\leftskip\z@\rightskip\z@\hbadness\@M\hfuzz\maxdimen
 \noindent}
\newwrite\auxwrite@
\newread\bbl@
\def\UseBibTeX{\immediate\openout\auxwrite@=\jobname.aux
 \let\cite\BTcite@
 \def\nocite##1{\immediate\write\auxwrite@{\string\citation{##1}}}%
 \def\bibliographystyle##1{\immediate\write\auxwrite@{\string
  \bibstyle{##1}}}%
 \def\bibliography@W{Bibliography}%
 \def\bibliography##1{\immediate\write\auxwrite@{\string\bibdata{##1}}%
  \immediate\openin\bbl@=\jobname.bbl
  \ifeof\bbl@
   \W@{No .bbl file}%
  \else
   \immediate\closein\bbl@
   \begingroup\input bibtex \input\jobname.bbl \endgroup
  \fi}%
 }
\def\BTcite@{%
 \DNii@(##1)##2{{\rm[}\BTcite@@##2,\BTcite@@{\rm, }{##1\/}{\rm]}%
  \immediate\write\auxwrite@{\string\citation{##2}}}%
 \def\nextiii@##1{{\rm[}\BTcite@@##1,\BTcite@@\/{\rm]}%
  \immediate\write\auxwrite@{\string\citation{##1}}}%
 \DN@{\ifx\next(\expandafter\nextii@\else\expandafter\nextiii@\fi}%
 \FN@\next@}%
\def\BTcite@@#1,{\BTcite@@@{#1}\FN@\BTcite@@@@}
\def\BTcite@@@@{\ifx\next\BTcite@@
 \expandafter\eat@\else{\rm, }\expandafter\BTcite@@\fi}
\catcode`\~=11
\def\BTcite@@@#1{\nolabel@\cite{#1}\relax
 \DNii@##1~##2\nextii@{##1}%
 \csL@{#1}\expandafter\nextii@\Next@\nextii@\fi}
\catcode`\~=\active

\def\beginthebibliography@#1{\rm\setboxz@h{#1\ }\bibindent@\wdz@
 \bigbreak\centerline{\smc\bibliography@W}\nobreak\medskip
 \sfcode`\.=\@m\everypar{}\parindent\z@}

\newif\iffigproofing@
\def\Figureproofing{\figproofing@true}
\def\noFigureproofing{\figproofing@false}
\newif\ifHby@
\def\Hbyw#1{\global\Hby@true\hbyw\vsize{#1}}
\def\hbyw#1#2{%
 \hbox{%
  \ifHby@
  \else
   \iffigproofing@
    \setbox\z@\vbox{\hrule\width5\p@}\ht\z@\z@
    \vbox to#1{\hrule\height5\p@\width.4\p@\vfil\hrule\height5\p@\width.4\p@}%
    \kern-.4\p@\rlap{\copy\z@}\raise#1\hbox{\rlap{\copy\z@}}%
   \fi
  \fi
  \vbox to#1{\hbox to#2{}\vfil}%
  \ifHby@
  \else
   \iffigproofing@
    \vbox to#1{\hrule\height5\p@\width.4\p@\vfil\hrule\height5\p@\width.4\p@}%
    \kern-.4\p@\llap{\copy\z@}\raise#1\hbox{\llap{\boxz@}}%
   \fi
  \fi}}
\newcount\island@C
\let\island@P\empty
\let\island@Q\empty
\def\island@S#1{#1\null.}
\let\island@N\arabic
\def\island@F{\rm}
\def\island@@@P{\csname\exxx@\islandtype@ @P\endcsname}
\def\island@@@Q{\csname\exxx@\islandtype@ @Q\endcsname}
\def\island@@@S{\csname\exxx@\islandtype@ @S\endcsname}
\def\island@@@N{\csname\exxx@\islandtype@ @N\endcsname}
\def\island@@@F{\csname\exxx@\islandtype@ @F\endcsname}
\def\island@@@C{\csname island@C\islandclass@\endcsname}
\newif\ifplace@
\newif\ifisland@
\def\island{%
 \ifplace@
  \DN@{\let\islandclass@\empty\def\islandtype@{\island}\FN@\island@}%
 \else
  \long\DN@##1\endisland{\Err@{\noexpand\island must be used after some
   type of \string\...place}}%
 \fi
 \next@}
\def\island@{\ifx\next\c\let\next@\island@c\else
 \DN@{\FN@\island@@}\fi\next@}
\def\island@@{\ifcat\bgroup\noexpand\next\let\next@\island@@@\else
 \DN@{\Err@{\noexpand\island must be followed by a {prefix} for
 \string\caption's}}\fi\next@}
\newbox\islandbox@
\newcount\captioncount@
\def\island@@@#1{\def\captionprefix@{#1}\captioncount@\z@
 \global\setbox\islandbox@\vbox\bgroup}
\def\island@c\c#1{%
 \ifplace@
 \DN@{\def\islandclass@{#1}%
  \expandafter\ifx\csname island@C#1\endcsname\relax
  \expandafter\newcount@\csname island@C#1\endcsname
   \global\csname island@C#1\endcsname\z@\fi
  \FNSS@\island@c@}%
 \else
 \DN@{\edef\next@{\long\def\noexpand\next@########1\expandafter\noexpand
  \csname end\exxx@\islandtype@\endcsname{\noexpand\Err@{\noexpand\noexpand
  \expandafter\noexpand
  \islandtype@ must be used after some type of \noexpand\string
   \noexpand\...place}}}\next@\next@}%
 \fi
 \next@}
\def\island@c@{%
 \ifcat\bgroup\noexpand\next
  \let\next@\island@c@@
 \else
  \DN@{\Err@{\noexpand\island\string\c{\expandafter\string\islandclass@} must
   be followed by a {prefix} for \string\caption's}}%
 \fi\next@}
\def\island@c@@#1{\def\captionprefix@{#1}%
 \captioncount@\z@\global\setbox\islandbox@\vbox\bgroup}
\rightadd@\caption\to\nofrillslist@
\newbox\captionbox@
\newbox\Captionbox@
\def\caption{%
 \ifnum\captioncount@=\z@
  \ifnopunct@
   \DN@{\egroup\nopunct@true}%
  \else
   \let\next@\egroup
  \fi
 \else
  \let\next@\relax
 \fi
 \next@
 \advance\captioncount@\@ne
 \FN@\caption@}
\def\caption@{\ifx\next"\expandafter\caption@q\else\expandafter\caption@@\fi}
\def\caption@q"#1"{\quoted@true
 {\noexpands@
 \let\pre\island@@@P\let\post\island@@@Q
 \let\style\island@@@S\let\numstyle\island@@@N
 \Qlabel@{#1}\let\style\relax\xdef\Qlabel@@@@{#1}}%
 \finishcaption@}
\def\caption@@{\quoted@false
 \global\advance\island@@@C\@ne
 {\noexpands@
 \xdef\Thelabel@@@{\number\island@@@C}%
 \xdefThelabel@\island@@@N
 \xdef\Thelabel@@@@{\island@@@P\Thelabel@\island@@@Q}%
 \xdefThelabel@@\island@@@S
 \xdef\Thepref@{\Thelabel@@@@}}%
 \finishcaption@}
\long\def\captionformat@#1#2#3{\rm\strut#1 {\island@@@F#2} #3%
 \punct@.\strut}
\long\def\widerthanisland@#1#2#3{\test@true\setbox\z@\vbox{\hsize\maxdimen
 \noindent@@\captionformat@{#1}{#2}{#3}\par\setboxzl@}%
 \ifdim\wdz@=\z@
  \global\setbox\captionbox@\hbox{\noset@\unlabel@
   \captionformat@{#1}{#2}{#3}}%
  \ifdim\wd\captionbox@>\wd\islandbox@\else\test@false\fi
 \fi}
\long\def\captionformat@@#1#2#3{\widerthanisland@{#1}{#2}{#3}%
 \iftest@
  \global\setbox\captionbox@\vbox{\hsize\wd\islandbox@
   \vskip-\parskip\noindent@@\noset@\unlabel@
   \captionformat@{#1}{#2}{#3}\par}%
 \else
  \global\setbox\captionbox@
   \hbox to\wd\islandbox@{\hfil\box\captionbox@\hfil}%
 \fi}
\long\def\finishcaption@#1{\def\entry@{#1}%
 {\locallabel@
 \captionformat@@
  {\expandafter\ignorespaces\captionprefix@\unskip}%
  {\ifx\thelabel@@\empty\unskip\else\thelabel@@\fi}%
  {\ignorespaces#1\unskip}%
 \ifnum\captioncount@=\@ne
  \global\setbox\islandbox@\vbox{\ticwrite@\vbox{\box\islandbox@}}%
  \global\setbox\Captionbox@\vbox{\box\captionbox@}%
 \else
  \global\setbox\islandbox@\vbox{\unvbox\islandbox@\setboxzl@
   \ticwrite@\boxz@}%
  \global\setbox\Captionbox@\vbox{\unvbox\Captionbox@
   \smallskip\box\captionbox@}%
 \fi}%
 \nopunct@false\nospace@false\ignorespaces}
\def\Sixtic@{\ifx\macdef@\empty\else
 \DN@##1##2\next@{\def\macdef@{##1##2}}%
 \expandafter\next@\macdef@\next@
 \edef\next@
  {\noexpand\six@\tic@\macdef@
  \space\space\space\space\space\space\space\space\space\space\space\space
  \noexpand\six@}%
 \next@\let\macdef@\relax\fi}
\def\ticwrite@{%
 \iftoc@
  {\noexpands@\let\style\relax
  \DN@{\island}%
  \edef\next@{\write\tic@{%
   \ifnopunct@\noexpand\noexpand\noexpand\nopunct\fi
   \ifx\islandtype@\next@\noexpand\noexpand\noexpand\island
    \noexpand\string\noexpand\c{\islandclass@}{\captionprefix@}%
     {\QorThelabel@@@@}\else\noexpand\noexpand\expandafter\noexpand
     \islandtype@{\QorThelabel@@@@}}\fi}%
  \next@}%
  \expandafter\unmacro@\meaning\entry@\unmacro@
  \Sixtic@
  \write\tic@{\noexpand\Page{\number\pageno}{\page@N}{\page@P}{\page@Q}^^J}%
 \fi}
\def\Htrim@#1{%
 \ifHby@
  \dimen@\vsize
  \ifnum\captioncount@=\z@
  \else
   \advance\dimen@-\ht\Captionbox@
   \advance\dimen@-#1%
  \fi
  \global\Hby@false
  \dimen@ii\wd\islandbox@
  \global\setbox\islandbox@\vbox
   {\unvbox\islandbox@\setboxzl@
   \vbox to\z@{\vss\boxz@}\nointerlineskip\hbyw\dimen@\dimen@ii}%
  \global\Hby@true
 \fi}
\newif\ifdata@
\def\iclasstest@#1{\DN@{#1}\ifx\next@\islandclass@
 \test@true\else\test@false\fi}
\skipdef\skipi@=1
\def\endisland{\ifnum\captioncount@=\z@\expandafter\egroup\fi
 \ifdata@
 \else
  \iclasstest@{T}%
  \iftest@
   {\rm\global\skipi@-\dp\strutbox}\global\advance\skipi@\bigskipamount
   \Htrim@\skipi@
   \global\setbox\islandbox@\vbox
    {\ifnum\captioncount@=\z@\else
     \box\Captionbox@
     \nointerlineskip
     \vskip\skipi@\fi
     \box\islandbox@}%
  \else
   {\rm\global\skipi@\dp\strutbox}\global\advance\skipi@\medskipamount
   \Htrim@\skipi@
   \global\setbox\islandbox@\vbox
    {\box\islandbox@
     \ifnum\captioncount@=\z@\else
     \nointerlineskip
     \vskip\skipi@
     \box\Captionbox@
     \fi}%
  \fi
  \ifHby@
  \else
   \dimen@\ht\islandbox@\advance\dimen@\dp\islandbox@
   \ifdim\dimen@>\vsize
    \DN@{\island}%
    \Err@{%
     \ifx\islandtype@\next@\noexpand\island\else
      \expandafter\noexpand\islandtype@\fi
     \ifnum\captioncount@=\z@\else
       with \noexpand\caption\fi
      is larger than page}%
     \ht\islandbox@=\vsize
   \fi
  \fi
 \fi
 \global\Hby@false\island@true}
\def\newisland#1\c#2#3{\define#1{}%
 \iftoc@\immediate\write\tic@{\noexpand\newisland\noexpand#1%
  \string\c{#2}{#3}^^J}\fi
 \expandafter\def\csname\exstring@#1@S\endcsname{\island@S}%
 \expandafter\def\csname\exstring@#1@N\endcsname{\island@N}%
 \expandafter\def\csname\exstring@#1@P\endcsname{\island@P}%
 \expandafter\def\csname\exstring@#1@Q\endcsname{\island@Q}%
 \expandafter\def\csname\exstring@#1@F\endcsname{\island@F}%
 \expandafter\def\csname end\exstring@#1\endcsname{\endisland}%
 \expandafter
 \ifx\csname island@C#2\endcsname\relax
  \expandafter\newcount@\csname island@C#2\endcsname
  \global\csname island@C#2\endcsname\z@
 \fi
 \edef\next@{\noexpand\expandafter\noexpand\let\noexpand
  \csname\exstring@#1@C\noexpand\endcsname
  \csname island@C#2\endcsname}%
 \next@
 \def#1{\def\islandtype@{#1}\island@c\c{#2}{#3}}}
\newisland\Figure\c{F}{Figure}
\newisland\Table\c{T}{Table}
\newbox\islandboxi
\newbox\islandboxii
\newbox\islandboxiii
\newbox\captionboxi
\newbox\captionboxii
\newbox\captionboxiii
\long\def\islandpairdata#1#2{{\data@true
 \place@true
 #1%
 \global\setbox\islandboxi\box\islandbox@
 \global\setbox\captionboxi\box\Captionbox@
 #2%
 \global\setbox\islandboxii\box\islandbox@
 \global\setbox\captionboxii\box\Captionbox@
 }}
\long\def\islandpairbox#1#2{\islandpairdata{#1}{#2}%
 \dimen@\ht\captionboxi
 \ifdim\ht\captionboxii>\dimen@\dimen@\ht\captionboxii\fi
 \ifdim\dimen@>\z@
  \ifdim\ht\captionboxi<\dimen@
   \global\setbox\captionboxi\vbox to\dimen@{\unvbox\captionboxi\vfil}\fi
  \ifdim\ht\captionboxii<\dimen@
   \global\setbox\captionboxii\vbox to\dimen@{\unvbox\captionboxii\vfil}\fi
 \fi
 \global\setbox\islandbox@\vbox
 {\hbox to\hsize{\hfil\box\islandboxi\hfil\box\islandboxii\hfil}%
 \ifdim\dimen@>\z@\nointerlineskip
 {\rm\global\skipi@\dp\strutbox}\global\advance\skipi@\medskipamount
  \vskip\skipi@
  \hbox to\hsize{\hfil\box\captionboxi\hfil\box\captionboxii\hfil}\fi}}	
\long\def\islandpairboxa#1#2{\islandpairdata{#1}{#2}%
 \dimen@\ht\captionboxi
 \ifdim\ht\captionboxii>\dimen@\dimen@\ht\captionboxii\fi
 \ifdim\dimen@>\z@
  \ifdim\ht\captionboxi<\dimen@
   \global\setbox\captionboxi\vbox to\dimen@{\vfil\unvbox\captionboxi}\fi
  \ifdim\ht\captionboxii<\dimen@
   \global\setbox\captionboxii\vbox to\dimen@{\vfil\unvbox\captionboxii}\fi
 \fi
 \dimen@ii\ht\islandboxi
 \ifdim\ht\islandboxii>\dimen@ii \dimen@ii\ht\islandboxii\fi
 \ifdim\dimen@ii>\z@
  \ifdim\ht\islandboxi<\dimen@ii
   \global\setbox\islandboxi\vbox to\dimen@ii{\box\islandboxi\vfil}\fi
  \ifdim\ht\islandboxii<\dimen@ii
   \global\setbox\islandboxii\vbox to\dimen@ii{\box\islandboxii\vfil}\fi
 \fi
 \global\setbox\islandbox@\vbox{\ifdim\dimen@>\z@
  \hbox to\hsize{\hfil\box\captionboxi\hfil\box\captionboxii\hfil}%
  \nointerlineskip{\rm\global\skipi@-\dp\strutbox}%
  \global\advance\skipi@\bigskipamount\vskip\skipi@\fi
  \hbox to\hsize{\hfil\box\islandboxi\hfil\box\islandboxii\hfil}}}
\long\def\islandtripledata#1#2#3{{\data@true\place@true
 #1%
 \global\setbox\islandboxi\box\islandbox@
 \global\setbox\captionboxi\box\Captionbox@
 #2%
 \global\setbox\islandboxii\box\islandbox@
 \global\setbox\captionboxii\box\Captionbox@
 #3%
 \global\setbox\islandboxiii\box\islandbox@
 \global\setbox\captionboxiii\box\Captionbox@
 }}
\long\def\islandtriplebox#1#2#3{\islandtripledata{#1}{#2}{#3}%
 \dimen@\ht\captionboxi
 \ifdim\ht\captionboxii>\dimen@ \dimen@\ht\captionboxii\fi
 \ifdim\ht\captionboxiii>\dimen@ \dimen@\ht\captionboxiii\fi
 \ifdim\dimen@>\z@
  \ifdim\ht\captionboxi<\dimen@
   \global\setbox\captionboxi\vbox to\dimen@{\unvbox\captionboxi\vfil}\fi
  \ifdim\ht\captionboxii<\dimen@
   \global\setbox\captionboxii\vbox to\dimen@{\unvbox\captionboxii\vfil}\fi
  \ifdim\ht\captionboxiii<\dimen@
   \global\setbox\captionboxiii\vbox to\dimen@{\unvbox\captionboxiii\vfil}\fi
 \fi
 \global\setbox\islandbox@\vbox
  {\hbox to\hsize{\hfil\box\islandboxi\hfil\box\islandboxii\hfil
   \box\islandboxiii\hfil}%
 \ifdim\dimen@>\z@\nointerlineskip
  {\rm\global\skipi@\dp\strutbox}\global\advance\skipi@\medskipamount
  \vskip\skipi@
  \hbox to\hsize{\hfil\box\captionboxi\hfil\box\captionboxii\hfil
   \box\captionboxiii\hfil}\fi}}
\def\islandtripleboxa#1#2#3{\islandtripledata{#1}{#2}{#3}%
 \dimen@\ht\captionboxi
 \ifdim\ht\captionboxii>\dimen@ \dimen@\ht\captionboxii\fi
 \ifdim\ht\captionboxiii>\dimen@ \dimen@\ht\captionboxiii\fi
 \ifdim\dimen@>\z@
  \ifdim\ht\captionboxi<\dimen@
   \global\setbox\captionboxi\vbox to\dimen@{\vfil\unvbox\captionboxi}\fi
  \ifdim\ht\captionboxii<\dimen@
   \global\setbox\captionboxii\vbox to\dimen@{\vfil\unvbox\captionboxii}\fi
  \ifdim\ht\captionboxiii<\dimen@
   \global\setbox\captionboxiii\vbox to\dimen@{\vfil\unvbox\captionboxiii}\fi
 \fi
 \dimen@ii\ht\islandboxi
 \ifdim\ht\islandboxii>\dimen@ii \dimen@ii\ht\islandboxii\fi
 \ifdim\ht\islandboxiii>\dimen@ii \dimen@ii\ht\islandboxiii\fi
 \ifdim\dimen@ii>\z@
  \ifdim\ht\islandboxi<\dimen@ii
   \global\setbox\islandboxi\vbox to\dimen@ii{\box\islandboxi\vfil}\fi
  \ifdim\ht\islandboxii<\dimen@ii
   \global\setbox\islandboxii\vbox to\dimen@ii{\box\islandboxii\vfil}\fi
  \ifdim\ht\islandboxiii<\dimen@ii
   \global\setbox\islandboxiii\vbox to\dimen@ii{\box\islandboxiii\vfil}\fi
 \fi
 \global\setbox\islandbox@\vbox
  {\ifdim\dimen@>\z@
  \hbox to\hsize{\hfil\box\captionboxi\hfil\box\captionboxii\hfil
   \box\captionboxiii\hfil}%
  \nointerlineskip{\rm\global\skipi@-\dp\strutbox}%
  \global\advance\skipi@\bigskipamount\vskip\skipi@\fi
  \hbox to\hsize{\hfil\box\islandboxi\hfil\box\islandboxii\hfil
   \box\islandboxiii\hfil}}}
\def\Figurepair#1\and#2\endFigurepair{\island@true
 \islandpairbox{\Figure#1\endFigure}{\Figure#2\endFigure}}
\def\Figuretriple#1\and#2\and#3\endFiguretriple{\island@true
 \islandtriplebox{\Figure#1\endFigure}{\Figure#2\endFigure}%
  {\Figure#3\endFigure}}
\def\Tablepair#1\and#2\endTablepair{\island@true
 \islandpairboxa{\Table#1\endTable}{\Table#2\endTable}}
\def\Tabletriple#1\and#2\and#3\endTabletriple{\island@true
 \islandtripleboxa{\Table#1\endTable}{\Table#2\endTable}%
 {\Table#3\endTable}}
\def\place#1{\place@true\island@false
 #1%
 \ifisland@
  \box\islandbox@
 \else
  \Err@{Whoa ... there's no \string\Figure, \string\Table,
   etc., here}%
 \fi
 \place@false}
\newskip\belowtopfigskip
\belowtopfigskip 15\p@ plus 5\p@ minus5\p@
\newskip\abovebotfigskip
\abovebotfigskip 18\p@ plus 6\p@ minus6\p@
\newdimen\minpagesize
\minpagesize 5pc
\dimen@\belowtopfigskip
\advance\dimen@-\abovebotfigskip
\skip\topins\dimen@
\dimen\topins\z@
\newcount\topinscount@
\newbox\topinsdims@
\def\storedim@{\global\setbox\topinsdims@
 \vbox{\hbox to\dimen@{}\unvbox\topinsdims@}}
\def\advancedimtopins@{%
 \ifnum\pageno=\@ne
 \else
   \advance\dimen@\dimen\topins
   \global\dimen\topins\dimen@
 \fi}
\newcount\flipcount@
\def\fliptopins@{%
 \global\flipcount@\z@
 \ifvoid\topins\else
 \setbox\z@\vbox
  {\vskip\p@
   \unvbox\topins
   \global\setbox\topins\vbox{}%
   \loop
    \test@false
    \ifdim\lastskip=\z@\unskip
     \ifdim\lastskip=\z@
      \test@true\fi\fi
    \iftest@
    \global\advance\flipcount@\@ne
    \setboxzl@
    \global\setbox\topins\vbox{\unvbox\topins\boxz@}%
    \unpenalty
   \repeat}\fi}
\newif\ifPar@
\newcount\Parcount@
\newbox\Parbox@
\expandafter\newbox\csname Parfigbox1\endcsname
\expandafter\newbox\csname Parfigbox2\endcsname
\expandafter\newbox\csname Parfigbox3\endcsname
\expandafter\newbox\csname Parfigbox4\endcsname
\expandafter\newbox\csname Parfigbox5\endcsname
\expandafter\newdimen\csname Parprev1\endcsname
\expandafter\newdimen\csname Parprev2\endcsname
\expandafter\newdimen\csname Parprev3\endcsname
\expandafter\newdimen\csname Parprev4\endcsname
\expandafter\newdimen\csname Parprev5\endcsname
\expandafter\newdimen\csname Parprev6\endcsname
\def\Par{\par\global\csname Parprev1\endcsname\prevdepth
 \global\Parcount@\@ne
 \global\Par@true\global\let\Parlist@\empty
 \global\setbox\Parbox@\vbox\bgroup\break}
\def\place@#1#2{%
 \ifisland@
  \ifhmode
   \ifPar@
    \ifnum\Parcount@>5
     \Err@{Only 5 \string\place's allowed per
      \string\Par...\noexpand\endPar paragraph}%
    \else
     \expandafter\expandafter\expandafter
      \global\expandafter\setbox
       \csname Parfigbox\number\Parcount@\endcsname\box\islandbox@
     \global\advance\Parcount@\@ne
     \xdef\Parlist@{\Parlist@#1}%
    \fi
   \else
    \vadjust{#2}%
   \fi
  \else
   #2%
  \fi
 \else
  \Err@{Whoa ... there's no \string\Figure,
   \string\Table, etc., here}%
 \fi
 \place@false}
\long\def\Aplace#1{\prevanish@
 \place@true\island@false
 #1%
 \place@ a\Aplace@
 \postvanish@}
\long\def\AAplace#1{\prevanish@\place@true\island@false
 #1%
 \place@ A\AAplace@
 \postvanish@}
\newif\ifAA@
\def\AAplace@{\AA@true\Aplace@\AA@false}
\let\AAlist@\empty
\def\Aplace@{\allowbreak
 \dimen@=\ht\islandbox@
 \advance\dimen@\abovebotfigskip
 \ht\islandbox@\dimen@
 \advance\dimen@\dp\islandbox@
 \storedim@
 \ifAA@
  \xdef\AAlist@{\AAlist@1}%
  \advancedimtopins@
 \else
  \xdef\AAlist@{\AAlist@0}%
  \ifnum\topinscount@>\@ne\else\advancedimtopins@\fi
 \fi
 \insert\topins{\penalty\z@\splittopskip\z@\floatingpenalty\z@
  \box\islandbox@}%
 \global\advance\topinscount@\@ne}
\long\def\Bplace#1{\prevanish@\place@true\island@false
 #1%
 \place@ b\Bplace@
 \postvanish@}
\def\Bplace@{\allowbreak
 \ifnum\topinscount@=\z@
  \setbox\z@\vbox{\vbox to-\belowtopfigskip{}}%
  \dimen@-\skip\topins
  \ht\z@\dimen@
  \storedim@
  \advancedimtopins@
  \insert\topins{\boxz@}%
  \global\advance\topinscount@\@ne
  \xdef\AAlist@{\AAlist@0}%
 \fi
 \dimen@\ht\islandbox@
 \advance\dimen@\abovebotfigskip
 \ht\islandbox@\dimen@
 \advance\dimen@\dp\islandbox@
 \storedim@
 \xdef\AAlist@{\AAlist@0}%
 \ifnum\topinscount@>\@ne\else\advancedimtopins@\fi
 \insert\topins{\penalty\z@\splittopskip\z@
  \floatingpenalty\z@
  \box\islandbox@}%
 \global\advance\topinscount@\@ne}
\def\breakisland@{\global\setbox\@ne\lastbox\global\skipi@\lastskip\unskip
 \global\setbox\thr@@\lastbox}%
\def\printisland@{\centerline{\box\thr@@}\nobreak\nointerlineskip
 \vskip\skipi@
 \ifdim\ht\@ne<\z@\box\@ne\else\centerline{\box\@ne}\fi}
\def\bottomfigs@{%
 \count@\@ne
 \loop
  \ifnum\count@<\flipcount@
  \nointerlineskip
  \vskip\abovebotfigskip
  \global\setbox\topins\vbox{\unvbox\topins\setboxzl@
   \unvbox\z@
   \breakisland@}%
  \printisland@
  \advance\count@\@ne
  \repeat}
\def\resetdimtopins@{%
 \global\advance\topinscount@-\flipcount@
 \global\setbox\topinsdims@\vbox
  {\unvbox\topinsdims@
   \count@\z@
   \DN@##1##2\next@{\gdef\AAlist@{##2}}%
   \loop
    \ifnum\count@<\flipcount@\setboxzl@
    \expandafter\next@\AAlist@\next@
    \advance\count@\@ne
    \repeat
   \dimen@\z@
   \count@\z@
   \setbox\tw@\vbox{}%
   \edef\nextiii@{\AAlist@}%
   \DN@##1##2\next@{\DNii@{##1}\def\nextiii@{##2}}%
   \loop
    \test@false
    \ifnum\count@<\topinscount@
    \expandafter\next@\nextiii@\next@
     \ifnum\count@<\tw@
      \test@true
     \else
      \if\nextii@ 1\test@true\fi
     \fi
    \fi
    \iftest@
     \setboxzl@
     \advance\dimen@\wdz@
     \setbox\tw@\vbox{\boxz@\unvbox\tw@}%
     \advance\count@\@ne
    \repeat
    \unvbox\tw@
    \global\dimen\topins\dimen@}}
\def\Place@#1#2{%
 \ifisland@
  \ifhmode
   \ifPar@
    \ifnum\Parcount@>5
     \Err@{Only 5 \string\place's allowed per
       \string\Par...\noexpand\endPar paragraph}%
    \else
     \expandafter\expandafter\expandafter\global\expandafter\setbox
      \csname Parfigbox\number\Parcount@\endcsname\box\islandbox@
     \global\advance\Parcount@\@ne
     \xdef\Parlist@{\Parlist@#1}%
     \vadjust{\break}%
    \fi
   \else
    \Err@{\noexpand#2allowed only in a \string\Par...\noexpand\endPar
     paragraph}%
   \fi
  \else
   #2%
  \fi
 \else
  \Err@{Who ... there's no \string\Figure, \string\Table,
   etc., here}%
 \fi
 \place@false}
\newif\ifC@
\newdimen\Cdim@
\long\def\Cplace#1{\prevanish@\place@true\island@false
 #1%
 \Place@ c\Cplace@
 \postvanish@}
\def\Cplace@{\allowbreak
 \ifnum\topinscount@>\z@\else
  \global\C@true\global\Cdim@\pagetotal\fi
 \Aplace@}
\long\def\Mplace#1{\prevanish@\place@true\island@false
 #1%
 \Place@ m\Mplace@
 \postvanish@}
\long\def\MXplace#1{\prevanish@\place@true\island@false
 #1%
 \Place@ M\MXplace@
 \postvanish@}
\newif\ifMX@
\def\MXplace@{\MX@true\Mplace@\MX@false}
\def\Mplace@{\allowbreak
 \dimen@\ht\islandbox@\advance\dimen@\dp\islandbox@
 \ifdim\pagetotal=\z@\else
  \ifdim\lastskip<\abovebotfigskip\advance\dimen@\abovebotfigskip
  \advance\dimen@-\lastskip\fi
 \fi
 \advance\dimen@\pagetotal
 \ifdim\dimen@>\pagegoal
  \Aplace@
 \else
  \nointerlineskip
  \ifdim\lastskip<\abovebotfigskip\removelastskip\vskip\abovebotfigskip\fi
  \setbox\z@\vbox{\unvbox\islandbox@
   \breakisland@}%
  \printisland@
  \ifnum\topinscount@=\z@
   \setbox\z@\vbox{\vbox to-\belowtopfigskip{}}%
   \dimen@-\skip\topins
   \ht\z@\dimen@
   \storedim@
   \advancedimtopins@
   \insert\topins{\boxz@}%
   \global\advance\topinscount@\@ne
   \xdef\AAlist@{\AAlist@0}%
  \fi
  \ifMX@
   \ifnum\topinscount@=\@ne
    \setbox\z@\vbox{\vbox to-\abovebotfigskip{}}%
    \ht\z@\z@
    \dimen@\z@
    \storedim@
    \advancedimtopins@
    \insert\topins{\boxz@}%
    \global\advance\topinscount@\@ne
    \xdef\AAlist@{\AAlist@0}%
   \fi
  \fi
  \nointerlineskip
  \vskip\belowtopfigskip
 \fi}
\expandafter\newbox\csname Parbox1\endcsname
\expandafter\newbox\csname Parbox2\endcsname
\expandafter\newbox\csname Parbox3\endcsname
\expandafter\newbox\csname Parbox4\endcsname
\expandafter\newbox\csname Parbox5\endcsname
\def\endPar{\egroup
 \count@\@ne
 {\vbadness\@M\vfuzz\maxdimen\splitmaxdepth\maxdimen\splittopskip\ht\strutbox
 \setbox\z@\vsplit\Parbox@ to\ht\Parbox@
 \loop
  \ifnum\count@<\Parcount@
  \expandafter\expandafter\expandafter\global\expandafter\setbox
   \csname Parbox\number\count@\endcsname\vsplit\Parbox@ to\ht\Parbox@
  \count@@\count@\advance\count@@\@ne
  \global\csname Parprev\number\count@@\endcsname
   \dp\csname Parbox\number\count@\endcsname
  \advance\count@\@ne
  \repeat}%
 \vskip\parskip
 \count@\@ne
 \def\nextv@##1##2\nextv@{\DN@{##1}\gdef\Parlist@{##2}}%
 \loop
  \ifnum\count@<\Parcount@
   \dimen@\csname Parprev\number\count@\endcsname
   \advance\dimen@\ht\strutbox
   \ifdim\dimen@<\baselineskip
    \advance\dimen@-\baselineskip\vskip-\dimen@
   \else
    \vskip\lineskip
   \fi
   \unvbox\csname Parbox\number\count@\endcsname
   \global\setbox\islandbox@\box\csname Parfigbox\number\count@\endcsname
   \expandafter\nextv@\Parlist@\nextv@
   \if a\next@\Aplace@\else
   \if A\next@\AAplace@\else
   \if b\next@\Bplace@\else
   \if c\next@\Cplace@\else
   \if m\next@\Mplace@\else
   \if M\next@\MXplace@\fi\fi\fi\fi\fi\fi
  \advance\count@\@ne
  \repeat
 \global\Par@false
 \ifvoid\Parbox@
  \prevdepth\csname Parprev\number\count@\endcsname
 \else
  \dimen@\csname Parprev\number\count@\endcsname\advance\dimen@\ht\strutbox
  \ifdim\dimen@<\baselineskip
    \advance\dimen@-\baselineskip\vskip-\dimen@
  \else
    \vskip\lineskip
  \fi
  \dimen@\dp\Parbox@
  \unvbox\Parbox@
  \prevdepth\dimen@
 \fi}
\def\folio{{\page@F\page@S{\page@P\page@N{\number\page@C}\page@Q}}}
\def\advancepageno{\global\advance\pageno\@ne}
\newif\ifspecialsplit@
\newbox\outbox@
\let\shipout@\shipout
\def\plainoutput{\specialsplit@false\ifvoid\topins\else\ifdim\ht\topins=\z@
 \specialsplit@true\advance\minpagesize-\skip\topins\fi\fi
 \fliptopins@
 \setbox\outbox@\vbox{\makeheadline\pagebody\makefootline}%
 {\noexpands@\let\style\relax
 \shipout@\box\outbox@}%
 \advancepageno
 \resetdimtopins@
 \ifvoid\@cclv\else\unvbox\@cclv\penalty\outputpenalty\fi
 \ifnum\outputpenalty>-\@MM\else\dosupereject\fi}
\def\pagebody{\vbox to\vsize{\boxmaxdepth\maxdepth
 \ifvoid\margin@\else
 \rlap{\kern\hsize\vbox to\z@{\kern4\p@\box\margin@\vss}}\fi
 \pagecontents}}
\newif\ifonlytop@
\def\pagecontents{%
 \onlytop@false
 \ifdim\ht\@cclv<\minpagesize\ifnum\flipcount@<\tw@\ifvoid\footins
  \onlytop@true\fi\fi\fi
 \test@false
 \ifC@
  \ifnum\flipcount@=\@ne
   \global\multiply\Cdim@\tw@
   \ifdim\Cdim@>\ht\@cclv
    \test@true
   \fi
  \fi
 \fi
 \global\C@false
 \iftest@
  \dimen@\ht\@cclv
  \advance\dimen@\skip\topins
  {\vfuzz\maxdimen\vbadness\@M
  \splitmaxdepth\maxdepth\splittopskip\topskip
  \setbox\z@\vsplit\@cclv to\dimen@
  \unvbox\z@}%
  \global\setbox\topins\vbox{\unvbox\topins
   \global\setbox\@ne\lastbox}%
  \setbox\z@\vbox{\unvbox\@ne
   \breakisland@}%
  \nointerlineskip
  \vskip\abovebotfigskip
  \printisland@
 \else
  \ifnum\flipcount@>\z@
   \global\setbox\topins\vbox{\unvbox\topins\global\setbox\@ne\lastbox}%
   \setbox\z@\vbox{\unvbox\@ne
    \breakisland@}%
   \printisland@
   \ifonlytop@\kern-\prevdepth\vfill\else\vskip\belowtopfigskip\fi
  \fi
 \fi
 \ifdim\ht\@cclv<\minpagesize
  \ifonlytop@\else\vfill\fi
 \else
  \ifspecialsplit@
   {\vfuzz\maxdimen\vbadness\@M
   \splitmaxdepth\maxdepth\splittopskip\topskip
   \dimen@ii\ht\@cclv \advance\dimen@ii\skip\topins
   \setbox\z@\vsplit\@cclv to\dimen@ii
   \unvbox\z@}%
  \else
   \unvbox\@cclv
  \fi
 \fi
 \bottomfigs@
 \ifvoid\footins\else\vskip\skip\footins\footnoterule\unvbox\footins\fi}
\newread\readdata@
\def\readthedata@#1{\expandafter
 \ifx\csname#1@D\endcsname\relax
  \immediate\openin\readdata@=#1.dat
  \ifeof\readdata@
   \Err@{No file #1.dat}%
  \else
   {\endlinechar\m@ne\gdef\Next@{}%
   \DNii@##1 ##2 ##3pt{\global\data@ht##1\global\data@dp##2%
    \global\data@wd##3pt}%
   \loop
    \ifeof\readdata@
    \else
    \read\readdata@ to\next@
    \ifx\next@\empty\else
     \edef\next@{\expandafter\nextii@\next@}%
     \expandafter\rightadd@\next@\to\Next@
    \fi
    \repeat}%
   \immediate\closein\readdata@
   \expandafter\expandafter\expandafter\global\expandafter
    \let\csname#1@D\endcsname\Next@\global\let\Next@\relax
  \fi
 \fi}
\newdimen\data@ht
\newdimen\data@dp
\newdimen\data@wd
\newif\ifgetdata@
\def\getdata@#1#2{\global\getdata@true\count@#2\relax
 {\let\\\or\xdef\Next@{\ifcase\number\count@#1\else
 \global\noexpand\getdata@false\fi}}\Next@}
\def\paste#1#2{\readthedata@{#1}%
 \getdata@{\csname#1@D\endcsname}{#2}%
 \ifgetdata@
 \dimen@\data@ht \advance\dimen@\data@dp
  \hbox{\special{dvipaste: #1 #2}%
   \lower\data@dp\vbox to\dimen@{\hbox to\data@wd{}\vfil}}%
 \else
  {\lccode`\Z=`\#\lccode`\N=`\N\lccode`\F=`\F%
   \lowercase{\Err@{No data for File [#1], Z#2}}}%
 \fi}
\newdimen\httable
\newdimen\dptable
\newdimen\wdtable
\def\measuretable#1#2{\readthedata@{#1}%
 \getdata@{\csname#1@D\endcsname}{#2}%
 \ifgetdata@
  \httable\data@ht \dptable\data@dp \wdtable\data@wd
 \else
  {\lccode`\Z=`\#\lccode`\N=`\N\lccode`\F=`\F%
  \lowercase{\Err@{No data for File [#1], Z#2}}}%
 \fi}
\def\East#1#2{\setboxz@h{$\m@th\ssize\;{#1}\;\;$}%
 \setbox\tw@\hbox{$\m@th\ssize\;{#2}\;\;$}\setbox4=\hbox{$\m@th#2$}%
 \dimen@\minaw@
 \ifdim\wdz@>\dimen@\dimen@\wdz@\fi\ifdim\wd\tw@>\dimen@\dimen@\wd\tw@\fi
 \ifdim\wd4 >\z@
  \mathrel{\mathop{\hbox to\dimen@{\rightarrowfill}}\limits^{#1}_{#2}}%
 \else
  \mathrel{\mathop{\hbox to\dimen@{\rightarrowfill}}\limits^{#1}}%
 \fi}
\def\West#1#2{\setboxz@h{$\m@th\ssize\;\;{#1}\;$}%
 \setbox\tw@\hbox{$\m@th\ssize\;\;{#2}\;$}\setbox4=\hbox{$\m@th#2$}%
 \dimen@\minaw@
 \ifdim\wdz@>\dimen@\dimen@\wdz@\fi\ifdim\wd\tw@>\dimen@\dimen@\wd\tw@\fi
 \ifdim\wd4 >\z@
  \mathrel{\mathop{\hbox to\dimen@{\leftarrowfill}}\limits^{#1}_{#2}}%
 \else
  \mathrel{\mathop{\hbox to\dimen@{\leftarrowfill}}\limits^{#1}}%
 \fi}
\font\arrow@i=lams1
\font\arrow@ii=lams2
\font\arrow@iii=lams3
\font\arrow@iv=lams4
\font\arrow@v=lams5
\newdimen\standardcgap
\standardcgap40\p@
\newdimen\hunit
\hunit\tw@\p@
\newdimen\standardrgap
\standardrgap32\p@
\newdimen\vunit
\vunit1.6\p@
\def\Cgaps#1{\RIfM@
 \standardcgap#1\standardcgap\relax\hunit#1\hunit\relax
 \else\nonmatherr@\Cgaps\fi}
\def\Rgaps#1{\RIfM@
 \standardrgap#1\standardrgap\relax\vunit#1\vunit\relax
 \else\nonmatherr@\Rgaps\fi}
\newdimen\getdim@
\def\getcgap@#1{\ifcase#1\or\getdim@\z@\else\getdim@\standardcgap\fi}
\def\getrgap@#1{\ifcase#1\getdim@\z@\else\getdim@\standardrgap\fi}
\def\cgaps{\RIfM@\expandafter\cgaps@\else\expandafter\nonmatherr@
 \expandafter\cgaps\fi}
\def\cgaps@{\ifnum\catcode`\;=\active\expandafter\cgapsA@\else
 \expandafter\cgapsO@\fi}
\def\cgapsO@#1{\toks@{\ifcase\i@\or\getdim@=\z@}%
 \gaps@@\standardcgap#1;\gaps@@\gaps@@
 \edef\next@{\the\toks@\noexpand\else\noexpand\getdim@\noexpand\standardcgap
  \noexpand\fi}%
 \toks@=\expandafter{\next@}%
 \edef\getcgap@##1{\i@##1\relax\the\toks@}\toks@{}}
{\catcode`\;=\active
 \gdef\cgapsA@#1{\toks@{\ifcase\i@\or\getdim@=\z@}%
 \gaps@@\standardcgap#1;\gaps@@\gaps@@
 \edef\next@{\the\toks@\noexpand\else\noexpand\getdim@\noexpand\standardcgap
  \noexpand\fi}%
 \toks@=\expandafter{\next@}%
 \edef\getcgap@##1{\i@##1\relax\the\toks@}\toks@{}}
}
\def\Gaps@@{\gaps@@}
\def\gaps@@#1#2;#3{\mgaps@#1#2\mgaps@
 \edef\next@{\the\toks@\noexpand\or\noexpand\getdim@
  \noexpand#1\the\mgapstoks@@}%
 \toks@\expandafter{\next@}%
 \DN@{#3}%
 \ifx\next@\Gaps@@\def\next@##1\gaps@@{}\else
  \def\next@{\gaps@@#1#3}\fi\next@}
{\catcode`\;=\active
 \gdef\rgaps#1{\RIfM@{\ifnum\catcode`\;=\active\def;{\string;}\fi
   \xdef\Next@{\noexpand\rgaps@{#1}}}%
  \Next@\edef\getrgap@##1{\i@##1\relax\the\toks@}\toks@{}\else
  \nonmatherr@\rgaps\fi}
}
\def\rgaps@#1{\toks@{\ifcase\i@\getdim@=\z@}%
 \gaps@@\standardrgap#1;\gaps@@\gaps@@
 \edef\next@{\the\toks@\noexpand\else\noexpand\getdim@\noexpand\standardrgap
  \noexpand\fi}%
 \toks@=\expandafter{\next@}}
\newbox\ZER@
\def\mgaps@#1{\let\mgapsnext@#1\FNSS@\mgaps@@}
\def\mgaps@@{\ifx\next\w\expandafter\mgaps@@@\else
 \expandafter\mgaps@@@@\fi}
\newtoks\mgapstoks@@
\def\mgaps@@@@#1\mgaps@{\getdim@\mgapsnext@\getdim@#1\getdim@
 \edef\next@{\noexpand\getdim@\the\getdim@}%
 \mgapstoks@@\expandafter{\next@}}
\def\mgaps@@@\w#1#2\mgaps@{\mgaps@@@@#2\mgaps@
 \setbox\ZER@\hbox{$\m@th\hskip15\p@\tsize@#1$}%
 \dimen@\wd\ZER@
 \ifdim\dimen@>\getdim@\getdim@\dimen@\fi
 \edef\next@{\noexpand\getdim@\the\getdim@}%
 \mgapstoks@@\expandafter{\next@}}
\def\changewidth#1#2{\setbox\ZER@{$\m@th#2}%
 \hbox to\wd\ZER@{\hss$\m@th#1$\hss}}
\atdef@({\FN@\ARROW@}
\def\ARROW@{\ifx\next)\let\next@\OPTIONS@\else
 \DN@{\csname\string @(\endcsname}\fi\next@}
\newif\ifoptions@
\def\OPTIONS@){\ifoptions@\let\next@\relax\else
 \DN@{\global\options@true\begingroup\optioncodes@}\fi\next@}
\newif\ifN@
\newif\ifE@
\newif\ifNESW@
\newif\ifH@
\newif\ifV@
\newif\ifHshort@
\expandafter\def\csname\string @(\endcsname #1,#2){%
 \ifoptions@\expandafter\endgroup\fi
 \N@false\E@false\H@false\V@false\Hshort@false
 \ifnum#1>\z@\E@true\fi
 \ifnum#1=\z@\V@true\global\tX@false\global\tY@false\global\a@false\fi
 \ifnum#2>\z@\N@true\fi
 \ifnum#2=\z@\H@true\global\tX@false\global\tY@false\global\a@false
  \ifshort@\Hshort@true\fi\fi
 \NESW@false
 \ifN@\ifE@\NESW@true\fi\else\ifE@\else\NESW@true\fi\fi
 \arrow@{#1}{#2}%
 \global\options@false
 \global\scount@\z@\global\tcount@\z@\global\arrcount@\z@
 \global\s@false\global\sxdimen@\z@\global\sydimen@\z@
 \global\tX@false\global\tXdimen@i\z@\global\tXdimen@ii\z@
 \global\tY@false\global\tYdimen@i\z@\global\tYdimen@ii\z@
 \global\a@false\global\exacount@\z@
 \global\x@false\global\xdimen@\z@
 \global\X@false\global\Xdimen@\z@
 \global\y@false\global\ydimen@\z@
 \global\Y@false\global\Ydimen@\z@
 \global\p@false\global\pdimen@\z@
 \global\label@ifalse\global\label@iifalse
 \global\dl@ifalse\global\ldimen@i\z@
 \global\dl@iifalse\global\ldimen@ii\z@
 \global\short@false\global\unshort@false}
\newif\iflabel@i
\newif\iflabel@ii
\newcount\scount@
\newcount\tcount@
\newcount\arrcount@
\newif\ifs@
\newdimen\sxdimen@
\newdimen\sydimen@
\newif\iftX@
\newdimen\tXdimen@i
\newdimen\tXdimen@ii
\newif\iftY@
\newdimen\tYdimen@i
\newdimen\tYdimen@ii
\newif\ifa@
\newcount\exacount@
\newif\ifx@
\newdimen\xdimen@
\newif\ifX@
\newdimen\Xdimen@
\newif\ify@
\newdimen\ydimen@
\newif\ifY@
\newdimen\Ydimen@
\newif\ifp@
\newdimen\pdimen@
\newif\ifdl@i
\newif\ifdl@ii
\newdimen\ldimen@i
\newdimen\ldimen@ii
\newif\ifshort@
\newif\ifunshort@
\def\zero@#1{\ifnum\scount@=\z@
 \if#1e\global\scount@\m@ne\else
 \if#1t\global\scount@\tw@\else
 \if#1h\global\scount@\thr@@\else
 \if#1'\global\scount@6 \else
 \if#1`\global\scount@7 \else
 \if#1(\global\scount@8 \else
 \if#1)\global\scount@9 \else
 \if#1s\global\scount@12 \else
 \if#1H\global\scount@13 \else
 \Err@{\Invalid@@ option \string\0}\fi\fi\fi\fi\fi\fi\fi\fi\fi
 \fi}
\def\one@#1{\ifnum\tcount@=\z@
 \if#1e\global\tcount@\m@ne\else
 \if#1h\global\tcount@\tw@\else
 \if#1t\global\tcount@\thr@@\else
 \if#1'\global\tcount@4 \else
 \if#1`\global\tcount@5 \else
 \if#1(\global\tcount@\ten@ \else
 \if#1)\global\tcount@11 \else
 \if#1s\global\tcount@12 \else
 \if#1H\global\tcount@13 \else
 \Err@{\Invalid@@ option \string\1}\fi\fi\fi\fi\fi\fi\fi\fi\fi
 \fi}
\def\a@#1{\ifnum\arrcount@=\z@
 \if#10\global\arrcount@\m@ne\else
 \if#1+\global\arrcount@\@ne\else
 \if#1-\global\arrcount@\tw@\else
 \if#1=\global\arrcount@\thr@@\else
 \Err@{\Invalid@@ option \string\a}\fi\fi\fi\fi
 \fi}
\def\ds@{\ifnum\catcode`\;=\active\expandafter\dsA@\else
 \expandafter\dsO@\fi}
\def\dsO@(#1;#2){\ds@@{#1}{#2}}
\def\ds@@#1#2{\ifs@\else
 \global\s@true
 \global\sxdimen@\hunit\global\sxdimen@#1\sxdimen@\relax
 \global\sydimen@\vunit\global\sydimen@#2\sydimen@\relax
 \fi}
\def\dtX@{\ifnum\catcode`\;=\active\expandafter\dtXA@\else
 \expandafter\dtXO@\fi}
\def\dtXO@(#1;#2){\dtX@@{#1}{#2}}
\def\dtX@@#1#2{\iftX@\else
 \global\tX@true
 \global\tXdimen@i\hunit\global\tXdimen@i#1\tXdimen@i\relax
 \global\tXdimen@ii\vunit\global\tXdimen@ii#2\tXdimen@ii\relax
 \fi}
\def\dtY@{\ifnum\catcode`\;=\active\expandafter\dtYA@\else
 \expandafter\dtYO@\fi}
\def\dtYO@(#1;#2){\dtY@@{#1}{#2}}
\def\dtY@@#1#2{\iftY@\else
 \global\tY@true
 \global\tYdimen@i\hunit\global\tYdimen@i#1\tYdimen@i\relax
 \global\tYdimen@ii\vunit\global\tYdimen@ii#2\tYdimen@ii\relax
 \fi}
{\catcode`\;=\active
 \gdef\dsA@(#1;#2){\ds@@{#1}{#2}}
 \gdef\dtXA@(#1;#2){\dtX@@{#1}{#2}}
 \gdef\dtYA@(#1;#2){\dtY@@{#1}{#2}}
}
\def\da@#1{\ifa@\else\global\a@true\global\exacount@#1\relax\fi}
\def\dx@#1{\ifx@\else
 \global\x@true
 \global\xdimen@\hunit\global\xdimen@#1\xdimen@\relax
 \fi}
\def\dX@#1{\ifX@\else
 \global\X@true
 \global\Xdimen@\hunit\global\Xdimen@#1\Xdimen@\relax
 \fi}
\def\dy@#1{\ify@\else
 \global\y@true
 \global\ydimen@\vunit\global\ydimen@#1\ydimen@\relax
 \fi}
\def\dY@#1{\ifY@\else
 \global\Y@true
 \global\Ydimen@\vunit\global\Ydimen@#1\Ydimen@\relax
 \fi}
\def\p@@#1{\ifp@\else
 \global\p@true
 \global\pdimen@\hunit\global\divide\pdimen@\tw@
 \global\pdimen@#1\pdimen@\relax
 \fi}
\def\L@#1{\iflabel@i\else
 \global\label@itrue\gdef\label@i{#1}%
 \fi}
\def\l@#1{\iflabel@ii\else
 \global\label@iitrue\gdef\label@ii{#1}%
 \fi}
\def\dL@#1{\ifdl@i\else
 \global\dl@itrue\global\ldimen@i\hunit\global\ldimen@i#1\ldimen@i\relax
 \fi}
\def\dl@#1{\ifdl@ii\else
 \global\dl@iitrue\global\ldimen@ii\hunit\global\ldimen@ii#1\ldimen@ii\relax
 \fi}
\def\s@{\ifunshort@\else\global\short@true\fi}
\def\uns@{\ifshort@\else\global\unshort@true\global\short@false\fi}
\def\optioncodes@{\let\0\zero@\let\1\one@\let\a\a@\let\ds\ds@\let\dtX\dtX@
 \let\dtY\dtY@\let\da\da@\let\dx\dx@\let\dX\dX@\let\dY\dY@\let\dy\dy@
 \let\p\p@@\let\L\L@\let\l\l@\let\dL\dL@\let\dl\dl@\let\s\s@\let\uns\uns@}
\def\slopes@{\\161\\152\\143\\134\\255\\126\\357\\238\\349\\45{10}\\56{11}%
 \\11{12}\\65{13}\\54{14}\\43{15}\\32{16}\\53{17}\\21{18}\\52{19}\\31{20}%
 \\41{21}\\51{22}\\61{23}}
\newcount\tan@i
\newcount\tan@ip
\newcount\tan@ii
\newcount\tan@iip
\newdimen\slope@i
\newdimen\slope@ip
\newdimen\slope@ii
\newdimen\slope@iip
\newcount\angcount@
\newcount\extracount@
\def\slope@{{\slope@i\secondy@\advance\slope@i-\firsty@
 \ifN@\else\multiply\slope@i\m@ne\fi
 \slope@ii\secondx@\advance\slope@ii-\firstx@
 \ifE@\else\multiply\slope@ii\m@ne\fi
 \ifdim\slope@ii<\z@
  \global\tan@i6 \global\tan@ii\@ne\global\angcount@23
 \else
  \dimen@\slope@i\multiply\dimen@6
  \ifdim\dimen@<\slope@ii
   \global\tan@i\@ne\global\tan@ii6 \global\angcount@\@ne
  \else
   \dimen@\slope@ii\multiply\dimen@6
   \ifdim\dimen@<\slope@i
    \global\tan@i6 \global\tan@ii\@ne\global\angcount@23
   \else
    \global\tan@ip\z@\global\tan@iip\@ne
    \def\\##1##2##3{\global\angcount@##3\relax
     \slope@ip\slope@i\slope@iip\slope@ii
     \multiply\slope@iip##1\relax\multiply\slope@ip##2\relax
     \ifdim\slope@iip<\slope@ip
      \global\tan@ip##1\relax\global\tan@iip##2\relax
     \else
      \global\tan@i##1\relax\global\tan@ii##2\relax
      \def\\####1####2####3{}%
     \fi}%
    \slopes@
    \slope@i\secondy@\advance\slope@i-\firsty@
    \ifN@\else\multiply\slope@i\m@ne\fi
    \multiply\slope@i\tan@ii\multiply\slope@i\tan@iip\multiply\slope@i\tw@
    \count@\tan@i\multiply\count@\tan@iip
    \extracount@\tan@ip\multiply\extracount@\tan@ii
    \advance\count@\extracount@
    \slope@ii\secondx@\advance\slope@ii-\firstx@
    \ifE@\else\multiply\slope@ii\m@ne\fi
    \multiply\slope@ii\count@
    \ifdim\slope@i<\slope@ii
     \global\tan@i\tan@ip\global\tan@ii\tan@iip
     \global\advance\angcount@\m@ne
    \fi
   \fi
  \fi
 \fi}%
}
\def\slope@a#1{{\def\\##1##2##3{\ifnum##3=#1\global\tan@i##1\relax
 \global\tan@ii##2\relax\fi}\slopes@}}
\newcount\i@
\newcount\j@
\newcount\colcount@
\newcount\Colcount@
\newcount\tcolcount@
\newdimen\rowht@
\newdimen\rowdp@
\newcount\rowcount@
\newcount\Rowcount@
\newcount\maxcolrow@
\newtoks\colwidthtoks@
\newtoks\Rowheighttoks@
\newtoks\Rowdepthtoks@
\newtoks\widthtoks@
\newtoks\Widthtoks@
\newtoks\heighttoks@
\newtoks\Heighttoks@
\newtoks\depthtoks@
\newtoks\Depthtoks@
\newif\iffirstCDcr@
\def\dotoks@i{%
 \global\widthtoks@\expandafter{\the\widthtoks@\else\getdim@\z@\fi}%
 \global\heighttoks@\expandafter{\the\heighttoks@\else\getdim@\z@\fi}%
 \global\depthtoks@\expandafter{\the\depthtoks@\else\getdim@\z@\fi}}
\def\dotoks@ii{%
 \global\widthtoks@{\ifcase\j@}%
 \global\heighttoks@{\ifcase\j@}%
 \global\depthtoks@{\ifcase\j@}}
\def\preCD@#1\endCD{\setbox\ZER@
 \vbox{%
  \def\arrow@##1##2{{}}%
  \global\rowcount@\m@ne\global\colcount@\z@\global\Colcount@\z@
  \global\firstCDcr@true\toks@{}%
  \global\widthtoks@{\ifcase\j@}%
  \global\Widthtoks@{\ifcase\i@}%
  \global\heighttoks@{\ifcase\j@}%
  \global\Heighttoks@{\ifcase\i@}%
  \global\depthtoks@{\ifcase\j@}%
  \global\Depthtoks@{\ifcase\i@}%
  \global\Rowheighttoks@{\ifcase\i@}%
  \global\Rowdepthtoks@{\ifcase\i@}%
  \Let@
  \everycr{%
   \noalign{%
    \global\advance\rowcount@\@ne
    \ifnum\colcount@<\Colcount@
    \else
     \global\Colcount@\colcount@\global\maxcolrow@\rowcount@
    \fi
    \global\colcount@\z@
    \iffirstCDcr@
     \global\firstCDcr@false
    \else
     \edef\next@{\the\Rowheighttoks@\noexpand\or\noexpand\getdim@\the\rowht@}%
      \global\Rowheighttoks@\expandafter{\next@}%
     \edef\next@{\the\Rowdepthtoks@\noexpand\or\noexpand\getdim@\the\rowdp@}%
      \global\Rowdepthtoks@\expandafter{\next@}%
     \global\rowht@\z@\global\rowdp@\z@
     \dotoks@i
     \edef\next@{\the\Widthtoks@\noexpand\or\the\widthtoks@}%
      \global\Widthtoks@\expandafter{\next@}%
     \edef\next@{\the\Heighttoks@\noexpand\or\the\heighttoks@}%
      \global\Heighttoks@\expandafter{\next@}%
     \edef\next@{\the\Depthtoks@\noexpand\or\the\depthtoks@}%
      \global\Depthtoks@\expandafter{\next@}%
     \dotoks@ii
    \fi}}%
  \tabskip\z@
  \halign{&\setbox\ZER@\hbox{\vrule\height\ten@\p@\width\z@\depth\z@     
   $\m@th\displaystyle{##}$}\copy\ZER@
   \ifdim\ht\ZER@>\rowht@\global\rowht@\ht\ZER@\fi
   \ifdim\dp\ZER@>\rowdp@\global\rowdp@\dp\ZER@\fi
   \global\advance\colcount@\@ne
   \edef\next@{\the\widthtoks@\noexpand\or\noexpand\getdim@\the\wd\ZER@}%
    \global\widthtoks@\expandafter{\next@}%
   \edef\next@{\the\heighttoks@\noexpand\or\noexpand\getdim@\the\ht\ZER@}%
    \global\heighttoks@\expandafter{\next@}%
   \edef\next@{\the\depthtoks@\noexpand\or\noexpand\getdim@\the\dp\ZER@}%
    \global\depthtoks@\expandafter{\next@}%
   \cr#1\crcr}}%
 \Rowcount@\rowcount@
 \global\Widthtoks@\expandafter{\the\Widthtoks@\fi\relax}%
 \edef\Width@##1##2{\i@##1\relax\j@##2\relax\the\Widthtoks@}%
 \global\Heighttoks@\expandafter{\the\Heighttoks@\fi\relax}%
 \edef\Height@##1##2{\i@##1\relax\j@##2\relax\the\Heighttoks@}%
 \global\Depthtoks@\expandafter{\the\Depthtoks@\fi\relax}%
 \edef\Depth@##1##2{\i@##1\relax\j@##2\relax\the\Depthtoks@}%
 \edef\next@{\the\Rowheighttoks@\noexpand\fi\relax}%
 \global\Rowheighttoks@\expandafter{\next@}%
 \edef\Rowheight@##1{\i@##1\relax\the\Rowheighttoks@}%
 \edef\next@{\the\Rowdepthtoks@\noexpand\fi\relax}%
 \global\Rowdepthtoks@\expandafter{\next@}%
 \edef\Rowdepth@##1{\i@##1\relax\the\Rowdepthtoks@}%
 \global\colwidthtoks@{\fi}%
 \setbox\ZER@\vbox{%
  \unvbox\ZER@
  \count@\rowcount@
  \loop
   \unskip\unpenalty
   \setbox\ZER@\lastbox
   \ifnum\count@>\maxcolrow@\advance\count@\m@ne
   \repeat
  \hbox{%
   \unhbox\ZER@
   \count@\z@
   \loop
    \unskip
    \setbox\ZER@\lastbox
    \edef\next@{\noexpand\or\noexpand\getdim@\the\wd\ZER@\the\colwidthtoks@}%
     \global\colwidthtoks@\expandafter{\next@}%
    \advance\count@\@ne
    \ifnum\count@<\Colcount@
    \repeat}}%
 \edef\next@{\noexpand\ifcase\noexpand\i@\the\colwidthtoks@}%
  \global\colwidthtoks@\expandafter{\next@}%
 \edef\Colwidth@##1{\i@##1\relax\the\colwidthtoks@}%
 \global\colwidthtoks@{}\global\Rowheighttoks@{}\global\Rowdepthtoks@{}%
 \global\widthtoks@{}\global\Widthtoks@{}\global\heighttoks@{}%
 \global\Heighttoks@{}\global\depthtoks@{}\global\Depthtoks@{}%
}
\newcount\xoff@
\newcount\yoff@
\newcount\endcount@
\newcount\rcount@
\newdimen\firstx@
\newdimen\firsty@
\newdimen\secondx@
\newdimen\secondy@
\newdimen\tocenter@
\newdimen\charht@
\newdimen\charwd@
\def\outside@{\Err@{This arrow points outside the \string\CD}}
\newif\ifsvertex@
\newif\iftvertex@
\def\arrow@#1#2{\global\xoff@#1\relax\global\yoff@#2\relax
 \count@\rowcount@\advance\count@-\yoff@
 \ifnum\count@<\@ne\outside@\else\ifnum\count@>\Rowcount@\outside@\fi\fi
 \count@\colcount@\advance\count@\xoff@
 \ifnum\count@<\@ne\outside@\else\ifnum\count@>\Colcount@\outside@\fi\fi
 \tcolcount@\colcount@\advance\tcolcount@\xoff@
 \Width@\rowcount@\colcount@\divide\getdim@\tw@\tocenter@-\getdim@
 \ifdim\getdim@=\z@
  \firstx@\z@\firsty@\mathaxis@\svertex@true
 \else
  \svertex@false
  \ifHshort@
   \Colwidth@\colcount@\divide\getdim@\tw@
   \ifE@ \firstx@\getdim@ \else \firstx@-\getdim@ \fi
  \else
   \ifE@ \firstx@\getdim@ \else \firstx@-\getdim@ \fi
  \fi
  \ifE@
   \ifH@ \advance\firstx@\thr@@\p@ \else \advance\firstx@-\thr@@\p@ \fi  
  \else
   \ifH@ \advance\firstx@-\thr@@\p@ \else \advance\firstx@\thr@@\p@ \fi  
  \fi
  \ifN@
   \Height@\rowcount@\colcount@ \firsty@\getdim@                         
   \ifV@ \advance\firsty@\thr@@\p@ \fi                                   
  \else
   \ifV@
    \Depth@\rowcount@\colcount@ \firsty@-\getdim@                        
    \advance\firsty@-\thr@@\p@                                           
   \else
    \firsty@\z@                                                          
   \fi
  \fi
 \fi
 \ifV@
 \else
  \Colwidth@\colcount@\divide\getdim@\tw@
  \ifE@\secondx@\getdim@\else\secondx@-\getdim@\fi
  \ifE@\else\getcgap@\colcount@\advance\secondx@-\getdim@\fi
  \endcount@\colcount@\advance\endcount@\xoff@
  \count@\colcount@
  \ifE@
   \advance\count@\@ne
   \loop
    \ifnum\count@<\endcount@
    \Colwidth@\count@\advance\secondx@\getdim@
    \getcgap@\count@\advance\secondx@\getdim@
    \advance\count@\@ne
    \repeat
  \else
   \advance\count@\m@ne
   \loop
    \ifnum\count@>\endcount@
    \Colwidth@\count@\advance\secondx@-\getdim@
    \getcgap@\count@\advance\secondx@-\getdim@
    \advance\count@\m@ne
    \repeat
  \fi
  \Colwidth@\count@\divide\getdim@\tw@
  \ifHshort@
  \else
   \ifE@\advance\secondx@\getdim@\else\advance\secondx@-\getdim@\fi
  \fi
  \ifE@\getcgap@\count@\advance\secondx@\getdim@\fi
  \rcount@\rowcount@\advance\rcount@-\yoff@
  \Width@\rcount@\count@\divide\getdim@\tw@
  \tvertex@false
  \ifH@\ifdim\getdim@=\z@\tvertex@true\Hshort@false\fi\fi
  \ifHshort@
  \else
   \ifE@\advance\secondx@-\getdim@\else\advance\secondx@\getdim@\fi
  \fi
  \iftvertex@
   \advance\secondx@.4\p@
  \else
   \ifE@\advance\secondx@-\thr@@\p@\else\advance\secondx@\thr@@\p@\fi    
  \fi
 \fi
 \ifH@
 \else
  \ifN@
   \Rowheight@\rowcount@\secondy@\getdim@
  \else
   \Rowdepth@\rowcount@\secondy@-\getdim@
   \getrgap@\rowcount@\advance\secondy@-\getdim@
  \fi
  \endcount@\rowcount@\advance\endcount@-\yoff@
  \count@\rowcount@
  \ifN@
   \advance\count@\m@ne
   \loop
    \ifnum\count@>\endcount@
    \Rowheight@\count@\advance\secondy@\getdim@
    \Rowdepth@\count@\advance\secondy@\getdim@
    \getrgap@\count@\advance\secondy@\getdim@
    \advance\count@\m@ne
    \repeat
  \else
   \advance\count@\@ne
   \loop
    \ifnum\count@<\endcount@
    \Rowheight@\count@\advance\secondy@-\getdim@
    \Rowdepth@\count@\advance\secondy@-\getdim@
    \getrgap@\count@\advance\secondy@-\getdim@
    \advance\count@\@ne
    \repeat
  \fi
  \tvertex@false
  \ifV@\Width@\count@\colcount@\ifdim\getdim@=\z@\tvertex@true\fi\fi
  \ifN@
   \getrgap@\count@\advance\secondy@\getdim@
   \Rowdepth@\count@\advance\secondy@\getdim@
   \iftvertex@
    \advance\secondy@\mathaxis@
   \else
    \Depth@\count@\tcolcount@\advance\secondy@-\getdim@
    \advance\secondy@-\thr@@\p@                                          
   \fi
  \else
   \Rowheight@\count@\advance\secondy@-\getdim@
   \iftvertex@
    \advance\secondy@\mathaxis@
   \else
    \Height@\count@\tcolcount@\advance\secondy@\getdim@
    \advance\secondy@\thr@@\p@                                           
   \fi
  \fi
 \fi
 \ifV@\else\advance\firstx@\sxdimen@\fi
 \ifH@\else\advance\firsty@\sydimen@\fi
 \iftX@
  \advance\secondy@\tXdimen@ii
  \advance\secondx@\tXdimen@i
  \slope@
 \else
  \iftY@
   \advance\secondy@\tYdimen@ii
   \advance\secondx@\tYdimen@i
   \slope@
   \secondy@\secondx@\advance\secondy@-\firstx@
   \ifNESW@\else\multiply\secondy@\m@ne\fi
   \multiply\secondy@\tan@i\divide\secondy@\tan@ii\advance\secondy@\firsty@
  \else
   \ifa@
    \slope@
    \ifNESW@\global\advance\angcount@\exacount@\else
     \global\advance\angcount@-\exacount@\fi
    \ifnum\angcount@>23 \global\angcount@23 \fi
    \ifnum\angcount@<\@ne\global\angcount@\@ne\fi
    \slope@a\angcount@
    \ifY@
     \advance\secondy@\Ydimen@
    \else
     \ifX@
      \advance\secondx@\Xdimen@
      \dimen@\secondx@\advance\dimen@-\firstx@
      \ifNESW@\else\multiply\dimen@\m@ne\fi
      \multiply\dimen@\tan@i\divide\dimen@\tan@ii
      \advance\dimen@\firsty@\secondy@\dimen@
     \fi
    \fi
   \else
    \ifH@\else\ifV@\else\slope@\fi\fi
   \fi
  \fi
 \fi
 \ifH@\else\ifV@\else\ifsvertex@\else
  \dimen@6\p@\multiply\dimen@\tan@ii
  \count@\tan@i\advance\count@\tan@ii\divide\dimen@\count@
  \ifE@\advance\firstx@\dimen@\else\advance\firstx@-\dimen@\fi
  \multiply\dimen@\tan@i\divide\dimen@\tan@ii
  \ifN@\advance\firsty@\dimen@\else\advance\firsty@-\dimen@\fi
 \fi\fi\fi
 \ifp@
  \ifH@\else\ifV@\else
   \getcos@\pdimen@\advance\firsty@\dimen@\advance\secondy@\dimen@
   \ifNESW@\advance\firstx@-\dimen@ii\else\advance\firstx@\dimen@ii\fi
  \fi\fi
 \fi
 \ifH@\else\ifV@\else
  \ifnum\tan@i>\tan@ii
   \charht@\ten@\p@\charwd@\ten@\p@
   \multiply\charwd@\tan@ii\divide\charwd@\tan@i
  \else
   \charwd@\ten@\p@\charht@\ten@\p@
   \divide\charht@\tan@ii\multiply\charht@\tan@i
  \fi
  \ifnum\tcount@=\thr@@
   \ifN@\advance\secondy@-.3\charht@\else\advance\secondy@.3\charht@\fi
  \fi
  \ifnum\scount@=\tw@
   \ifE@\advance\firstx@.3\charht@\else\advance\firstx@-.3\charht@\fi
  \fi
  \ifnum\tcount@=12
   \ifN@\advance\secondy@-\charht@\else\advance\secondy@\charht@\fi
  \fi
  \iftY@
  \else
   \ifa@
    \ifX@
    \else
     \secondx@\secondy@\advance\secondx@-\firsty@
     \ifNESW@\else\multiply\secondx@\m@ne\fi
     \multiply\secondx@\tan@ii\divide\secondx@\tan@i
     \advance\secondx@\firstx@
    \fi
   \fi
  \fi
 \fi\fi
 \ifH@\harrow@\else\ifV@\varrow@\else\arrow@@\fi\fi}
\newdimen\mathaxis@
\mathaxis@90\p@\divide\mathaxis@36
\def\harrow@b{\ifE@\hskip\tocenter@\hskip\firstx@\fi}
\def\harrow@bb{\ifE@\hskip\xdimen@\else\hskip\Xdimen@\fi}
\def\harrow@e{\ifE@\else\hskip-\firstx@\hskip-\tocenter@\fi}
\def\harrow@ee{\ifE@\hskip-\Xdimen@\else\hskip-\xdimen@\fi}
\def\harrow@{\dimen@\secondx@\advance\dimen@-\firstx@
 \ifE@\let\next@\rlap\else\multiply\dimen@\m@ne\let\next@\llap\fi
 \next@{%
  \harrow@b
  \smash{\raise\pdimen@\hbox to\dimen@
   {\harrow@bb\arrow@ii
    \ifnum\arrcount@=\m@ne\else\ifnum\arrcount@=\thr@@\else
     \ifE@
      \ifnum\scount@=\m@ne
      \else
       \ifcase\scount@\or\or\char118 \or\char117 \or\or\or\char119 \or
       \char120 \or\char121 \or\char122 \or\or\or\arrow@i\char125 \or
       \char117 \hskip\thr@@\p@\char117 \hskip-\thr@@\p@\fi
      \fi
     \else
      \ifnum\tcount@=\m@ne
      \else
       \ifcase\tcount@\char117 \or\or\char117 \or\char118 \or\char119 \or
       \char120 \or\or\or\or\or\char121 \or\char122 \or\arrow@i\char125
       \or\char117 \hskip\thr@@\p@\char117 \hskip-\thr@@\p@\fi
      \fi
     \fi
    \fi\fi
    \dimen@\mathaxis@\advance\dimen@.2\p@
    \dimen@ii\mathaxis@\advance\dimen@ii-.2\p@
    \ifnum\arrcount@=\m@ne
     \let\leads@\null
    \else
     \ifcase\arrcount@
      \def\leads@{\hrule\height\dimen@\depth-\dimen@ii}\or
      \def\leads@{\hrule\height\dimen@\depth-\dimen@ii}\or
      \def\leads@{\hbox to\ten@\p@{%
       \leaders\hrule\height\dimen@\depth-\dimen@ii\hfil
       \hfil
      \leaders\hrule\height\dimen@\depth-\dimen@ii\hskip\z@ plus2fil\relax
       \hfil
       \leaders\hrule\height\dimen@\depth-\dimen@ii\hfil}}\or
     \def\leads@{\hbox{\hbox to\ten@\p@{\dimen@\mathaxis@\advance\dimen@1.2\p@
       \dimen@ii\dimen@\advance\dimen@ii-.4\p@
       \leaders\hrule\height\dimen@\depth-\dimen@ii\hfil}%
       \kern-\ten@\p@
       \hbox to\ten@\p@{\dimen@\mathaxis@\advance\dimen@-1.2\p@
       \dimen@ii\dimen@\advance\dimen@ii-.4\p@
       \leaders\hrule\height\dimen@\depth-\dimen@ii\hfil}}}\fi
    \fi
    \cleaders\leads@\hfil
    \ifnum\arrcount@=\m@ne\else\ifnum\arrcount@=\thr@@\else
     \arrow@i
     \ifE@
      \ifnum\tcount@=\m@ne
      \else
       \ifcase\tcount@\char119 \or\or\char119 \or\char120 \or\char121 \or
       \char122 \or \or\or\or\or\char123 \or\char124 \or
       \char125 \or\char119 \hskip-\thr@@\p@\char119 \hskip\thr@@\p@\fi
      \fi
     \else
      \ifcase\scount@\or\or\char120 \or\char119 \or\or\or\char121 \or\char122
      \or\char123 \or\char124 \or\or\or\char125 \or
      \char119 \hskip-\thr@@\p@\char119 \hskip\thr@@\p@\fi
     \fi
    \fi\fi
    \harrow@ee}}%
  \harrow@e}%
 \iflabel@i
  \dimen@ii\z@\setbox\ZER@\hbox{$\m@th\tsize@@\label@i$}%
  \ifnum\arrcount@=\m@ne
  \else
   \advance\dimen@ii\mathaxis@
   \advance\dimen@ii\dp\ZER@\advance\dimen@ii\tw@\p@
   \ifnum\arrcount@=\thr@@\advance\dimen@ii\tw@\p@\fi
  \fi
  \advance\dimen@ii\pdimen@
  \next@{\harrow@b\smash{\raise\dimen@ii\hbox to\dimen@
   {\harrow@bb\hskip\tw@\ldimen@i\hfil\box\ZER@\hfil\harrow@ee}}\harrow@e}%
 \fi
 \iflabel@ii
  \ifnum\arrcount@=\m@ne
  \else
   \setbox\ZER@\hbox{$\m@th\tsize@\label@ii$}%
   \dimen@ii-\ht\ZER@\advance\dimen@ii-\tw@\p@
   \ifnum\arrcount@=\thr@@\advance\dimen@ii-\tw@\p@\fi
   \advance\dimen@ii\mathaxis@\advance\dimen@ii\pdimen@
   \next@{\harrow@b\smash{\raise\dimen@ii\hbox to\dimen@
    {\harrow@bb\hskip\tw@\ldimen@ii\hfil\box\ZER@\hfil\harrow@ee}}\harrow@e}%
  \fi
 \fi}
\let\tsize@\tsize
\def\tsizeCDlabels{\let\tsize@\tsize}
\def\ssizeCDlabels{\let\tsize@\ssize}
\def\tsize@@{\ifnum\arrcount@=\m@ne\else\tsize@\fi}
\def\varrow@{\dimen@\secondy@\advance\dimen@-\firsty@
 \ifN@\else\multiply\dimen@\m@ne\fi
 \setbox\ZER@\vbox to\dimen@
  {\ifN@\vskip-\Ydimen@\else\vskip\ydimen@\fi
   \ifnum\arrcount@=\m@ne\else\ifnum\arrcount@=\thr@@\else
    \hbox{\arrow@iii
     \ifN@
      \ifnum\tcount@=\m@ne
      \else
       \ifcase\tcount@\char117 \or\or\char117 \or\char118 \or\char119 \or
       \char120 \or\or\or\or\or\char121 \or\char122 \or\char123 \or
       \vbox{\hbox{\char117}\nointerlineskip\vskip\thr@@\p@
       \hbox{\char117}\vskip-\thr@@\p@}\fi
      \fi
     \else
      \ifcase\scount@\or\or\char118 \or\char117 \or\or\or\char119 \or
      \char120 \or\char121 \or\char122 \or\or\or\char123 \or
      \vbox{\hbox{\char117}\nointerlineskip\vskip\thr@@\p@
      \hbox{\char117}\vskip-\thr@@\p@}\fi
     \fi}%
    \nointerlineskip
   \fi\fi
   \ifnum\arrcount@=\m@ne
    \let\leads@\null
   \else
    \ifcase\arrcount@\let\leads@\vrule\or\let\leads@\vrule\or
    \def\leads@{\vbox to\ten@\p@{%
     \hrule\height1.67\p@\depth\z@\width.4\p@
     \vfil
     \hrule\height3.33\p@\depth\z@\width.4\p@
     \vfil
     \hrule\height1.67\p@\depth\z@\width.4\p@}}\or
    \def\leads@{\hbox{\vrule\height\p@\hskip\tw@\p@\vrule}}\fi
   \fi
  \cleaders\leads@\vfill\nointerlineskip
   \ifnum\arrcount@=\m@ne\else\ifnum\arrcount@=\thr@@\else
    \hbox{\arrow@iv
     \ifN@
      \ifcase\scount@\or\or\char118 \or\char117 \or\or\or\char119 \or
      \char120 \or\char121 \or\char122 \or\or\or\arrow@iii\char123 \or
      \vbox{\hbox{\char117}\nointerlineskip\vskip-\thr@@\p@
      \hbox{\char117}\vskip\thr@@\p@}\fi
     \else
      \ifnum\tcount@=\m@ne
      \else
       \ifcase\tcount@\char117 \or\or\char117 \or\char118 \or\char119 \or
       \char120 \or\or\or\or\or\char121 \or\char122 \or\arrow@iii\char123 \or
       \vbox{\hbox{\char117}\nointerlineskip\vskip-\thr@@\p@
       \hbox{\char117}\vskip\thr@@\p@}\fi
      \fi
     \fi}%
   \fi\fi
   \ifN@\vskip\ydimen@\else\vskip-\Ydimen@\fi}%
 \ifN@
  \dimen@ii\firsty@
 \else
  \dimen@ii-\firsty@\advance\dimen@ii\ht\ZER@\multiply\dimen@ii\m@ne
 \fi
 \rlap{\smash{\hskip\tocenter@\hskip\pdimen@\raise\dimen@ii\box\ZER@}}%
 \iflabel@i
  \setbox\ZER@\vbox to\dimen@{\vfil
   \hbox{$\m@th\tsize@@\label@i$}\vskip\tw@\ldimen@i\vfil}%
  \rlap{\smash{\hskip\tocenter@\hskip\pdimen@
  \ifnum\arrcount@=\m@ne\let\next@\relax\else\let\next@\llap\fi
  \next@{\raise\dimen@ii\hbox{\ifnum\arrcount@=\m@ne\hskip-.5\wd\ZER@\fi
   \box\ZER@\ifnum\arrcount@=\m@ne\else\hskip\tw@\p@\fi}}}}%
 \fi
 \iflabel@ii
  \ifnum\arrcount@=\m@ne
  \else
   \setbox\ZER@\vbox to\dimen@{\vfil
    \hbox{$\m@th\tsize@\label@ii$}\vskip\tw@\ldimen@ii\vfil}%
   \rlap{\smash{\hskip\tocenter@\hskip\pdimen@
   \rlap{\raise\dimen@ii\hbox{\ifnum\arrcount@=\thr@@\hskip4.5\p@\else
    \hskip2.5\p@\fi\box\ZER@}}}}%
  \fi
 \fi
}
\newdimen\goal@
\newdimen\shifted@
\newcount\Tcount@
\newcount\Scount@
\newbox\shaft@
\newcount\slcount@
\def\getcos@#1{%
 \ifnum\tan@i<\tan@ii
  \dimen@#1%
  \ifnum\slcount@<8 \count@9 \else \ifnum\slcount@<12 \count@8 \else
   \count@7 \fi\fi
  \multiply\dimen@\count@\divide\dimen@\ten@
  \dimen@ii\dimen@\multiply\dimen@ii\tan@i\divide\dimen@ii\tan@ii
 \else
  \dimen@ii#1%
  \count@-\slcount@\advance\count@24
  \ifnum\count@<8 \count@9 \else \ifnum\count@<12 \count@8
   \else\count@7 \fi\fi
  \multiply\dimen@ii\count@\divide\dimen@ii\ten@
  \dimen@\dimen@ii\multiply\dimen@\tan@ii\divide\dimen@\tan@i
 \fi}
\newdimen\adjust@
\def\Nnext@{\ifN@\let\next@\raise\else\let\next@\lower\fi}
\def\arrow@@{\slcount@\angcount@
 \ifNESW@
  \ifnum\angcount@<\ten@
   \let\arrowfont@\arrow@i\global\advance\angcount@\m@ne
   \global\multiply\angcount@13
  \else
   \ifnum\angcount@<19
    \let\arrowfont@\arrow@ii\global\advance\angcount@-\ten@
    \global\multiply\angcount@13
   \else
    \let\arrowfont@\arrow@iii\global\advance\angcount@-19
    \global\multiply\angcount@13
  \fi\fi
  \Tcount@\angcount@
 \else
  \ifnum\angcount@<5
   \let\arrowfont@\arrow@iii\global\advance\angcount@\m@ne
   \global\multiply\angcount@13 \global\advance\angcount@65
  \else
   \ifnum\angcount@<14
    \let\arrowfont@\arrow@iv\global\advance\angcount@-5
    \global\multiply\angcount@13
   \else
    \ifnum\angcount@<23
     \let\arrowfont@\arrow@v\global\advance\angcount@-14
     \global\multiply\angcount@13
    \else
     \let\arrowfont@\arrow@i\global\angcount@117
  \fi\fi\fi
  \ifnum\angcount@=117 \Tcount@115 \else\Tcount@\angcount@\fi
 \fi
 \Scount@\Tcount@
 \ifE@
  \ifnum\tcount@=\z@\advance\Tcount@\tw@\else\ifnum\tcount@=13
   \advance\Tcount@\tw@\else\advance\Tcount@\tcount@\fi\fi
  \ifnum\scount@=\z@\else\ifnum\scount@=13 \advance\Scount@\thr@@\else
   \advance\Scount@\scount@\fi\fi
 \else
  \ifcase\tcount@\advance\Tcount@\thr@@\or\or\advance\Tcount@\thr@@\or
  \advance\Tcount@\tw@\or\advance\Tcount@6 \or\advance\Tcount@7
  \or\or\or\or\or\advance\Tcount@8 \or\advance\Tcount@9 \or
  \advance\Tcount@12 \or\advance\Tcount@\thr@@\fi
  \ifcase\scount@\or\or\advance\Scount@\thr@@\or\advance\Scount@\tw@\or
  \or\or\advance\Scount@4 \or\advance\Scount@5 \or\advance\Scount@\ten@
  \or\advance\Scount@11 \or\or\or\advance\Scount@12 \or\advance
  \Scount@\tw@\fi
 \fi
 \ifcase\arrcount@\or\or\global\advance\angcount@\@ne\else\fi
 \ifN@\shifted@\firsty@\else\shifted@-\firsty@\fi
 \ifE@\else\advance\shifted@\charht@\fi
 \goal@\secondy@\advance\goal@-\firsty@
 \ifN@\else\multiply\goal@\m@ne\fi
 \setbox\shaft@\hbox{\arrowfont@\char\angcount@}%
 \ifnum\arrcount@=\thr@@
  \getcos@{1.5\p@}%
  \setbox\shaft@\hbox to\wd\shaft@{\arrowfont@
   \rlap{\hskip\dimen@ii
    \smash{\ifNESW@\let\next@\lower\else\let\next@\raise\fi
     \next@\dimen@\hbox{\arrowfont@\char\angcount@}}}%
   \rlap{\hskip-\dimen@ii
    \smash{\ifNESW@\let\next@\raise\else\let\next@\lower\fi
      \next@\dimen@\hbox{\arrowfont@\char\angcount@}}}\hfil}%
 \fi
 \rlap{\smash{\hskip\tocenter@\hskip\firstx@
  \ifnum\arrcount@=\m@ne
  \else
   \ifnum\arrcount@=\thr@@
   \else
    \ifnum\scount@=\m@ne
    \else
     \ifnum\scount@=\z@
     \else
      \setbox\ZER@\hbox{\ifnum\angcount@=117 \arrow@v\else\arrowfont@\fi
       \char\Scount@}%
      \ifNESW@
       \ifnum\scount@=\tw@
        \dimen@\shifted@\advance\dimen@-\charht@
        \ifN@\hskip-\wd\ZER@\fi
        \Nnext@
        \next@\dimen@\copy\ZER@
        \ifN@\else\hskip-\wd\ZER@\fi
       \else
        \Nnext@
        \ifN@\else\hskip-\wd\ZER@\fi
        \next@\shifted@\copy\ZER@
        \ifN@\hskip-\wd\ZER@\fi
       \fi
       \ifnum\scount@=12
        \advance\shifted@\charht@\advance\goal@-\charht@
        \ifN@\hskip\wd\ZER@\else\hskip-\wd\ZER@\fi
       \fi
       \ifnum\scount@=13
        \getcos@{\thr@@\p@}%
        \ifN@\hskip\dimen@\else\hskip-\wd\ZER@\hskip-\dimen@\fi
        \adjust@\shifted@\advance\adjust@\dimen@ii
        \Nnext@
        \next@\adjust@\copy\ZER@
        \ifN@\hskip-\dimen@\hskip-\wd\ZER@\else\hskip\dimen@\fi
       \fi
      \else
       \ifN@\hskip-\wd\ZER@\fi
       \ifnum\scount@=\tw@
        \ifN@\hskip\wd\ZER@\else\hskip-\wd\ZER@\fi
        \dimen@\shifted@\advance\dimen@-\charht@
        \Nnext@
        \next@\dimen@\copy\ZER@
        \ifN@\hskip-\wd\ZER@\fi
       \else
        \Nnext@
        \next@\shifted@\copy\ZER@
        \ifN@\else\hskip-\wd\ZER@\fi
       \fi
       \ifnum\scount@=12
        \advance\shifted@\charht@\advance\goal@-\charht@
        \ifN@\hskip-\wd\ZER@\else\hskip\wd\ZER@\fi
       \fi
       \ifnum\scount@=13
        \getcos@{\thr@@\p@}%
        \ifN@\hskip-\wd\ZER@\hskip-\dimen@\else\hskip\dimen@\fi
        \adjust@\shifted@\advance\adjust@\dimen@ii
        \Nnext@
        \next@\adjust@\copy\ZER@
        \ifN@\hskip\dimen@\else\hskip-\dimen@\hskip-\wd\ZER@\fi
       \fi	
      \fi
  \fi\fi\fi\fi
  \ifnum\arrcount@=\m@ne
  \else
   \loop
    \ifdim\goal@>\charht@
    \ifE@\else\hskip-\charwd@\fi
    \Nnext@
    \next@\shifted@\copy\shaft@
    \ifE@\else\hskip-\charwd@\fi
    \advance\shifted@\charht@\advance\goal@-\charht@
    \repeat
   \ifdim\goal@>\z@
    \dimen@\charht@\advance\dimen@-\goal@
    \divide\dimen@\tan@i\multiply\dimen@\tan@ii
    \ifE@\hskip-\dimen@\else\hskip-\charwd@\hskip\dimen@\fi
    \adjust@\shifted@\advance\adjust@-\charht@\advance\adjust@\goal@
    \Nnext@
    \next@\adjust@\copy\shaft@
    \ifE@\else\hskip-\charwd@\fi
   \else
    \adjust@\shifted@\advance\adjust@-\charht@
   \fi
  \fi
  \ifnum\arrcount@=\m@ne
  \else
   \ifnum\arrcount@=\thr@@
   \else
    \ifnum\tcount@=\m@ne
    \else
     \setbox\ZER@
      \hbox{\ifnum\angcount@=117 \arrow@v\else\arrowfont@\fi\char\Tcount@}%
     \ifnum\tcount@=\thr@@
      \advance\adjust@\charht@
      \ifE@\else\ifN@\hskip-\charwd@\else\hskip-\wd\ZER@\fi\fi
     \else
      \ifnum\tcount@=12
       \advance\adjust@\charht@
       \ifE@\else\ifN@\hskip-\charwd@\else\hskip-\wd\ZER@\fi\fi
      \else
       \ifE@\hskip-\wd\ZER@\fi
     \fi\fi
     \Nnext@
     \next@\adjust@\copy\ZER@
     \ifnum\tcount@=13
      \hskip-\wd\ZER@
      \getcos@{\thr@@\p@}%
      \ifE@\hskip-\dimen@\else\hskip\dimen@\fi
      \advance\adjust@-\dimen@ii
      \Nnext@
      \next@\adjust@\box\ZER@
     \fi
  \fi\fi\fi}}%
 \iflabel@i
  \rlap{\hskip\tocenter@
  \dimen@\firstx@\advance\dimen@\secondx@\divide\dimen@\tw@
  \advance\dimen@\ldimen@i
  \dimen@ii\firsty@\advance\dimen@ii\secondy@\divide\dimen@ii\tw@
  \global\multiply\ldimen@i\tan@i\global\divide\ldimen@i\tan@ii
  \ifNESW@\advance\dimen@ii\ldimen@i\else\advance\dimen@ii-\ldimen@i\fi
  \setbox\ZER@\hbox{\ifNESW@\else\ifnum\arrcount@=\thr@@\hskip4\p@\else
   \hskip\tw@\p@\fi\fi
   $\m@th\tsize@@\label@i$\ifNESW@\ifnum\arrcount@=\thr@@\hskip4\p@\else
   \hskip\tw@\p@\fi\fi}%
  \ifnum\arrcount@=\m@ne
   \ifNESW@\advance\dimen@.5\wd\ZER@\advance\dimen@\p@\else
    \advance\dimen@-.5\wd\ZER@\advance\dimen@-\p@\fi
   \advance\dimen@ii-.5\ht\ZER@
  \else
   \advance\dimen@ii\dp\ZER@
   \ifnum\slcount@<6 \advance\dimen@ii\tw@\p@\fi
  \fi
  \hskip\dimen@
  \ifNESW@\let\next@\llap\else\let\next@\rlap\fi
  \next@{\smash{\raise\dimen@ii\box\ZER@}}}%
 \fi
 \iflabel@ii
  \ifnum\arrcount@=\m@ne
  \else
   \rlap{\hskip\tocenter@
   \dimen@\firstx@\advance\dimen@\secondx@\divide\dimen@\tw@
   \ifNESW@\advance\dimen@\ldimen@ii\else\advance\dimen@-\ldimen@ii\fi
   \dimen@ii\firsty@\advance\dimen@ii\secondy@\divide\dimen@ii\tw@
   \global\multiply\ldimen@ii\tan@i\global\divide\ldimen@ii\tan@ii
   \advance\dimen@ii\ldimen@ii
   \setbox\ZER@\hbox{\ifNESW@\ifnum\arrcount@=\thr@@\hskip4\p@\else
    \hskip\tw@\p@\fi\fi
    $\m@th\tsize@\label@ii$\ifNESW@\else\ifnum\arrcount@=\thr@@\hskip4\p@
    \else\hskip\tw@\p@\fi\fi}%
   \advance\dimen@ii-\ht\ZER@
   \ifnum\slcount@<9 \advance\dimen@ii-\thr@@\p@\fi
   \ifNESW@\let\next@\rlap\else\let\next@\llap\fi
   \hskip\dimen@\next@{\smash{\raise\dimen@ii\box\ZER@}}}%
  \fi
 \fi
}
\def\outCD@#1{\def#1{\Err@{\noexpand#1must not be used within \string\CD}}}
\newskip\preCDskip@
\newskip\postCDskip@
\preCDskip@\z@
\postCDskip@\z@
\def\preCDspace#1{\RIfMIfI@
 \onlydmatherr@\preCDspace\else\advance\preCDskip@#1\relax\fi\else
 \onlydmatherr@\preCDspace\fi}
\def\postCDspace#1{\RIfMIfI@
 \onlydmatherr@\postCDspace\else\advance\postCDskip@#1\relax\fi\else
 \onlydmatherr@\postCDspace\fi}
\def\predisplayspace#1{\RIfMIfI@
 \onlydmatherr@\predisplayspace\else
 \advance\abovedisplayskip#1\relax
 \advance\abovedisplayshortskip#1\relax\fi
 \else\onlydmatherr@\preCDspace\fi}
\def\postdisplayspace#1{\RIfMIfI@
 \onlydmatherr@\postdisplayspace\else
 \advance\belowdisplayskip#1\relax
 \advance\belowdisplayshortskip#1\relax\fi
 \else\onlydmatherr@\postdisplayspace\fi}
\def\PreCDSpace#1{\global\preCDskip@#1\relax}
\def\PostCDSpace#1{\global\postCDskip@#1\relax}
\def\CD#1\endCD{%
 \outCD@\cgaps\outCD@\rgaps\outCD@\Cgaps\outCD@\Rgaps
 \preCD@#1\endCD
 \advance\abovedisplayskip\preCDskip@
 \advance\abovedisplayshortskip\preCDskip@
 \advance\belowdisplayskip\postCDskip@
 \advance\belowdisplayshortskip\postCDskip@
 \vcenter{\offinterlineskip
  \vskip\preCDskip@\Let@\global\colcount@\@ne\global\rowcount@\z@
  \everycr{%
   \noalign{%
    \ifnum\rowcount@=\Rowcount@
    \else
     \getrgap@\rowcount@\vskip\getdim@
     \global\advance\rowcount@\@ne\global\colcount@\@ne
    \fi}}%
  \tabskip\z@
  \halign{&\global\xoff@\z@\global\yoff@\z@
   \getcgap@\colcount@\hskip\getdim@
   \hfil\vrule\height\ten@\p@\width\z@\depth\z@
   $\m@th\displaystyle{##}$\hfil
   \global\advance\colcount@\@ne\cr
   #1\crcr}\vskip\postCDskip@}%
 \preCDskip@\z@\postCDskip@\z@
 \def\getcgap@##1{\ifcase##1\or\getdim@\z@\else\getdim@\standardcgap\fi}%
 \def\getrgap@##1{\ifcase##1\getdim@\z@\else\getdim@\standardrgap\fi}%
 \let\Width@\relax\let\Height@\relax\let\Depth@\relax\let\Rowheight@\relax
 \let\Rowdepth@\relax\let\Colwidth@\relax
}
\let\enddocument\bye
\def\alloc@#1#2#3#4#5{\global\advance\count1#1by\@ne
  \ch@ck#1#4#2%
  \allocationnumber=\count1#1%
  \global#3#5=\allocationnumber
  \wlog{\string#5=\string#2\the\allocationnumber}}
\catcode`\@=\active

\loadmsam
\loadmsbm
\loadeufm

\litdelimiter*
\litbackslash"

\input paper.st\relax

\def\Gm{{\Bbb G}_{\text{m}}}

\def\II{{\Bbb I}}

\def\QQ{{\Bbb Q}}

\def\ZZ{{\Bbb Z}}

\def\into{\hookrightarrow}
\newsymbol\onto 1310
\def\iso{\simeq}
\def\to{\rightarrow}

\def\bar{\overline}

\def\invlim{\varprojlim}
\def\dirlim{\varinjlim}
\def\tensor{\otimes}

\def\Ž{\'e}
\def\{\`e}
\def\ˆ{\`a}
\def\{\c c}
\def\"{\^\i}
\def\Ÿ{\"u}
\def\Š{\"a}

\newsymbol\nmid 232D

\def\divides{\mid}

\def\Ext{\operatorname{Ext}}

\def\Gal{\operatorname{Gal}}

\def\Hom{\operatorname{Hom}}

\def\Pic{\operatorname{Pic}}

\def\ord{\operatorname{ord}}

\def\res{\operatorname{res}}
\def\cores{\operatorname{cores}}

\def\Spec{\operatorname{Spec}}
\def\rank{\operatorname{rank}}

\catcode`!=11

\define\!dlog[#1]{\operatorname{dlog}^{(#1)}}
\catcode`!=12

\newclaim\Theorem\c{}{Theorem}
\newfontstyle\Theorem{\bf}
\shortenclaim\Theorem\theorem
\newclaim\Lemma\c{}{Lemma}
\newfontstyle\Lemma{\bf}
\shortenclaim\Lemma\lemma
\newclaim\Proposition\c{}{Proposition}
\newfontstyle\Proposition{\bf}
\shortenclaim\Proposition\proposition
\newclaim\Corollary\c{}{Corollary}
\newfontstyle\Corollary{\bf}
\shortenclaim\Corollary\corollary
\newclaim\Definition\c{}{Definition}
\shortenclaim\Definition\definition
\newfontstyle\Definition{\bf}
\define\proof{\demo{\it Proof}}
\define\endproof{\qed\enddemo}

\def\notation#1{\noindent{\it Notation\/}.
\smallskip
\halign\bgroup##\hfil&\quad
\vtop{\setbox0=\hbox{#1\quad}
\dimen0=\hsize
\advance\dimen0 by -\wd0
\parindent=0pt\hsize=\dimen0\hangindent1em\strut##\strut}\cr
}
\def\endnotation{\egroup\medskip}

\def\display#1{
\smallskip
\halign\bgroup##\hfil&\quad
\vtop{\setbox0=\hbox{#1\quad}
\dimen0=\hsize
\advance\dimen0 by -\wd0
\parindent=0pt\hsize=\dimen0\hangindent1em\strut##\strut}\cr
}
\def\enddisplay{\egroup\medskip}

\catcode`\@=11

\ifx\undefined\square
  \def\square{\vrule width.6em height.5em depth.1em\relax}\fi
\def\qed{\ifhmode\unskip\nobreak\fi\quad
  \ifmmode\square\else$\m@th\square$\fi}

\def\widestkey#1{\rm\setboxz@h{[{\bf #1}]\ }\bibindent@\wdz@}

\catcode`\@=\active

\catcode`!=11

\def\!ifnextchar#1#2#3{%
  \let\!testchar=#1%
  \def\!first{#2}%
  \def\!second{#3}%
  \futurelet\!nextchar\!testnext}
\def\!testnext{%
  \ifx \!nextchar \!spacetoken 
    \let\!next=\!skipspacetestagain
  \else
    \ifx \!nextchar \!testchar
      \let\!next=\!first
    \else 
      \let\!next=\!second 
    \fi 
  \fi
  \!next}
\def\\{\!skipspacetestagain} 
  \expandafter\def\\ {\futurelet\!nextchar\!testnext} 
\def\\{\let\!spacetoken= } \\
 
\catcode`!=12

\catcode`!=11

\define\N{\!ifnextchar[{\!N}{\!N[\null]}}
\define\!N[#1]{@()\L{#1}@(0,1)}

\define\!E[#1]{@()\L{#1}@(1,0)}
\predefine\Ssymbol{\S}
\redefine\S{\!ifnextchar[{\!S}{\!S[\null]}}
\define\!S[#1]{@()\L{#1}@(0,-1)}

\define\!W[#1]{@()\L{#1}@(-1,0)}

\define\!NE[#1]{@()\L{#1}@(1,1)}

\define\!SE[#1]{@()\L{#1}@(1,-1)}

\define\!NW[#1]{@()\L{#1}@(-1,1)}

\define\!SW[#1]{@()\L{#1}@(-1,-1)}

\catcode`!=12

\catcode`\"=\active

\document

\newclaim\Conjecture\c{}{Conjecture}
\shortenclaim\Conjecture\conjecture
\newfontstyle\Conjecture{\bf}

\def\tor{{\hbox{\rm tor}}}
\def\p{{\frak p}}
\def\II{{\bold I}}
\def\m{\frak m}
\def\a{\frak a}
\def\seq{{\bold x}}

\def\K#1{K_{#1}}
\def\Free#1{\Lambda^{#1}}
\def\NN{{\Bbb N}}
\def\Un{E_{F}^{u}}
\def\No{E_{F}^{n}}
\def\Lo{E_{F}^{\text{loc}}}

\title \smc 
GREENBERG'S CONJECTURE AND \\
\smc UNITS IN MULTIPLE $\ZZ_p$-EXTENSIONS
\endtitle

\author William G. McCallum \endauthor

\abstract 
Let $A$ be the inverse limit of the $p$-part of the ideal class 
groups in a $\ZZ_{p}^{r}$-extension $K_{\infty}/K$.
Greenberg conjectures that if 
$r$ is maximal, then $A$ is pseudo-null as a module over the Iwasawa 
algebra $\Lambda$ (that is, has codimension at least 2). 
We prove this conjecture in the case that $K$ is the 
field of $p$-th roots of unity, $p$ has index of irregularity 1, 
satisfies Vandiver's conjecture, and satisfies a mild additional 
 hypothesis on units. We also show that if $K$ is the field of $p$-th 
roots of unity and $r$ is maximal, Greenberg's conjecture for 
$K$
implies that the maximal $p$-ramified pro-$p$-extension of $K$ cannot have a free pro-$p$ 
quotient of rank $r$ unless $p$ is regular (see also \cite{LN}). 
Finally, we prove a generalization  of a theorem of 
Iwasawa in the case $r=1$ concerning the Kummer extension of 
$K_{\infty}$ generated by $p$-power roots of $p$-units. We show that 
the Galois group of this extension is torsion-free as a 
$\Lambda$-module if there is only one prime of $K$ above $p$ and 
$K_{\infty}$ contains all the $p$-power roots of unity.

\endabstract

\date{June, 2000}

\thanks{This research was supported in part by an AMS
Bicentennial Fellowship and a National Science Foundation grant
DMS-9624219.}

\address{Department of Mathematics, University of Arizona,
Tucson, AZ, 85721}

\maketitle

Let $K$ be a number field,  let $p$ be an odd prime number, and let 
$A(K)$ be the $p$-Sylow subgroup
of the ideal class group of $K$. In 1956, Iwasawa introduced
the idea of studying the behaviour of $A(F)$
as $F$ varies over  all intermediate fields in a $\ZZ_p$-extension 
$K_\infty/K$. Greenberg \cite{G2} gives a comprehensive account of the 
subsequent 
development of Iwasawa theory; we recall here some of the basic ideas.
The inverse limit  $A = \invlim_F A(F)$ is a finitely generated torsion
module over the Iwasawa algebra $\Lambda =
\ZZ_p[[\Gal(K_\infty/K)]] = \invlim_F  \ZZ_p[\Gal(F/K)]$. This algebra
has a particularly simple structure: given a 
topological
generator $\gamma$ of $\Gamma$, there is an isomorphism $\Lambda \iso 
\ZZ_{p}[[T]]$ taking $1+ \gamma$ to $T$. Iwasawa showed that 
there is a {\it pseudo-isomorphism}
(a map with finite kernel and cokernel) from $A$ to an elementary module $E =
\sum_i \Lambda/(f_i)$; the power series $f = \prod_i f_i$ is the {\it characteristic
power series} of $A$. In \cite{G1}, Greenberg conjectured that $f=1$ when $K$ 
is totally
real.

It is natural to consider the generalization of these ideas to a
compositum of
$\ZZ_p$-extensions of
$K$, that is, an extension
$K_\infty/K$  whose galois group
$\Gamma = \Gal(K_\infty/K)$ is isomorphic to $\ZZ_p^r$ for some positive
integer
$r$. The corresponding Iwasawa algebra, $\Lambda = \invlim_F
\ZZ_p[\Gal(F/K)]$, 
is isomorphic to a power series ring in $r$
variables over $\ZZ_p$.
If $K$ is totally real and satisfies Leopoldt's conjecture, then the 
only $\ZZ_{p}$-extension is the cyclotomic one
and Greenberg's  conjecture that $f=1$ in this case is equivalent to the statement 
that $A$ is finite, and hence has annihilator 
of height 2.
More generally, Greenberg has made the following conjecture \cite{G2}.

\claim""{Conjecture}{(Greenberg)}
Let $K_\infty/K$ be the compositum of all $\ZZ_p$-extensions of $K$. Then
the annilator of $A$  
has height at least two.
\endclaim

A module whose annihilator has height at least 2 is said to be {\it 
pseudo-null}. 
Let $E = {\Cal O}_{K}^{\times}$ and $U =\prod_{\p\divides 
p}{\Cal O}_{K_\p}^{\times}$, and denote by $\overline E$ the closure 
of $E$ in $U$. The following theorem is a consequence of the main 
theorem of this paper, Theorem~\ref{ReallyMainTheorem}.

\theorem\label{MainTheorem}
Let $K = \QQ(e^{2\pi i/p})$. Suppose that 
\list
\item  $A(K) \iso \ZZ/p\ZZ$
\item $(U/\overline E)[p^{\infty}] \iso \ZZ/p\ZZ$.
\endlist
Then
$K$ satisfies Greenberg's conjecture.
\endtheorem

We note that Greenberg's conjecture is trivially true for regular 
primes, since in that case $A = 0$. Condition (1) of the theorem 
implies that $p$ satisfies Vandiver's conjecture, by the reflection 
principle.
The conditions of the theorem are satisfied often (heuristically, by about 
3/4 of irregular primes), for example by  37, 59, 67, 101, 
103, 131, and 149. More specifically, 
the conditions are equivalent to the following:
\list 
\item $p$ satisfies Vandiver's conjecture, and divides exactly one of the Bernoulli numbers $B_2,B_4, 
\dots , B_{p-1}$
\item if we write the characteristic  series of $A$ for the cyclotomic
$\ZZ_{p}$-extension in the form 
$f = (T-cp)u$, where $u$ is a unit power series  and  $\gamma$ is chosen
to satisfy  $\zeta^{\gamma} = 
\zeta^{1+p}$ for all $p$-power roots of unity  $\zeta$, then $c \not\equiv 1 
\pmod{p}$.
\endlist
(See \cite{W} Theorem 10.16 
and Theorem 8.25.) The computations in \cite{BCEMS} verify 
Vandiver's conjecture for all primes 
less than 12,000,000, and condition (1) is 
satisfied by about 30\%
of those primes (about 61\%
of them are regular, in which case $A=0$). Iwasawa and 
Sims \cite{IS} tabulate the congruence classes mod $p$ of the 
$p$-adic integer $c$ for $1< p < 400$ and $3600 < p \le 4001$.
Condition (2) is satisfied for all primes in their  tables
that satisfy
condition (1).

Greenberg's conjecture may be regarded as a statement about the 
structure of  Galois groups.  We call an
extension of number fields {\it $p$-ramified} if it is unramified 
outside all the primes above $p$, and {\it $p$-split} if all the 
primes above $p$ split completely.
 Let
$$
\align
M_\infty &= \text{maximal abelian $p$-ramified pro-$p$-extension of $K_\infty$ }\\
L_\infty &= \text{maximal abelian unramified pro-$p$-extension of
$K_\infty$}\\
L'_\infty &= \text{maximal $p$-split subextension of
$L_\infty/K_\infty$}
\endalign
$$
Then  the three galois groups
$$
Y =
\Gal(M_\infty/K_\infty), \quad X = \Gal(L_\infty/K_\infty) \quad \text{and}
\quad X' =
\Gal(L'_\infty/K_\infty)\tag\label{StandardIwasawaModules}
$$
may be regarded as $\Lambda$-modules via the conjugation action of $\Gamma$.
From class field theory we have a canonical isomorphism of
$\Lambda$-modules
$
A \iso X
$. Thus Greenberg's conjecture 
asserts that $X$ is pseudo-null. This turns to be 
equivalent to the statement that $Y$ is $\Lambda$-torsion-free (see 
Theorem~\ref{YInterpretation} and Corollary~\ref{greenbergcorollary}),
and it is this latter formulation that we prove. The equivalence
was communicated to the author by Greenberg, and has been proved 
independently by Lannuzel and 
Nguyen Quang Do in \cite{LN}.

Greenberg's conjecture has a connection with the 
structure of ${\Cal G} = \Gal(\Omega/K)$, where $\Omega$ is the 
maximal $p$-ramified
pro-$p$ extension of $K$. This connection has been studied by 
Lannuzel and Nguyen Quang Do in \cite{LN}. 
It is well-known that ${\Cal G}$ is the 
quotient of a free pro-$p$-group on $g$ generators by $s$ relations, where 
$$
\align 
g &= \dim_{\ZZ/p\ZZ} H^1({\Cal G}, \ZZ/p\ZZ)\\
s &=\dim_{\ZZ/p\ZZ}H^2({\Cal G}, \ZZ/p\ZZ),
\endalign
$$
and that $g$ and $s$ are minimal in this regard. 
Thus ${\Cal G}$ cannot have a
free pro-$p$ quotient of rank greater than $g-s$.  On the other hand, ${\Cal G}$ 
is free of rank $g-s$ if $s=0$, that is, in the case 
$K=\QQ(\zeta_{p})$, if $p$ is regular.

\theorem\label{WingbergContradictionTheorem}
Suppose that $\QQ(\zeta_p)$ satisfies Greenberg's conjecture.
Then ${\Cal G}$ has a pro-$p$-free quotient
of rank $g-s$ if and only if $p$ is regular.
\endtheorem

This is proved in Section~\ref{YInterpretation}. Lannuzel and 
Nguyen Quang Do prove a similar result for a  general ground field $K$,
but under the more restrictive 
hypothesis that all finite abelian $p$-ramified extensions of $K$
satisfy Leopoldt's conjecture (\cite{LN}, Theorem 5.4.)

In proving that $Y$ is torsion-free, we will make use of 
the following  generalization of a theorem of 
Iwasawa on units \cite{Iw}. Let
$$
N_\infty = K_\infty(\epsilon^{1/p^{n}}: n \in {\Bbb N}, \epsilon 
\in {\Cal O}_{K_{\infty}}[1/p]^\times).
$$
\theorem\label{torsionfreetheorem}
Suppose that $K_\infty/K$ is a multiple $\ZZ_p$-extension satisfying:
\list
\item $K_\infty$
contains all the
$p^n$-th roots of unity,
$n>0$ 
\item there is only one prime of $K$ above
$p$.
\endlist
Then the $\Lambda$-module $Y' = \Gal(N_\infty/K_\infty)$ 
is torsion-free. In particular, $Y_{\tor}$ fixes $N_{\infty}$.
\endtheorem

\medskip

{\it Acknowledgements} The essential ideas of this work were developed in 
1995--96, at the Institute for Advanced Study. The author would like to thank
them for their hospitality,   Ralph Greenberg 
and David Marshall for useful comments and conversations, and the 
referee for suggesting numerous improvements. 

\section{Ext Groups and Local Cohomology\label{Extsection}}

In this section we collect together various results from the 
literature (see, for example, \cite{B-H} and \cite{H}).
Let $N = r+1$, the dimension of $\Lambda$. Let $\m$ be the maximal
ideal of $\Lambda$. 
Let ${\seq} = x_1, \dots , x_N$ be a 
finite sequence in $\Lambda$ that generates an $\m$-primary ideal. 
The Koszul complex 
$$
0 \East d{} \K N(\seq) \East d{} \K {N-1}(\seq) \East d{}
\cdots \East d{} \K 1(\seq)  \East d{}\K 0(\seq) \to 0
$$
is  defined as follows. Let  
$\Free N$ be the free $\Lambda$-module with basis 
$\{e_i : 1 \le i \le N\}$. Then
$$
\K i(\seq) = 
\bigwedge^{i} \Free N
$$
and 
$$
d(e_{j_1}\wedge \dots \wedge e_{j_i}) =
\sum_{s=1}^i(-1)^{s-1}x_{j_s} e_{j_1} \wedge \cdots \wedge \hat
e_{j_s} \wedge \cdots \wedge e_{j_i},
$$
where the hat indicates that a term is omitted. 
If $\seq$ is a regular sequence, then the Koszul complex is
a free resolution of $\Lambda/(\seq)$ (\cite{Ma}, Theorem
43).

Given a finitely generated $\Lambda$-module $X$, we define a complex
$\K
\null^*(\seq,X)$ by
$$
\K \null^*(\seq,X) = \Hom_\Lambda(\K *(\seq),X)
$$ 
and define
 $H^*(\seq,X)$ to be the cohomology of this complex. 

Now let $(\seq_n)$ be a 
sequence of sequences
$$\seq_n = (x_{n,1},\dots,x_{n,{N}}),$$ such that
 $x_{n,i} \divides x_{m,{i}}$, $n \le m$, $1 \le i \le N$.
If $n \le m$, there is a natural 
map of complexes
$$
\K .(\seq_n) \to \K .(\seq_m),
$$
which multiplies 
$e_{j_1}\wedge \dots \wedge e_{j_i}$ by
$$ 
\left({x_{m,{j_1}}\over x_{n,{j_1}}}\right)\cdots  
\left({{x_{m,{j_i}}\over x_{n,{j_i}}}}\right).
$$

Denote by $I(\seq)$ the ideal of $\Lambda$ generated by the components of 
$\seq$. Suppose that 
for each $k \in \NN$ there are $i,j\in\NN$ such that $I(\seq_i)
\subset \m^k$ and $\m^j \subset I(\seq_k)$. 
Then, for a finitely generated $\Lambda$-module $X$,
$$
\dirlim_n H^i(\seq_n, X)
$$
is the local cohomology group $H^i_\m(X)$ defined by
Grothendieck. This is shown in Sect.~3.5 of \cite{B-H} for the case
where $\seq = x_1, \dots , x_N$ is a sequence in $\Lambda$ that
generates an $\m$-primary ideal and $x_{n,i} = x_i^n$. The
construction there generalizes easily to our situation. 

\Theorem{(Grothendieck Duality for $\Lambda$)}\label{OG}
Suppose that for each positive integer $n$ we have a sequence 
$\seq_n = (x_{n,1}, \dots , x_{n, N})$ generating an $\m$-primary ideal,
such that 
 $x_{n,i} \divides x_{m,{i}}$, $n \le m$, $1 \le i \le N$, and
such that 
for each $k \in \NN$ there are $i,j\in\NN$ such that $I(\seq_i)
\subset \m^k$ and $\m^j \subset I(\seq_k)$. Then there is an
isomorphism of functors of finitely generated $\Lambda$-modules
$$
\Ext^{N-i}_\Lambda(\cdot, \Lambda)  \iso \Hom_{\ZZ_p}(\dirlim_n H^i(\seq_n, \cdot),
\QQ/\ZZ). 
$$
\endtheorem

\proof
This is the version of Grothendieck's local duality
theorem suggested by Exercise 3.5.14 of \cite{B-H}. 
The local duality theorem is 
\cite{B-H}, Theorem 3.5.8. We note that $\Lambda$ is a regular local 
ring, so it is {\it a fortiori} a Cohen-Macaulay local ring, and therefore
satisfies the hypotheses of the results in \cite{B-H} cited here.
\endproof

Now let $I \subset \Lambda$  be the augmentation 
ideal, let $g_1, \dots, g_r$ be a set of topological generators for
$\Gamma$, and let
$(T_1,
\dots , T_r)$ be the corresponding set of generators for $I$ ($g_i =
1+T_i$). Let
 $$
\omega_n(T_i)  = (1+T_i)^{p^n} - 1, \quad \nu_{n}(T_{i}) = 
{\omega_{n}(T_{i})\over T_{i}} \qquad n\ge0.
$$
There are two particularly useful choices of sequence in applying
Theorem \Ref{OG}:
$$
\seq_n = (p^n, \omega_n(T_1), \dots , \omega_n(T_r)) \quad 
\text{and}\quad \seq_n = (p^n, \nu_n(T_1), \dots , \nu_n(T_r)).
$$
Let 
$$
\seq'_n = (\omega_n(T_1), \dots , \omega_n(T_r)) \quad 
\text{(respectively $(\nu_n(T_1), \dots , \nu_n(T_r))$).}
$$
One easily deduces from the definitions an exact sequence
$$
0 \to H^{i-1}(\seq'_n, X)/p^n \to 
H^i(\seq_n, X) \to  H^i(\seq'_n, X)[p^n] \to
0.\tag\label{UsefulExactSequence}
$$
In the case $i=r$, 
the map on the right is induced by $(y_0, \dots , y_r) \mapsto y_0$.
Note that, in the case $\seq_{n} = (p^n, \omega_n(T_1), \dots , 
\omega_n(T_r))$,  $H^r(\seq'_n, X) = X/(\omega_n(T_1), \dots ,
\omega_n(T_r))$. 
Define a functor
$$
E(X) = \Hom_{\ZZ_p}(\dirlim H^r(\seq'_n, X), \QQ_p/\ZZ_p).
$$
In the case where $r=1$ and $X_n$ is torsion, $E(X)$ is Iwasawa's
original definition of the adjoint module of $X$.
The following proposition is proved in \cite{Mc} (Proposition 2).

\proposition\label{extinjection}
There is a functorial injective map $E(X) \to 
\Ext_{\Lambda}^1(X,\Lambda)$. This map is an isomorphism 
if $H^r(\seq'_n, X)$ is finite for all sufficiently large $n$.
\endproposition

We conclude with two lemmas on pseudo-null modules which will be 
needed later.

\lemma\label{pseudonulllemma}
Let $X$ be a finitely generated $\Lambda$-module. Then $X$ is pseudo-null if
and only if $X$ is torsion and $\Ext_\Lambda^1(X, \Lambda) = 0$.
\endlemma

\proof
This is also proved in \cite{H}, Theorem 11.
Since $\Lambda$ has dimension $r+1$, it follows from Theorem~\ref{OG}
that $X$ is torsion and $\Ext_\Lambda^1(X, \Lambda) = 0$ if and 
only if $H_{\m}^{r+1}(X) = H_{\m}^{r}(X) = 0$, where $\m \subset 
\Lambda$ is the maximal ideal and $H_{\m}^{*}(\cdot)$ is the local 
cohomology. By \cite{B-H} Theorem 3.5.7, these two local cohomology 
groups vanish if and only if the dimension of $X$ as a 
$\Lambda$-module is less than or equal to $r-1$. Since $R$ is a regular local 
ring, the dimension of $X$ is $r+1$ minus the height of its annihilator 
(\cite{B-H}, Corollary 2.1.4), hence the dimension of $X$ is less 
than or equal to $r-1$ if and only if the height of the annihilator of $X$ is at 
least 2, that is, if and only if $X$ is pseudo-null.
\endproof

\lemma\label{torsionext1}
Let $X$ be a finitely-generated torsion $\Lambda$-module. Then $X$ is 
pseudo-isomorphic to $\Ext^{1}(X)$. Further, 
$\Ext^{1}(X)$ is pseudo-null if and only if it is zero.
\endlemma

\proof
There is a pseudo-isomorphism $X \to E$, where 
$E$ is an elementary module, 
that is, a product of modules of the form $\Lambda/(f)$, $f \in 
\Lambda$.
Taking $\Ext^{*}(\cdot)$ 
of the sequence
$$
0 \to \Lambda \East{f}{} \Lambda \to \Lambda/(f) \to 0,
$$
we see that
$\Ext^{1}(\Lambda/(f)) = \Lambda/(f)$ if $f \neq 0$. Hence, for any 
elementary torsion-module $E$, $\Ext^{1}(E) \iso E$. Now, a 
pseudo-isomorphism 
$$
X \to E
$$
leads to a map 
$$
\Ext^{1}(E) \to \Ext^{1}(X),
$$
which is a pseudo-isomorphism, since $\Ext^{i}$ of a pseudo-null 
module is pseudo-null for any $i$ (the annihilator of $X$ also annihilates 
its $\Ext$ 
groups). Thus $\Ext^{1}(X,\Lambda)$ is pseudo-isomorphic to 
$X$. The second statement of the lemma now follows immediately 
from Lemma~\Ref{pseudonulllemma}.
\endproof

\section{The Direct Limit of the Ideal Class 
Group\label{DirectLimitSection}}

Let $X'$ be as in \Ref{StandardIwasawaModules}. For a positive 
integer $n$, let $K_{n}$ be the fixed field of $p^n \Gamma$. 
Let $X'_n = X'_{K_{n}}$. Define ideals $\omega_{n}$ and $\nu_{n,m}$ 
in $\Lambda$ by
$$
\align
\omega_{n}& = (\omega_{n}(T_1), \dots , \omega_{n}(T_{r})) \\
\nu_{n,m}& = (\nu_{n,m}(T_1), \dots , \nu_{n,m}(T_{r})), \quad \nu_{n,m}(T) = 
{\omega_n(T) \over \omega_m(T)}.\\
\endalign
$$
From \Ref{UsefulExactSequence} we get a surjective map
$$
H^r(\seq_n, X')  \onto X'/\omega_n X'[p^n]
$$
and, composing this with the natural map
$$
X'/\omega_n X' \to X'_{n},\tag\label{justanothermap}
$$
we get a map
$$
H^r(\seq_n, X') \to X'_n[p^{n}].
$$

\theorem\label{DirectLimitTheorem}
Suppose  
that $K$ has only one prime above $p$. 
Then the natural map 
$$
H^r(\seq_n, X') \to X'_{n}[p^{n}]
$$
induces an isomorphism 
$$
\dirlim_{n} H^{r}(\seq_{n}, X') \iso \dirlim X'_{n}[p^{n}] = \dirlim X'_{n}.
$$
In particular, $\Ext^{1}(X')$ is the Pontryagin dual of $\dirlim 
X'_{n}$.
\endtheorem

\proof
Since there is only one prime of $K$ above $p$, its decomposition 
group in $\Gamma$ has finite index.
Choose $n_0$ large enough so 
that all the primes above $p$ are totally ramified in $K_{\infty}/K_{n_0}$.  
Then \Ref{justanothermap} is surjective, hence the 
direct limit of $H^r(\seq_n, X') \to X'_{n}[p^{n}]$
 surjective also. Thus it 
suffices to show that the direct limit is injective.

Choose a fixed prime $\p$ of $L'$ above $p$.
Since $L'/K_{\infty}$ is $p$-split, the 
decomposition group of $\p$ in $\Gal(L'/K_{n})$ is isomorphic to 
$\Gamma_{K_{n}} = p^{n}\Gamma$; in particular, it is abelian. Denote 
this decomposition group by $\tilde \Gamma_{\p,n}$. 

Suppose that $n \ge n_{0}$. 
We claim that the kernel  of $X'/\omega_{n}X' \to X'_{n}$ is contained 
in $\nu_{n,n_0}X'/\omega_{n}X'$. First, since $\Gal(L'/K_{n}) = \tilde \Gamma_{\p,n} 
X'$ and $\tilde \Gamma_{\p,n}$ is abelian, the commutator subgroup of 
$\Gal(L'/K_{n})$ is $\omega_{n}X'$. Thus $X'/\omega_{n}X'$ is the 
maximal abelian quotient of $\Gal(L'/K_{n})$, and hence the kernel
of $X'/\omega_{n}X' \to X'_{n}$
is generated by the decomposition groups above $p$. 
Since there is only one prime of $K$ 
above $p$, all the decomposition groups are conjugate under the action 
of $\Gamma$. Thus, the subgroup of $X'$ generated by 
the decomposition groups  is generated by commutators
$[\gamma,g']$, where $\gamma\in\tilde\Gamma_{\p,n}$
and $g'\in \Gal(L'/K)$. Furthermore, we can restrict $\gamma$ to a set of 
generators of $\tilde\Gamma_{\p,n}$.
Now, as noted above, the image of $\tilde 
\Gamma_{\p,n}$ in $\Gamma$ is $p^{{n}}\Gamma$. 
So we may assume that the image of $\gamma$ in 
$\Gamma$ is $g_i^{p^n}$ for some $i$.
Choose $g$   so that
$$
g^{p^{n_0}} = \gamma_{0} x \qquad \gamma_{0}^{p^{n-n_0}} = \gamma
$$
for some $\gamma_{0} \in \tilde \Gamma_{\p,n_0}$ and some $x \in X'$. 
Let $T$ (resp. $T'$) be the image of $1-g$ (resp. $1-g'$) in $\Lambda$.
Then, using the commutator identity
$$
[a^N,b] = [a,b]^{a^{N-1}}[a,b]^{a^{N-2}}\dots[a,b],
$$
we find
$$
\nu_{n_0}(T)[g,g'] = [g^{p^{n_0}},g']= [\gamma_{0} x, g'] = [\gamma_{0},g'] + 
[x, g'] = [\gamma_{0},g'] + T'x.
$$
Thus
$$
\multline
\nu_{n}(T)[g,g'] = \nu_{n,n_0}(T)([\gamma_{0},g'] + T'x) = \\
[\gamma_{0}^{p^{n-n_0}},g'] + \nu_{n,n_0}(T)T'x = 
[\gamma, g'] + 
\nu_{n,n_0}(T)T'x.
\endmultline
$$
Hence, since $[g,g'] \in X'$ and $\nu_{n,n_0}$ divides $\nu_{n}$, 
$[\gamma,g'] \in \nu_{n,n_0}(T) X'$. This proves the claim.

We apply the results 
of Section~\ref{Extsection} using 
$$
\seq_n = (p^n, \omega_n(T_1), \dots , \omega_n(T_r)).
$$
and
$$
{\bold y}_n = (p^n, \nu_{n,n_0}(T_1), \dots , \nu_{n,n_0}(T_r)).
$$

The fact that  the kernel  of $X'/\omega_{n}X' \to X'_{n}$ is contained 
in $\nu_{n,n_0}X'/\omega_{n}X'$ implies that
$H^{r}({\bold y}'_n, X)$ is a quotient of $X'_{n}$. Hence it is 
finite, so by Proposition~\ref{extinjection}, $$\dirlim H^{r}({\bold y}'_n, X')
\iso D = \Hom(\Ext^{1}(X'),\QQ_{p}/\ZZ_{p}).$$ On the 
other hand, it  follows from Theorem~\ref{OG} that 
$$
D \iso \dirlim_{n} H^r(\seq_n, X').
$$
Hence
$$
\dirlim_n H^{r}(\seq_{n},X') \iso \dirlim_n H^{r}({\bold y}'_n, X').
$$
Since $H^{r}(\seq_{n},X') \to H^{r}({\bold y}'_n, X')$ 
factors through the  $X'_{n}\to H^{r}({\bold y}'_n, X')$, this 
implies 
that $\dirlim_n H^{r}(\seq_{n},X') \to \dirlim_n X'_{n}$ is
injective. This concludes the proof of the theorem.
\endproof

\section{Interpretation of Greenberg's Conjecture in Terms of 
$Y$\label{YInterpretation}}

For each prime $\p$ of $K$ lying above
$p$, let $r_\p$ be the integer such that the decomposition group $\Gamma_\p$
 of
$\p$ in
$\Gamma$ is isomorphic to  $\ZZ_p^{r_\p}$. Let $Y$ and $X'$ be as in 
\Ref{StandardIwasawaModules}. For any $\Lambda$-module $M$ and $i \in 
\ZZ$ we denote by $M(i)$ the $i$-th Tate twist of $M$.

\theorem\label{reinterpretation}
Suppose that $K_\infty$ contains the $p^n$-th roots of unity for every
positive integer $n$, and that
$r_\p
\ge 2$ for all primes $\p$ of $K$ lying above $p$. Then 
$\Ext^{1}_{\Lambda}(X',\Lambda)(1)$ is pseudo-isomorphic to $Y_{\tor}$, 
and is isomorphic to $Y_{\tor}$ if $r_\p \ge 3$ for all $\p$.
\endtheorem

In order to prove Theorem~\ref{reinterpretation}, we  relate
$\Ext^{1}(Y, \Lambda)$ to $X'$ using
a duality theorem of Jannsen \cite{J2}.  Then we use 
a basic structure theorem for
$Y$ to turn this around into a relation between $Y_{\tor}$ and 
$\Ext^{1}(X',\Lambda)$. The
structure theorem, due independently to Jannsen \cite{J1} and Nguyen-Quang-Do
\cite{Ng}, may be
 formulated conveniently in terms of an auxiliary Iwasawa-module
$Z$, defined by Nguyen-Quang-Do
\cite{Ng}. 

\Theorem{(\cite{J1}, \cite{Ng})}
Suppose that $K_\infty$ contains the $p^n$-th roots of unity for every
positive integer $n$. There is a $\Lambda$-module $Z$ with a presentation
$$
0  \to \Lambda^s \to \Lambda^g \to Z \to
0,\tag\label{StructureSequence}
$$
for some non-negative integers $s$ and $g$, such that $Y$ fits in a short 
exact sequence
$$
0 \to Y \to Z \to I \to 0,
$$
where $I \subset \Lambda$ is the augmentation ideal.
\endtheorem

We define $Z$ as follows. Let
$\Omega$ be the maximal $p$-ramified pro-$p$-extension of 
$K$,  ${\Cal G}= \Gal(\Omega/K)$,
and  ${\Cal I}$  the augmentation ideal in $\ZZ_p[[{\Cal G}]]$.
Define
$$
Z_F = H_0(\Omega/F,{\Cal I}), \qquad Z = Z_{K_{\infty}}.
$$
If $K \subset F\subset
K_\infty$, the
identification 
$$
H_0(\Gal(K_\infty/F),Z) = Z_F
$$ 
gives a natural surjective map
$$
Z \onto Z_F. 
$$
Let  $Y_F$ be the Galois group of the maximal abelian $p$-ramified
pro-$p$ extension of
$F$. Denote by $I_F$ the augmentation ideal in
$\ZZ[\Gal(F/K)]$. 
By taking $\Gal(\Omega/F)$ homology of the exact
sequence
$$
0 \to {\Cal I} \to \ZZ_p[[{\Cal G}]] \to \ZZ_p \to 0,
$$
one obtains an exact sequence
$$
0 \to Y_F \to Z_F \to I_F \to 0,\tag\label{ZLsequence}
$$
functorial in $F$.

In
particular, with $F = K_{\infty}$, we obtain the second sequence in the 
statement of the theorem. The integers $g$ and $s$ in the theorem are
$$
\align 
g &= \dim_{\ZZ/p\ZZ} H^1({\Cal G}, \ZZ/p\ZZ)\\
s &=\dim_{\ZZ/p\ZZ}H^2({\Cal G}, \ZZ/p\ZZ)
\endalign
$$
As explained in \cite{Ng}, \Ref{StructureSequence} may be derived 
from the
presentation of ${\Cal G}$ as a pro-$p$ group with $g$ generators and $s$
relations. (Injectivity on the left in \Ref{StructureSequence} depends on the assumption that $K_{\infty}$ 
 contains all $p$-power roots of unity, 
 which implies that $K_{\infty}$ satisfies the weak Leopoldt 
 conjecture, $H^{2}(K_{\infty},\QQ_{p}/\ZZ_{p}) = 0$; see \cite{Ng} for
details.) 

The presentation~\Ref{StructureSequence} 
allows us to compute $\Ext^i(Y) =
\Ext^i_\Lambda(Y,\Lambda)$ for $i \ge 2$. Although it is not 
needed in what follows, the answer is simple so we present it here.
\proposition\label{ExtY}
Let $r = \rank_{\ZZ_p} \Gamma$.
If $i \ge 2$, then 
$$
\Ext^i(Y) = \cases 0 & i \neq r-2 \\
\ZZ_p & i = r-2. \endcases
$$
\endproposition

\proof
Since $Z$ has a free resolution of length 2,
$\Ext^i(Z) = 0$ if $i \ge 2$.  
Thus, taking $\Ext$ of \Ref{ZLsequence} with $F = K_\infty$, we see that 
$\Ext^i(Y) \iso \Ext^{i+1}(I)$ if
$i\ge 2$. Now use the following lemma. 
\endproof
\lemma\label{ExtI}
If $i\ge 1$, then
$$
\Ext^i(I) = \cases 0 & i \neq r-1 \\ \ZZ_p & i = r-1.
\endcases
$$
\endlemma

\proof
Taking $\Ext$ of 
$$
0 \to I \to \Lambda \to \ZZ_p \to 0,
$$
we see that $\Ext^{i}(I) \iso
\Ext^{i+1}(\ZZ_p)$ if $i \ge 1$. But, using the Koszul resolution
of $\ZZ_p$, one can see that this latter group is zero unless $i+1 =
r$, in which case it is $\ZZ_p$.
\endproof

\medskip
The structure of $\Ext^i(Y)$ when $i=1$ is more subtle. A 
duality theorem  between $X'$ and $Y$, due to Jannsen,
relates it to 
the Iwasawa module
$$
H^2_\infty = \invlim_F H^2({\Cal O}'_F, \ZZ_p(1)), 
$$
which is related to $X'$ by the exact sequence
$$
0 \to X' \to H^2_\infty \to \sum_{\p\divides p}
\ZZ_p[[\Gamma/\Gamma_{\p}]] \to \ZZ_p \to 0,\tag\label{XprimetoH2}
$$
where the sum is over all primes of $K$ dividing $p$. We denote by 
$X(n)$ the usual Tate twist of a $\Lambda$-module $X$ by $\ZZ_p(1)^{\tensor n}$, 
where $\ZZ_p(1)= \invlim \mu_{p^n}$.

 \Theorem{(\cite{J2}, Theorem 5.4d)}\label{DualityTheorem}
If $r \neq 2, 3$, there is an isomorphism
$$
 H^2_\infty \iso  \Ext^1(Y, \Lambda)(1).
$$
If $r=2$ there is an exact sequence
$$
\ZZ_p(1) \to H^2_\infty \to  \Ext^1(Y, \Lambda)(1) \to 0,
$$
and if $r = 3$, there is an exact sequence
$$
0 \to H^2_\infty \to  \Ext^1(Y, \Lambda)(1) \to
\ZZ_p(1) .
$$
\endtheorem

\demo{\it Proof of Theorem \Ref{reinterpretation}}
Let $\p$ be a prime of $K$ lying above $p$. The annihilator of  $\Lambda_\p = \ZZ_p[[\Gamma/\Gamma_\p]]$ as a 
$\Lambda$-module is the augmentation ideal in $\ZZ_p[[\Gamma_\p]]$, which has height $r_\p$.
Thus, the hypothesis on $r_\p$ implies that $\Lambda_{\p}$ is pseudo-null, 
and hence has vanishing $\Ext^1$ by Lemma~\ref{pseudonulllemma}. 
Furthermore, if $r_\p \ge 3$ then its $\Ext^2$ vanishes as well---the argument
for this is entirely similar to that in the proof of Lemma~\ref{pseudonulllemma}.
There is a natural surjection
$X
\to X'$, whose kernel is generated as a
$\ZZ_p$-module by the Frobenius automorphisms corresponding to primes above
$p$, and is therefore a finitely generated module over $\bigoplus_\p \Lambda_\p$.
Hence the kernel is pseudo-null, and so 
$\Ext^{1}(X)$ is isomorphic to $\Ext^{1}(X')$.

From \Ref{XprimetoH2} we get an injection $X' \into H^{2}_{\infty}$ whose 
cokernel is pseudo-null, and has vanishing $\Ext^{2}$ if $r_\p \ge 3$ for all
$\p$. 
Thus $\Ext^{1}(X')$ is pseudo-isomorphic to 
$\Ext^{1}(H^{2}_{\infty})$, and is isomorphic if $r_\p \ge 3$ for all $\p$.

Since 
$\ZZ_p(1)$ has vanishing
$\Ext^1$ when $r\ge2$ and vanishing $\Ext^2$ when $r \ge 3$, 
Theorem~\Ref{DualityTheorem} implies
that
$\Ext^1(H^{2}_{\infty})(1)$ is isomorphic to $\Ext^1(\Ext^1(Y))$.
It follows from \Ref{ZLsequence} and
Lemma~\Ref{ExtI} that there is a pseudo-isomorphism 
$\Ext^1(Z)\to \Ext^1(Y)$. Furthermore, this is an injection with 
cokernel  $\ZZ_{p}$ if $r=3$, and an isomorphism if $r > 3$.
Hence from 
Lemma~\Ref{torsionext1}  that 
$\Ext^1(\Ext^1(Y))$ is pseudo-isomorphic to 
$\Ext^1(\Ext^1(Z))$, isomorphic if $r \ge 3$.
An elementary calculation with the two-step free resolution of $Z$
shows that
$$
\Ext^1(\Ext^1(Z)) \iso Z_\tor.
$$ 
Finally, it follows from
\Ref{StructureSequence} that $Z_\tor \iso Y_\tor$.
\endproof

\corollary\label{greenbergcorollary}
Suppose that $K_\infty$ contains the $p^n$-th roots of unity for every
positive integer $n$, and that
$r_\p
\ge 2$ for all primes $\p$ of $K$ lying above $p$. 
Then Greenberg's conjecture is satisfied if and only if 
$Y_{\tor}=0$.
\endcorollary

\proof
Theorem~\ref{reinterpretation}
 and Lemma~\Ref{torsionext1} imply the corollary directly 
 in the case $r_{\p} \ge 3$. In the case $r_{\p}=2$ they
  imply that 
 Greenberg's conjecture is equivalent to the pseudo-nullity of 
 $Y_{\tor}$. However, since $K_{\infty}$ contains $\mu_{p^{\infty}}$, 
 it verifies the weak Leopoldt conjecture, hence $Y$ has no non-trivial
 pseudo-null submodule (see \cite{Ng}), so $Y_{\tor}$ is pseudo-null 
 if and only if  it is zero.
\endproof

This result has also been shown in \cite{LN}.
We conclude this section by proving Theorem~\ref{WingbergContradictionTheorem},
which we restate here for convenience.

\theorem"2"
Suppose that $\QQ(\zeta_p)$ satisfies Greenberg's conjecture.
Then ${\Cal G}$ has a pro-$p$-free quotient
of rank $g-s$ if and only if $p$ is regular.
\endtheorem

\proof
Suppose that there is a surjective map ${\Cal G} \to {\Cal F}$, where
${\Cal F}$ is a free group of rank $g-s$. Then, from the definition of $Z$, 
we see that there is a surjective map 
$Z \to \Lambda^{g-s}$. On the other hand, $Z$ has rank $g-s$, so we must have
$$
Z \iso \Lambda^{g-s} \oplus Z_{\text{tor}}.
$$
By Corollary~\ref{greenbergcorollary}, $Y_{\text{tor}} = 0$, so $Z$ is free. Since $Z/IZ \iso Z_K \iso Y_K \iso 
{\Cal G}^{\text{ab}}$, this
implies that the maximal abelian $p$-ramified $p$-extension of $K$ has Galois 
group isomorphic to $\ZZ_p^{g-s}$; in particular, there is no torsion in the 
Galois group. It follows from the theory of cyclotomic fields that $p$ is 
regular (see, for example, the proof of Lemma~\ref{bla}).
\endproof

\section{A Generalization of Iwasawa's Theorem on Units\label{UnitSection}}
In this section we  prove Theorem~\Ref{torsionfreetheorem}, whose 
statement we recall here.
Recall
$$
N_\infty = K_\infty(\epsilon^{1/p^{n}}: n \in {\Bbb N}, \epsilon 
\in {\Cal O}_{K_{\infty}}[1/p]^\times).
$$
\theorem"3"
Suppose that $K_\infty/K$ is a multiple $\ZZ_p$-extension satisfying:
\list
\item $K_\infty$
contains all the
$p^n$-th roots of unity,
$n>0$ 
\item there is only one prime of $K$ above
$p$.
\endlist
Then the $\Lambda$-module $Y' = \Gal(N_\infty/K_\infty)$ 
is torsion-free.  In particular, $Y_{\tor}$ fixes $N_{\infty}$.
\endtheorem
For a finite extension $F/K$ contained in 
$K_\infty$, let
$E_F = {\Cal O}_F[1/p]^\times
$, and let $\Un \subset E_F$ be the group of
universal norms, that is, the elements of $E_F$ that are norms of elements
of $E_{L}$ for every extension $L/F$ contained in $K_\infty/F$.
Let
${\Cal E} = \dirlim_F E_F$, and consider  the subgroup ${\Cal E}^u  = \dirlim_F
\Un
\subset {\Cal E}$. Define $N^u_\infty$ to be the extension of $K_\infty$
obtained by adjoining all $p$-power roots of all elements of ${\Cal E}^u$, and
let
$Y^u = \Gal(N^u_\infty/K_\infty)$.

\lemma\label{UnitProposition}
Suppose that $K_\infty$ contains all $p$-power roots of unity. Then 
$Y^u$ is $\Lambda$-torsion free.
\endlemma

\proof
From Kummer theory, we have an isomorphism of $\Lambda$-modules
$$
Y^u \iso \Hom_{\ZZ_p}({\Cal E}^u \tensor \QQ_p/\ZZ_p,
\QQ_p/\ZZ_p(1)).
$$
Let $E = \invlim_F \Un = \invlim_F E_F$, regarded as a
$\Lambda$-module. We define a homomorphism
$$
\Hom_{\ZZ_p}({\Cal E}^u \tensor \QQ_p/\ZZ_p, \QQ_p/\ZZ_p(1)) \to
\Hom_\Lambda(E, \Lambda(1))\tag\label{UnitMap}
$$
by mapping $\chi$ to the homomorphism
$$
(u_F) \mapsto (\sum_{\sigma \in \Gal(F/K)} q(F)_{*}\chi(u_F^\sigma 
\tensor
q(F)^{-1}) \sigma),
$$
where $q(F)$ is the number $p$-power roots of unity in $F$, and 
$$
q(F)_{*}:  q(F)^{-1}\ZZ/\ZZ\to \ZZ/q(F)\ZZ
$$ is induced by multiplication by $q(F)$.
Since every universal norm $u \in \Un$ occurs as part of a norm
compatible sequence,  $\chi$
 is in the kernel of this map if and only if $\chi = 0$, hence
\Ref{UnitMap} is injective. The lemma now follows from the fact
that $\Hom_\Lambda(E, \Lambda(1))$ is torsion-free, since $\Lambda(1)$ is 
free of rank 1.
\endproof
It follows from the lemma that the torsion submodule of $Y'$ is contained in 
the kernel of
$Y'\to Y^u$, which, by Kummer theory, is the Pontriagin dual of 
$\dirlim_F (E_F/\Un) \tensor \QQ_p/\ZZ_p$. 
We  prove Theorem~\Ref{torsionfreetheorem} by showing that the direct
limit is zero. We define a filtration $\Un \subset \No \subset
\Lo
\subset E_F$ by
$$
\align
\No &=  \{ x \in E_F: x \in N_{L/F}L^\times \text{ for all finite 
$L/F$ in $K_\infty$}\}\\
\Lo & =  \{ x \in E_F:x_v \in N_{L_v/F_v} L_v^\times,
\text{  all finite $L/F$ in $K_\infty$, all  valuations
$v$ of
$L$}\}\\
\endalign
$$

We  show $\dirlim_{F} (E_{F}/\Un)\tensor(\QQ_{p}/\ZZ_{p}) = 0$ 
by considering each graded piece in turn. 
We  make use of the following  lemmas, whose proofs are 
elementary and left to the reader:
\lemma\label{keytool}
Suppose that $(M_j)$, $j \in J$, is a direct system of $\ZZ_p$-modules, and
that for every positive integer $n$ and every $j$ there exists $j'\ge j$ such
that the image of 
$$M_j
\to M_{j'}$$ is contained in $p^n
M_{j'}$. Then 
$$
\dirlim_j M_j \tensor \QQ_p/\ZZ_p = 0.
$$
\endlemma

For abelian groups  $H \subset \Gamma$ of finite index, define
$$
e(\Gamma,H) = {\text{index of $H$ in $\Gamma$} \over (\text{exponent of
$\Gamma/H$})}.
$$
\lemma\label{e(H)}
If $\Gamma \iso \ZZ_p^r$
with $r  >1$, then $$\lim_H e(\Gamma,H) = p^\infty,$$ where the limit is
over open neighborhoods $H$ of the identity. 
\endlemma

\topic{Units modulo local norms}
 Let
$\Gamma_F =
\Gal(K_\infty/F)$, and, for a prime
$\p$ of
$F$, denote by $\Gamma_{F,p}$
 the decomposition group of
$\p$ in $\Gamma_{F}$ (where the context removes
ambiguity, we  denote $\Gamma_{F,\p}$ simply by $\Gamma_{\p}$).

\proposition\label{lemma1}
Suppose that $r \ge 2$ and that  there is only one prime of $K$ lying above $p$.
Then 
$$\dirlim((E_F/\Lo)\tensor
\QQ_{p}/\ZZ_p) = 0.$$
\endlemma

\proof
Recall that  $Y_F$ and $X'_{F}$ are the Galois groups respectively 
of the maximal abelian 
$p$-ramified
pro-$p$-extension of
$F$, and the maximal abelian unramified pro-$p$-extension in which all primes 
above $p$ split completely. By class field theory we have an exact sequence
$$
E_{F} \tensor \ZZ_{p} \to \prod_{\p \divides p} \Gal(\bar F_\p^{\text{ab}}/F_\p)
\to Y_{F} \to X'_{F} \to 0,
$$
where the map on the left is the product of the local reciprocity maps
over all primes of $F$ above $p$. Let $H_{\p}$ be the kernel of 
$\Gal(\bar F_\p^{\text{ab}}/F_\p)\to 
\Gamma_{F,\p}$.
The reciprocity map takes an element to $H_{\p}$ if and only if it is 
a local universal norm at $\p$. Let 
$H_{F} \subset Y_{F}$ be the image of $\prod_{\p}H_{\p}$ and let 
$W_{F} = Y_{F}/H_{F}$.
Then we get an exact sequence
$$
0 \to E_F/\Lo\tensor \ZZ_{p} \to \prod_{\p}\Gamma_{F,\p} \to 
W_{F}  \to X'_{F} \to 0.
$$
Thus it suffices to show that 
$$
\dirlim_{F} \prod_{\p}\Gamma_{F,\p}\tensor \QQ_{p}/\ZZ_{p} = 0
$$
and 
$$
\dirlim_{F} W_{F,\tor} \into \dirlim_{F} X'_{F}.
$$
The first follows easily from Lemmas~\Ref{e(H)} and \Ref{keytool}, 
since, for an extension $F'/F$ contained in $K_{\infty}$, 
the transition map $\Gamma_{F,\p} \to \Gamma_{F',p}$ is the transfer, 
which is multiplication by the index $[\Gamma_{F,\p}:\Gamma_{F',p}]$. 
Hence its image lies in $e(\Gamma_{F,\p},\Gamma_{F',p})\Gamma_{F',p}$.

To prove the second, note that $\invlim \Gamma_{F,\p} = 0$, 
hence $$\invlim_{F}W_{F}=X' = 
\Gal(L'/K_{\infty}),$$ where $L'$ and $X'$ are as in 
\Ref{StandardIwasawaModules}. Let $K_{n}$ be as in 
Section~\ref{DirectLimitSection}, and let $\tilde \Gamma_{\p,n}$
be a decomposition group of $\p$ in $\Gal(L'/K_{n})$.
Choose $n$ large enough so 
that all the primes above $p$ are totally ramified in $K_{\infty}/K_{n}$. 
Then, as pointed out in the proof of 
Theorem~\ref{DirectLimitTheorem}, the commutator subgroup of $\Gal(L'/K_{n})$ is 
$\omega_{n}X'$. Hence we have an exact sequence
$$
0 \to X'/\omega_{n}X' \to W_{K_{n}} \to \Gamma_{K_{n}}.
$$
Since $\Gamma_{K_{n}}$ is torsion-free, 
$W_{K_{n},\tor} \subset X'/\omega_{n}X'$.

Thus it suffices to show that 
$$
\dirlim_{n}(X'/\omega_{n}X')_\tor \into \dirlim_{n} X'_{n},
$$
or, equivalently, that
$$
\dirlim_{n}(X'/\omega_{n}X')[p^n] \into \dirlim_{n} X'_{n}[p^n].
$$
Since the map 
$$
(X'/\omega_{n}X')[p^n] \to X'_{n}[p^n]
$$
is a factor of 
$$
H^r(\seq_n, X') \onto (X'/\omega_{n}X')[p^n] \to  X'_{n}[p^n],
$$
this follows from Theorem~\ref{DirectLimitTheorem}.

\endproof

\goodbreak

\topic{Local norms modulo global norms}

 For a $\ZZ_p$-module $M$, let $M^\wedge$ denote $\Hom_{\ZZ_p}(M,\QQ_p/\ZZ_p)$.
\lemma\label{PartTwo}
Suppose that no prime splits completely in $K_{\infty}/K$. Let $F/K$ 
be 
a finite subextension of $K_{\infty}/K$.
There is a surjective map of $\Gal(F/K)$-modules
$$
H^2(\Gamma_F, \QQ_p/\ZZ_p)^\wedge \onto { \Lo \over \No}.
$$
\endlemma

\proof
Let $\II_F$ be the id\`ele group of
$F$ and let
${\bold C}_F = \II_{F}/F^{\times}$ be the id\`ele class group. It 
follows from the hypothesis, and the fact that $K_\infty/K$ is
unramified outside $p$, that
any element of $F^\times$ that is a universal local norm must  be a
$p$-unit, and so
$\Lo/\No$ 
is
the kernel of 
$$
\invlim_{L} \hat H^0(L/F,
L^{\times}) \to \invlim_{L} \hat H^0(L/F,
\II_{L}),
$$
where the limit is taken over all finite extension $L/F$ in $K_\infty$.
Taking $\Gal(L/F)$-cohomology of the exact sequence
$$
0 \to L^\times \to \II_{L} \to {\bold C}_{L}\to 0
$$
yields an exact sequence
$$
\hat H^{-1}(L/F, {\bold C}_{L})
\to
\hat H^0(L/F, L^\times) \to \hat H^0(L/F,
\II_{L}).
$$
By class field theory,
$$
\hat H^{-1}(L/F, {\bold C}_{L}) = \hat H^{-3}(L/F, \ZZ) =
H^2(L/F,
\QQ/\ZZ)^\wedge = H^2(L/F, \QQ_p/\ZZ_p)^\wedge.
$$
Thus $\Lo/\No$ is
a quotient of 
$$
\invlim_L H^2(L/F, \QQ_p/\ZZ_p)^\wedge = (\dirlim H^2(L/F,\QQ_p/\ZZ_p))^\wedge 
= H^2(\Gamma_F, \QQ_p/\ZZ_p)^\wedge.
$$
\endproof

\lemma\label{technicallemma}
Let $\Gamma$ be an abelian group isomorphic to $\ZZ_p^r$   and let $\Gamma'$
be a subgroup of $\Gamma$ of finite index. Then the
image of $H^2(\Gamma, \QQ_p/\ZZ_p)^\wedge \to H^2(\Gamma',\QQ_p/\ZZ_p)^\wedge$
(dual to the corestriction  map)
 is
contained in
$e(\Gamma,\Gamma')H^2(\Gamma',\QQ_p,\ZZ_p)^\wedge$.
\endlemma

\proof
Dually, we prove that the corestriction map 
$$
H^2(\Gamma',\QQ_p/\ZZ_p) \to H^{2}(\Gamma,\QQ_{p}/\ZZ_{p})
$$
vanishes on 
$H^{2}(\Gamma',\QQ_{p}/\ZZ_{p})[{e(\Gamma,\Gamma')}]$.
It is shown in \cite{Br} that $H^{2}(\Gamma',\QQ_{p}/\ZZ_{p})$ is generated by 
cup products of characters $\chi'_{1}, \chi_{2}' \in H^{1}(\Gamma',\QQ_{p}/\ZZ_{p})$. 
Furthermore, if $e(\Gamma,\Gamma')$ kills $\chi'_{1}\cup\chi'_{2}$, 
then it kills one of the characters, say $\chi'_{2}$.
The corestriction map on characters is composition with the transfer
map $\Gamma \to \Gamma'$, which is just multiplication by 
$[\Gamma:\Gamma']$ in this case. If we 
write $\chi_{i}' = \res \chi_{i}$ for some characters $\chi_{i}$ of 
$\Gamma$,
$i=1,2$, 
and let
$p^{e} = \ord(\chi_{2})/\ord(\chi'_{2})$,
then $p^{e}$ divides the exponent of $\Gamma/\Gamma'$, and
$$
\cores(\chi_{1}'\cup\chi_{2}') = (\chi_{1} \cup \cores \chi_{2}') = 
[\Gamma:\Gamma'](\chi_{1}\cup\chi_{2}).
$$
Now, since $e_{2}(\Gamma,\Gamma')$ kills $\chi_{2}'$, 
$p^{e}e_{2}(\Gamma,\Gamma')$ kills 
 $\chi_{2}$. But $p^{e}$ divides the exponent of $\Gamma/\Gamma'$, and 
 hence $p^{e}e_{2}(\Gamma,\Gamma')$ divides $[\Gamma:\Gamma']$, 
 hence $[\Gamma:\Gamma']$ kills $\chi_{2}$.
Thus the right hand side 
 of the displayed equation is zero, as required.
\endproof

\proposition\label{lemma2}
Suppose that $K_{\infty}$ contains all $p$-power roots of unity and 
that $r\ge 2$. Then the direct system $(\Lo/\No) \tensor \ZZ_p$ satisfies the
conditions of Lemma~\Ref{keytool}, and hence $\dirlim 
(\Lo/\No)\tensor \QQ_{p}/\ZZ_{p} = 0$. 
\endproposition

\proof
This follows from  Lemma~\Ref{PartTwo} and Lemma~\Ref{technicallemma}.
\endproof

\topic{Global norms modulo norms of units}

\lemma\label{PartThree}
Let $F/K$ be a finite extension contained in $K_{\infty}$. There is a surjective map of $\Gal(F/K)$-modules
$$
H_1(\Gamma_F, X') \onto {\No
\over \Un}.
$$
\endlemma

\proof
Let 
$I_F$ be the group of ideals in 
${\Cal O}_{F}[1/p]$, $P_F$  the group of  principal ideals, and $C'_F =
I_F/P_F$. For a finite extension $L/F$ in $K_\infty$, consider the exact sequences 
$$ 0 \to E_L \to L^\times \to 
P_L \to 0 
$$
and 
$$ 0 \to P_L \to I_L \to C'_L \to 0. 
$$
Since $I_L$ is a direct sum of induced $\Gal(L/F)$-modules, it is
cohomologically trivial. Hence cohomology of the second
sequence yields
$$
H_1(L/F, C'_{L}) = \hat{H}^{-2}(L/F, C'_{L}) \iso
\hat{H}^{-1}(L/F, P_{L}).
$$
Cohomology of the first sequence gives
$$
\hat{H}^{-1}(L/F, P_{L}) \to \hat{H}^0(L/F, E_L) 
\to \hat{H}^0(L/F, L^\times).
$$
Splicing these together we get
$$
H_1(L/F, C'_{L}) \to {E_F \cap N_{L/F}L^\times \over
N_{L/F}E_L}
\to 0.
$$
Taking the inverse limit over $L$, and noting that, since $L/F$ is a
$p$-extension, $H_1(L/F, C'_L) \iso H_1(L/F, A'_{L})$, we obtain the
required map. 
\endproof

{\it Remark.\/} If ${K_{n}}$ is as in Section~\ref{DirectLimitSection}, then 
$$H_{1}(\Gamma_{K_n}, X') = H^{r-1}(\seq_{n}',X'),$$ where $\seq_{n}' = 
(\omega_{n}(T_{1}), \dots , \omega_{n}(T_{r}))$.

\proposition\label{lemma3}
Suppose $K$ has exactly one prime above $p$. Then 
$\dirlim (\No/\Un)\tensor \QQ_p/\ZZ_p = 0$.
\endproposition

\proof
From Lemma~\Ref{PartThree} and the remark above it suffices to prove 
that 
$$\dirlim_{n} H^{r-1}(\seq_{n}',X') \tensor \QQ_{p}/\ZZ_{p}= 0.$$ The exact sequence 
\Ref{UsefulExactSequence} with $i=r$ yields
$$
0 \to H^{r-1}(\seq'_n, X')/p^n \to 
H^r(\seq_n, X') \to (X'/\omega_{n}X')[p^n] \to
0.
$$
Theorem~\ref{DirectLimitTheorem} implies the direct limit of the 
right hand arrow is an isomorphism; hence the direct limit of the 
group 
on the left is zero.
\endproof

\demo{\it Proof of Theorem~\Ref{torsionfreetheorem}}
The case $r=1$ follows from \cite{Iw, Theorem 15}.
Hence we may assume $r \ge 2$. Recall from the discussion at the beginnning of this section that
$Y^u$ is torsion-free, and that the kernel of $Y'
\to Y^u$ is the Pontriagin dual of $\dirlim_F (E_F/\Un) \tensor \QQ_p/\ZZ_p$.
It follows from Propositions~\Ref{lemma1},~\Ref{lemma2}, and \Ref{lemma3} 
that  this direct limit is zero, and hence the theorem follows.
\enddemo

\section{A  Criterion for the Vanishing of $Y_\tor$}

We recall the notation of Section~\ref{YInterpretation}, in particular 
$$
\align 
g &= \dim_{\ZZ/p\ZZ} H^1({\Cal G}, \ZZ/p\ZZ)\\
s &=\dim_{\ZZ/p\ZZ}H^2({\Cal G}, \ZZ/p\ZZ)
\endalign
$$
Note that
if $s = 0$, then, by \Ref{StructureSequence}, $Z$ is a free
$\Lambda$-module, hence $Y_\tor = 0$. Our aim is to develop a test
in the next simplest case. Let $K^{\text{ab}}$ be the maximal 
abelian $p$-ramified
pro-$p$-extension of $K$.
\lemma\label{criterion}
Suppose that $K$ satisfies Leopoldt's conjecture and that $s = 1$. 
If 
$Y_\tor$  fixes $K^{\text{ab}}$, then $Y_{\tor} = 0$.
\endlemma

\proof
Since $s=1$, the resolution of $Z$ looks like 
$$
0 \to \Lambda \to \Lambda^g \to Z \to 0.\tag\label{rank1}
$$
Let $\lambda = (\lambda_1, \dots , \lambda_g)$ be the image of $1$ under the
left  hand map in this resolution. Choose an element $y \in Y_{\tor}$
such that if $y = f y'$, $f \in \Lambda$, $f, y' \neq 0$, then $f$ is a 
unit in $\Lambda$. Let
$(\mu_{1}, \dots , \mu_{g}) \in \Lambda^g$ map to $y$. Then 
the $\mu_{i}$ are relatively prime, and 
there exist
$f, h \in
\Lambda$,
$f
\neq 0$,  such that 
$$
f\cdot (\mu_1, \dots , \mu_g) = h\cdot (\lambda_1, \dots, \lambda_g).
$$
Since $\Lambda$ is a unique factorization domain and the $\mu_{i}$ 
are relatively prime, $h$ divides $f$; replacing $f$ by $f/h$ we 
get 
$$
f\cdot (\mu_1, \dots , \mu_g) =  (\lambda_1, \dots,
\lambda_g).\tag\label{multiple}
$$
(This last step depended on the assumption $s=1$ in a crucial way.)
Now take $H_0(\Gamma, \cdot)$ of
the sequence
\Ref{rank1}, noting that from \Ref{ZLsequence} we have $H_0(\Gamma,Z) = Z_K
= Y_K = \Gal(K^{\text{ab}}/K)$:
$$
0 \East {}{} \ZZ_{p} \East {i}{} \ZZ_{p}^g \East {}{} \Gal(K^{\text{ab}}/K)
\East{}{} 0
$$
The fact that $i$ is injective is a consequence of
Leopoldt's conjecture. Indeed, as explained in \cite{Ng}, the 
kernel, $H_{1}(K_{\infty}/K,Z)$, is isomorphic to the Pontriagin dual of 
$H^{2}(K,\QQ_{p}/\ZZ_{p})$, and the vanishing of this latter 
cohomology group is an equivalent formulation of Leopoldt's conjecture.
We have $i(1) = \lambda(0)$; since $i$ is injective, $\lambda(0) \neq 0$. 
It follows that $f(0)
\neq 0$. Since $f(0)\mu(0) = \lambda(0)$, the power of $p$ dividing 
$f(0)$ is the order of the image of  $\mu(0)$  in
$\Gal(K^{\text{ab}}/K)$. By hypothesis, that order is 1, hence $f(0)$ 
is a unit. Hence $f$ itself is a unit power series, and thus $y=0$.
\endproof

In view of Theorem~\Ref{torsionfreetheorem}, our first step in applying the
criterion developed in Lemma~\Ref{criterion} is to restrict ourselves to a
class of multiple
$\ZZ_p$-extensions
$K_\infty/K$ for which
$K^{\text{ab}} \subset N_\infty$. The following lemma gives one such class
of extensions. Recall that 
$A$ is the
$p$-primary part of the ideal class group of $K$, 
$E = {\Cal O}_{K}^{\times}$,
$U =\prod_{\p\divides 
p}{\Cal O}_{K_\p}^{\times}$, and that  $\overline E$ is the closure 
of $E$ in $U$.

\lemma\label{bla}
Suppose that $K=\QQ(\zeta_p)$.
\list
    \item Suppose that $A(K) \iso \ZZ/p\ZZ$. Then $s=1$.
    \item Suppose further that $(U/\overline E)[p^{\infty}] \iso \ZZ/p\ZZ$. 
    Then $G^{\text{ab}}[p^{\infty}] \iso \ZZ/p\ZZ$.
    \item Suppose further that $K_\infty$ contains the $p^n$-th roots of
unity for all $n>0$. Then $K^{\text{ab}} \subset N_\infty$.
\endlist
\endlemma

\proof
(1) 
We must show that $H^{2}({\Cal G}, \ZZ/p\ZZ) \iso \ZZ/p\ZZ$.
Since $\mu_{p} \subset K$, we have $H^{2}({\Cal G}, \ZZ/p\ZZ) \iso 
H^{2}({\Cal G},\mu_{p})$. The latter may be calculated by regarding it 
as an \'etale cohomology group over $S = \Spec({\Cal O}[1/p])$, and 
 using the Kummer sequence of \'etale sheaves on $S$
$$
1 \to \mu_{p} \to \Gm \to \Gm \to 1.
$$
Let $C$ be the ideal class group of $K$.
Since the prime of $K$ above $p$ is principal, we have $C \iso 
\Pic({\Cal O}_{K}[1/p]) \iso H^{1}_{\text{\'et}}(S, \Gm)$. 
Furthermore, $H^{2}_{\text{\'et}}(S, \Gm)=0$, since it is the Brauer group of $S$, and may be 
identified with the subgroup of the Brauer group of $K$ unramified 
outside $p$. Since there is only one prime above $p$ and $K$ is 
totally complex, this subgroup is 
zero. Thus $H^{1}_{\text{\'et}}(S, \mu_{p}) \iso C/pC \iso A/pA$. 

(2) 
From class field theory we have  an  exact sequence
$$
0 \to U^{(1)}/\overline{E^{(1)}} \to G^{\text{ab}} \to A \to 0,
$$
where $U$ is the group of units in the completion of $K$ at the prime 
above $p$, and $\overline{E}$ is the closure of the global units in 
$U$, and the superscript 1 indicates that we take units congruent to 1 
modulo $p$. 
Let $(\alpha) = \a^{p}$ where $\a$ is an ideal of $K$ representing a 
 generator $a$ for $A(K)$.  Write  
$\alpha = \alpha'u$ for some $u
\in \mu_{p-1}(K_{\p})$, $\alpha' \in U^{(1)}$. 
If $a$ lifts to $\overline{a} \in G^{\text{ab}}$, 
then $\overline{a}^{p} \in U^{(1)}/\overline{E^{(1)}}$ is represented by 
$\alpha'$.
Since $\alpha 
\neq 1$ and $\mu_{p-1}(K) = \{1\}$,  $\alpha' \neq 1$.
It follows that 
the image of $\alpha_{\p}$ in  $U^{(1)}/\overline{E^{(1)}}$ is not a 
torsion element. Indeed, 
$\a \in A(i)$ where $i$ is 
odd and $i \neq 1$, hence $\overline{E^{(1)}}(i) = 0$.
Therefore, all the torsion in $G^{\text{ab}}$ is 
contained in $U^{(1)}/\overline{E^{(1)}}$.

(3)
Let $\tilde K$ denote the compositum of all $\ZZ_{p}$-extensions of 
$K$.
We will first show that $\tilde K$ is contained in the field obtained 
by 
adjoining $p$-power roots of $p$-units to $K(\mu_{p^\infty}) \subset 
K_{\infty}$. 
Indeed, let $$\omega:\ZZ_{p}^{\times} \to \mu_{p-1}(\ZZ_{p})$$ be the 
Teichmuller character (characterized by $\omega(a) \equiv a \pmod{p}$). 
Let $(\zeta_{n})$ be a generator of $\ZZ_{p}(1) = \invlim 
\mu_{p^{n}}$, and let 
$$
\eta_{k,n} = \prod_{a \in 
(\ZZ/p^n\ZZ)^{\times}}(1-\zeta_n^{a})^{\omega^{k}(a)a^{-1}}.
$$
(Here $a^{-1}$ denotes any integer congruent to $a^{-1}$ modulo $p^{n}$.)
Then $\eta_{k,n}$ is a $p$-unit in $\QQ(\mu_{p^{n}})$. Further
$G_{n} = \Gal(\QQ(\mu_{p^{n}})/K)$ acts on the class of $\eta_{k,n}$ 
in $\QQ(\mu_{p^n})^{\times}/\QQ(\mu_{p^n})^{\times p^n}$ via the cyclotomic 
character, and $\Delta = \Gal(K/\QQ)$ acts via the 
character $\omega^{-k}$. Hence, by Kummer theory, the $p^n$-th root of 
$\eta_{{k,n}}$ 
generates an extension $L_{n}/\QQ(\mu_{p^{n}})$ on which $G_{n}$ acts trivially, 
and which is 
therefore abelian over $K$. Thus $L_{n} = K(\mu_{p^{n}})K_{n}$ for 
some $p$-extension $K_{n}/K$ such that $\Delta$ acts on 
$\Gal(K_{n}/K)$  via 
$\omega^{k+1}$.  Furthermore, if $m > n$,
$$N_{\QQ(\mu_{p^{m}})/\QQ(\mu_{p^{n}})}\eta_{k,m} = \eta_{k, n}^{p^{m-n}},$$ so 
$K_{n}\subset K_{m}$, and  the union of $K_{n}$, $n>1$, forms an 
extension of $K$ which is either trivial or is a 
$\ZZ_{p}$-extension on which $\Delta$ acts by $\omega^{k+1}$. In fact 
the extension is non-trivial precisely when $k=p-2$ (the cyclotomic 
$\ZZ_{p}$-extension) or when $0\le k \le p-1$ and $k$ is even. The 
compositum of all these extensions has galois group
$\ZZ_{p}^{(p+1)/2}$, the maximal allowable rank by 
class field theory.  Hence it is $\tilde K$, 
and is visibly  contained in $N_\infty$.

By definition, $\Gal(\tilde K/K)$ is the torsion-free quotient of 
$G^{\text{ab}}$, so it follows from (2) that $\Gal(K^{\text{ab}}/\tilde K)
\iso \ZZ/p\ZZ$.
On the other hand, we can generate a $p$-ramified 
extension of $\tilde K$
whose galois group is isomorphic to $\ZZ/p\ZZ$ simply by taking the 
extension $K(\alpha^{1/p})$, where $\alpha$ is as in 
(2). This 
extension is  a non-trivial $p$-ramified extension (it is not 
contained in $\tilde K$ because of eigenspace considerations) and is 
abelian over $K$ since $\alpha \in K$. 
Thus $K^{\text{ab}} = \tilde K (\alpha^{1/p})$. We claim 
that $\a$ capitulates in $K_{\infty}$.
This implies that
 $\alpha$ can be written as a $p$-th power times a unit in 
 $K_{\infty}$, and thus $\alpha^{1/p} \in N_{\infty}$. Since we have 
 already seen that $\tilde K \subset N_{\infty}$, this proves the 
 lemma.

 To prove the claim, we show that
 $K_{\infty}$ contains $H$, the $p$-Hilbert class field 
of $K$. It follows from well-known results in the theory of cyclotomic 
fields that
$H$ is generated 
by the $p^{\hbox{\rm th}}$-th root of a 
cyclotomic unit. (The reflection principle, \cite{W},  Theorem 10.9, 
implies that $p$ does not divide $h^+$. The claim now follows from 
Theorems 10.16, 8.2, and 8.25.) Since the $p$-th roots of all the 
cyclotomic units in $K$
are contained in $K_{\infty}$, this implies 
that $H \subset K_{\infty}$.   \endproof

\section{The Main Theorem}

\theorem\label{ReallyMainTheorem}
Let $K_{\infty}/K$ be a multiple 
$\ZZ_{p}$-extension of $K$, $\Gal(K_{\infty}/K) \iso \ZZ_{p}^{r}$, $r 
\ge 2$.
Suppose that $K$ and $K_{\infty}$ satisfy the following conditions:
\list
   \item\label{containsmu} $K_\infty$ contains all p-power roots of unity
     \item\label{onlyoneprime} there is only one prime of $K$ above $p$
    \item\label{satisfiesLeopoldt} $K$ satisfies Leopoldt's conjecture
     \item\label{sisone} $\dim_{\ZZ/p\ZZ}H^{2}({\Cal G},\ZZ/p\ZZ) = 1$
     \item\label{bonk} the maximal abelian $p$-ramified pro-$p$ extension of $K$ is 
    contained in $$N_{\infty} = (K_\infty(\epsilon^{1/p^{n}}: n \in {\Bbb N}, \epsilon 
\in {\Cal O}_{K_{\infty}}[1/p]^\times)$$
\endlist
Then $A$ is a pseudo-null $\Lambda$-module.
\endtheorem

\proof
It 
follows from Theorem~\Ref{torsionfreetheorem} that
$Y_\tor$ fixes $N_{\infty}$. Since $K^{\text{ab}} \subset N_{\infty}$, 
it follows from Lemma~\ref{criterion} that $Y_{\tor}=0$. 
The theorem follows from Corollary~\ref{greenbergcorollary} (note that
since $K$ has only one prime above $p$, $r_\p = r \ge 2$).
\endproof

\demo{\it Proof of Theorem \Ref{MainTheorem}}
Let $K = \QQ(\zeta_p)$ and suppose that $K$ satisfies conditions (1) 
and (2) in the statement of Theorem~\Ref{MainTheorem}. 
Let $K_{\infty}$ be the compositum of all 
$\ZZ_{p}$-extensions of  $K$. Trivially,
conditions~\Ref{containsmu} and \Ref{onlyoneprime} of 
Theorem~\Ref{ReallyMainTheorem} are satisfied. 
Condition~\Ref{satisfiesLeopoldt},  Leopoldt's conjecture, is 
well-known
(\cite{W}, Theorem 5.25).
Finally, by 
Lemma~\ref{bla}, conditions~\Ref{sisone} and~\Ref{bonk} of 
Theorem~\ref{ReallyMainTheorem} are satisfied. 
 \enddemo
 
\makebib

\widestkey{BCEMS}

\bib \key Br \by K. S. Brown \book Cohomology of Groups \publ
Springer-Verlag \publaddr New York \yr 1982 \endbib   

\bib\key {B-H} \by W. Bruns \& J. Herzog \book Cohen-Macaulay rings,
Cambridge studies in advanced mathematics {\bf 39} \publ Cambridge
University Press \publaddr Cambridge \yr 1993 \endbib

\bib\key BCEMS \by J. Buhler, R. Crandall, R. Ernvall, T. 
Mets\"ankyl\"a, M. A. Shokrollahi \paper Irregular Primes and 
Cyclotomic Invariants to 12 Million \jour J. Symbolic Computation 
\vol 31 \yr 2001 \pp 89--96 \endbib

\bib\key G1 \by R. Greenberg \paper 
On the Iwasawa Invariants of Totally Real Number Fields
\jour Amer. J. Math. \vol  98, No. 1 \yr 1976 \pp  263--284 \endbib

\bib\key G2 \bysame \paper Iwasawa Theory---Past and Present 
\jour 
\lit* http://www.math.washington.edu"linebreak"0/~greenber* \endbib

\bib\key H \by J. A. Hillman \book Alexander Ideals of Links, Lecture 
Notes in Mathematics {\bf 895} \publ Springer-Verlag \publaddr Berlin 
\yr 1981 \endbib

\bib \key Iw \by K. Iwasawa \paper On ${\ZZ}_\ell$-extensions
of algebraic number fields \jour Annals of Math. \vol 98
\yr 1973 \pp 246--326 \endbib

\bib \key IS \by K. Iwasawa and C. C. Sims \paper Computation of 
invariants in the theory of cyclotomic fields \jour J. Math. Soc. 
Japan \vol 18, No. 1 \yr 1966 \pp 86--96 \endbib

\bib \key J1 \by U. Jannsen \paper On the structure of Galois groups
as Galois modules \inbook in Number Theory Noordwijkerhout 1983,
Lecture Notes in Math. {\bf 1068} \publ Springer \publaddr
Berlin \yr 1984 \pp 109--126
\endbib

\bib \key J2 \bysame \paper Iwasawa Modules up to 
Isomorphism \inbook in Advanced Studies in Pure Mathematics {\bf 17}
\publ Academic Press \publaddr Orlando \yr 1989 \pp 171-207
\endbib

\bib \key LN \paper Conjectures de Greenberg et extensions 
pro-$p$-libres d'un corps de nombres \by A. Lannuzel \& T. 
Nguyen Quang Do \jour to appear in Manuscripta Mathematica \endbib

\bib \key Mc \by W. G. McCallum \paper Duality Theorems in the 
Multivariable Iwasawa Theory of Number Fields, 
\lit*http://www.math.arizona.edu/~wmc* \endbib

\bib \key Ma \by H. Matsumura \book Commutative Algebra
\publ Benjamin \publaddr Reading \yr 1980
\endbib

\bib \key Ng \by T. Nguyen-Quang-Do \paper
 Formations de classes et modules  d'Iwasawa
\inbook Number Theory Noordwijkerhout, 1983, Lecture Notes in Math.
{\bf 1068}
\publ Springer
\publaddr Berlin \yr 1984 \pp 167-185 
\endbib

\bib \key W \by L. C. Washington \book Introduction to Cyclotomic 
Fields \publ Springer-Verlag \publaddr New York \yr 1982 
\endbib

\endmakebib

\enddocument